\documentclass[12pt, leqno]{article}
\usepackage{amsmath,amsfonts,amsthm,amscd,amssymb,verbatim} 
\usepackage[all]{xy}

\numberwithin{equation}{section}

\theoremstyle{plain}
\newtheorem{theorem}{Theorem}[section]
\newtheorem{definition}[theorem]{Definition}
\newtheorem{proposition}[theorem]{Proposition}
\newtheorem{lemma}[theorem]{Lemma}
\newtheorem{corollary}[theorem]{Corollary}

\newtheorem{question}[theorem]{Question}

\theoremstyle{definition}
\newtheorem{remark}[theorem]{Remark}
\newtheorem{example}[theorem]{Example} 
\newtheorem{hypothesis}[theorem]{Hypothesis}

\renewcommand{\b}{\bullet}
\newcommand{\beast}{\begin{eqnarray*}}
\newcommand{\east}{\end{eqnarray*}}

\newcommand{\N}{{\Bbb N}}
\newcommand{\Z}{{\Bbb Z}}

\newcommand{\Q}{{\Bbb Q}}

\newcommand{\C}{{\Bbb C}}

\newcommand{\G}{{\Bbb G}}
\newcommand{\fG}{\wh{{\Bbb G}}}

\newcommand{\Spec}{{\mathrm{Spec}}\,}

\newcommand{\Spf}{{\mathrm{Spf}}\,}

\newcommand{\lra}{\longrightarrow}
\newcommand{\ra}{\rightarrow}
\newcommand{\hra}{\hookrightarrow}

\newcommand{\isom}{\overset{\sim}{=}}
\newcommand{\ti}[1]{\widetilde{#1}}

\newcommand{\ol}[1]{\overline{#1}}
\newcommand{\os}{\overset}

\newcommand{\crys}{{\mathrm{crys}}}

\newcommand{\Ker}{{\mathrm{Ker}}}
\newcommand{\im}{{\mathrm{Im}}}

\newcommand{\id}{{\mathrm{id}}}

\newcommand{\Inthom}{{\cal H}{\mathrm{om}}}

\newcommand{\cD}{{\cal D}}
\newcommand{\cE}{{\cal E}}
\newcommand{\cF}{{\cal F}}

\newcommand{\cH}{{\cal H}}

\newcommand{\cO}{{\cal O}}
\newcommand{\cP}{{\cal P}}

\newcommand{\wh}{\widehat}

\newcommand{\lara}[1]{\langle #1 \rangle}

\newcommand{\pa}{\partial}

\newcommand{\MIC}{{\mathrm{MIC}}}
\newcommand{\MIWC}{{\mathrm{MIWC}}}
\newcommand{\HIG}{{\mathrm{HIG}}}

\newcommand{\qn}{{\mathrm{qn}}}
\newcommand{\Hyp}{{\bold{Hyp}}}

\begin{document}
\title{Notes on generalizations of local Ogus-Vologodsky 
correspondence}
\author{Atsushi Shiho
\footnote{
Graduate School of Mathematical Sciences, 
University of Tokyo, 3-8-1 Komaba, Meguro-ku, Tokyo 153-8914, JAPAN. 
E-mail address: shiho@ms.u-tokyo.ac.jp \, 
Mathematics Subject Classification (2010): 12H25, 14F30, 14F40.}}
\date{}
\maketitle

\begin{abstract}
Given a smooth scheme over $\Z/p^n\Z$ with a lift of relative Frobenius 
to $\Z/p^{n+1}\Z$, 
we construct a functor from the category of Higgs modules to that 
of modules with integrable connections as the composite of 
the level raising inverse image functors from the category of 
modules with integrable $p^{m}$-connections to that of 
modules with integrable $p^{m-1}$-connections for $1 \leq m \leq n$. 
In the case $m=1$, we prove that the level raising inverse image 
functor is an equivalence  
when restricted to quasi-nilpotent objects, which 
generalizes a local result of Ogus-Vologodsky. 
We also prove that 
the above level raising inverse image functor for a smooth $p$-adic 
formal scheme induces an equivalence of $\Q$-linearized categories 
for general $m$ when restricted to nilpotent objects (in strong sense), 
under a strong condition on Frobenius lift. 
We also prove a similar result for 
 the category of modules with integrable $p^{m}$-Witt-connections. 
\end{abstract}

\tableofcontents

\section*{Introduction}
For a proper smooth algebraic variety $X$ over $\C$, the equivalence of 
the category of modules endowed with integrable 
connections on $X$ and the category of Higgs modules on $X$ 
(with semistability and vanishing Chern number condition) is 
established by Simpson \cite{simpson0}. 
In search of the analogue of it in chatacteristic 
$p>0$, Ogus and Vologodsky proved in \cite{ov} similar equivalences 
for a smooth scheme $X_1$ over a scheme $S_1$ of characteristic $p>0$. 
One of their results \cite[2.11]{ov} is described 
as follows: Let $F_{S_1}: S_1 \lra S_1$ be the 
Fronenius morphism, let us put 
$X_1^{(1)} := S_1 \times_{F_{S_1},S_1} X_1$, denote the 
projection $X_1^{(1)} \lra S_1$ by $f^{(1)}_1$ and let  
$F_{X_1/S_1}: X_1 \lra X^{(1)}_1$ be 
the relative Frobenius morphism. 
Assume that we are given 
smooth lifts $f_2: X_2 \lra S_2, f^{(1)}_2: X^{(1)}_2 \lra S_2$ 
of $f_1, f^{(1)}_1$ to morphisms 
of flat $\Z/p^2\Z$-schemes and 
a lift $F_2: X_2 \lra X^{(1)}_2$ of the morphism 
$F_{X_1/S_1}$ 
which is a morphism over $S_2$. Then there exists an equivalence 
between the category of quasi-nilpotent Higgs modules on 
$X^{(1)}_1$ and the category of modules with 
quasi-nilpotent integrable connections on $X_1$. 
(There is also a version \cite[2.8]{ov} which does not assume the existence of 
$f_2$ and $F_2$. In this case, the categories they treat are 
a little more restricted in some sense.) 
The proof is done by using Azumaya algebra property of the sheaf 
$\cD^{(0)}_{X_1/S_1}$ of differential operators of level $0$ on 
$X_1$ over $S_1$. Also, in the above situaion, one can give an 
explicit description of this functor as the inverse image by 
`divided Frobenius'. 
(See \cite[2.11.2]{ov}, \cite[5.9, 6.5]{glsq} or 
Remark \ref{qrem} in this paper.) \par 
The purpose of this paper is to construct a functor from Higgs modules 
to modules with integrable connections 
for smooth schemes over some 
flat $\Z/p^n\Z$-schemes and study the properties of this functor 
and related functors. 
Let us fix $n \in \N$ and 
let $S_{n+1}$ be a scheme flat over $\Z/p^{n+1}\Z$. 
Let us put $S_j := S_{n+1} \otimes \Z/p^j\Z \,(j \in \N, \leq n+1)$, 
let $f_1, F_{S_1}$ be as above and for $0 \leq m \leq n$, 
let us put 
$X_1^{(m)} := S_1 \times_{F^m_{S_1},S_1} X_1$, denote the 
projection $X_1^{(m)} \lra S_1$ by $f^{(m)}_1$ and for 
$1 \leq m \leq n$, let  
$F^{(m)}_{X_1/S_1}: X^{(m-1)}_1 \lra X^{(m)}_1$ be 
the relative Frobenius morphism for $f^{(m)}_1$. Moreover, 
assume that we are given 
a smooth lift $f_{n+1}: X_{n+1} \lra S_{n+1}$ of $f_1$, 
smooth lifts $f^{(m)}_{n+1}: X^{(m)}_{n+1} \lra S_{n+1}$ of $f_1^{(m)} \, 
(0 \leq m \leq n)$ with $f_1^{(0)} = f_1$ and 
lifts $F^{(m)}_{n+1}: X^{(m-1)}_{n+1} \lra X^{(m)}_{n+1}$ of the morphisms 
$F^{(m)}_{X_1/S_1} \,(1 \leq m \leq n)$ 
which is a morphism over $S_{n+1}$. Furthermore, 
let $f_{n}: X_{n} \lra S_{n}$, 
$f^{(m)}_{n}: X^{(m)}_{n} \lra S_{n} \, (0 \leq m \leq n)$, 
$F^{(m)}_{n}: X^{(m-1)}_{n} \lra X^{(m)}_{n}$ be 
$f_{n+1} \otimes \Z/p^n\Z, f^{(m)}_{n+1} \otimes \Z/p^n\Z, 
F^{(m)}_{n+1} \otimes \Z/p^n\Z$, respectively. 
Under this assumption, 
first we construct 
the functor $\Psi^{\qn}: \HIG(X^{(n)}_n)^{\qn} \lra \MIC(X_n)^{\qn}$ 
from the category $\HIG(X^{(n)}_n)^{\qn}$ 
of quasi-nilpotent Higgs modules on $X^{(n)}_n$ 
to the category $\MIC(X_n)^{\qn}$ of 
modules with quasi-nilpotent integrable connections on $X_n$. \par 
The construction of the functor $\Psi^{\qn}$ is done in the 
following way: 
For $0 \leq m \leq n$, let $\MIC^{(m)}(X^{(m)}_n)^{\qn}$ be the 
category of modules with quasi-nilpotent integrable $p^m$-connections 
(also called `quasi-nilpotent integrable connections of level $-m$') on 
$X^{(m)}_n$, that is, the category of pairs $(\cE,\nabla)$ 
consisting of an $\cO_{X^{(m)}_n}$-module $\cE$ endowed with 
an additive map $\nabla: \cE \lra \cE \otimes \Omega^1_{X^{(m)}_n/S_n}$ 
with the property $\nabla(fe) = f\nabla(e) + p^me \otimes df 
\,(e \in \cE, f \in \cO_{X^{(m)}_n})$, the integrability and 
the quasi-nilpotence 
condition. (The notion of $p^m$-connection
is the case $\lambda=p^m$ of the notion of $\lambda$-connection of 
Simpson \cite{simpson}.) 
Note that we have 
$\HIG(X^{(n)}_n)^{\qn} = \MIC^{(n)}(X^{(n)}_n)^{\qn}$ since 
$p^n=0$ on $X^{(n)}_n$. We construct the functor 
$\Psi^{\qn}$ as the composite of the functors 
$F^{(m),*,\qn}_{n+1}: \MIC^{(m)}(X^{(m)}_n)^{\qn} \lra 
\MIC^{(m-1)}(X^{(m-1)}_n)^{\qn}$ for $1 \leq m \leq n$, which 
are defined as the inverse image by divided Frobenius (we call 
it the level raising inverse image) associated to $F^{(m)}_{n+1}$, 
as in the case of \cite[2.11.2]{ov}, \cite[5.9]{glsq}. \par 
The first naive question might be whether the functor 
$\Psi^{\qn}$ is an equivalence or not. 
(Note that, 
since the sheaf of differential operators of level $0$ on 
$X_n$ over $S_n$ does not seem to have Azumaya algebra property, 
it would be hard to generarize the method of Ogus-Vologodsky 
in this case.) 
Unfortunately, the functors $F^{(m),*,\qn}$ are not equivalences 
(not full, not essentially surjective) for $m \geq 2$ and so 
the functor $\Psi^{\qn}$ is not an equivalence either. 
So an interesting question would be to construct nice functors from 
the functors $F^{(m),*,\qn}_{n+1}$. 
As a first answer to this question, 
we prove that the functor 
$F^{(1),*,\qn}_{n+1}$ is an equivalence, under the assumption that 
there exists a closed immersion $S_{n+1} \hra S$ into a $p$-adic 
formal scheme flat over $\Z_p$. This generalizes \cite[2.11]{ov} 
under the existence of $S$. 
To prove this, we may work locally and so we may assume the 
existence of 
a smooth lift $f: X \lra S$ of $f_{n+1}$, a smooth lift 
$f^{(1)}: X^{(1)} \lra S$ of $f_1^{(1)}$ and 
a lift $F^{(1)}: X \lra X^{(1)}$ of the morphism 
$F^{(1)}_{X_1/S_1}$. In this situation, we introduce the notion 
of the sheaf of $p$-adic differential operators $\cD^{(-1)}_{X^{(1)}/S}$ of 
level $-1$ of $X^{(1)}$ over $S$, which is a level $-1$ version of the 
sheaf of $p$-adic differential operators defined and studied by Berthelot 
\cite{berthelotI}, \cite{berthelotII}, and prove that 
the category 
$\MIC^{(1)}(X^{(1)}_n)^{\qn}$ is 
equivalent to the category of $p^n$-torsion quasi-nilpotent left 
$\cD_{X^{(1)}/S}^{(-1)}$-modules 
and re-define the functor $F^{(1),*,\qn}$ as the level raising 
inverse image functor from the category of $p^n$-torsion 
quasi-nilpotent left 
$\cD_{X^{(1)}/S}^{(-1)}$-modules to 
 the category of $p^n$-torsion quasi-nilpotent left 
$\cD_{X/S}^{(0)}$-modules. This is 
a negative level version of the level raising 
inverse image functor defined in \cite[2.2]{berthelotII}. 
(We also give a definition of the 
sheaf of $p$-adic differential operators of level $-m$ and 
give the interpretation of (quasi-nilpotent) integrable $p^m$-connections 
and (level raising) inverse image functors of them 
in terms of $\cD$-modules for general $m \in \N$.) 
Then, one can prove the equivalence of the functor $F^{(1),*,\qn}$
by following the proof of 
Frobenius descent by Berthelot \cite[2.3]{berthelotII}. 
\par 
To explain our next result, let us fix $m \in \N$ and 
assume that we are given 
smooth lifts $f: X \lra S$ of $f_1$, 
$f^{(1)}: X^{(1)} \lra S$ of $f_1^{(1)}$ and a `nice' 
lift $F: X \lra X^{(1)}$ of 
$F_{X_1/S_1}$ which is a morphism over $S$. 
Then, under certain assumption on nilpotence condition, 
we prove that the $i$-th de Rham cohomology of an object 
in the category $\MIC^{(m)}(X^{(1)})$ of integrable $p^m$-connections 
on $X^{(1)}$ is isomorphic to 
the $i$-th de Rham cohomology of its level raising inverse image by $F$ 
(which is an object in the category $\MIC^{(m-1)}(X)$ of 
integrable $p^{m-1}$-connections on $X$) when $i=0$, and 
isomorphic modulo torsion for general $i$. 
This implies that the level raising inverse image
 by $F$ induces a fully faithful 
functor between a full subcategory of 
$\MIC^{(m)}(X^{(1)})$ and that of $\MIC^{(m-1)}(X)$ 
satisfying certain nilpotence condition and that it induces 
an equivalence of $\Q$-linearized categories (again under 
nilpotence condition), which gives a 
second answer to the question we raised in the previous 
paragraph. \par 
Also, we prove a Witt version of the result in the previous 
paragraph: We introduce the notion of 
integrable $p^m$-Witt connection $\MIWC^{(m)}(X_1)$ on $X_1$ 
for a smooth morphism $X_1 \lra S_1$ of chatacteristic $p>0$ 
with $S_1$ perfect and 
the level raising inverse image functor 
$F_*: \MIWC^{(m)}(X_1^{(1)}) \lra \MIWC^{(m-1)}(X_1)$ 
for $X_1 \lra S_1$ as above and $X_1^{(1)} := 
S_1 \times_{F_{S_1},S_1} X_1$. We prove that the functor 
$F_*$ induces a fully faithful functor 
between a full subcategory of 
$\MIWC^{(m)}(X^{(1)}_1)$ and that of $\MIWC^{(m-1)}(X_1)$ 
satisfying certain nilpotence condition and that it induces 
an equivalence of $\Q$-linearized categories (again under 
nilpotence condition). This result has an advantage that 
we need no assumption on the liftability of objects and 
Frobenius morphisms. \par 
The content of each section is as follows: 
In Section 1, we introduce the notion of integrable 
$p^m$-connections, (level raising) inverse image functors 
and define the functor $\Psi$ ($\Psi^{\qn}$) 
from the category of (quasi-nilpotent) Higgs modules 
to the category of (quasi-nilpotent) integrable connections 
as the iteration of 
level raising inverse image functors. We also give an example 
which shows that the functor $\Psi, \Psi^{\qn}$ are not 
equivalences. \par 
In Section 2, 
we introduce the sheaf of $p$-adic differential operators of 
negative level and prove basic properties of it. 
In particular, 
we prove the equivalence between the category of (quasi-nilpotent) left 
$\cD$-modules of level $-m$ and the category of (quasi-nilpotent) 
modules with integrable $p^m$-connections ($m \in \N$). 
In the case of schemes over $\Z/p^n\Z$, we also introduce certain 
crystalized categories to describe the categories of 
(quasi-nilpotent) modules with integrable $p^m$-connections, and 
prove certain crystalline property of them. 
In Section 3, we prove that 
the level raising inverse image functors from the category of 
modules with integrable $p$-connections to that of 
modules with integrable connections defined in Section 1 is an 
equivalence of categories 
when restricted to quasi-nilpotent objects. 
In Section 4, 
we compare the de Rham cohomology of certain 
modules with integrable $p^m$-connections over smooth 
$p$-adic formal schemes and the de Rham cohomology of 
the pull-back of them by the level raising inverse image functor, 
and deduce a full-faithfulness 
(resp. an equivalence) between 
the category (resp. the $\Q$-linearized category) of 
modules with integrable $p^m$-connections to that of 
modules with integrable $p^{m-1}$-connections satisfying 
nilpotent conditions. 
In Section 4, 
we introduce the notion of modules with integrable 
$p^m$-Witt connections and the level raising inverse image 
functors from the category of modules with integrable 
$p^m$-Witt connections to that of modules with integrable 
$p^{m-1}$-Witt connections. 
We compare the de Rham cohomology of certain 
modules with integrable $p^m$-Witt-connections and 
the de Rham cohomology of the pull-back of them by the level raising 
inverse image functor, 
and deduce a full-faithfulness 
(resp. an equivalence) between 
the category (resp. the $\Q$-linearized category) of 
modules with integrable $p^m$-Witt-connections to that of 
modules with integrable $p^{m-1}$-Witt-connections satisfying 
nilpotent conditions. \par 
The author is partly supported by Grant-in-Aid for Young Scientists (B) 
21740003 (representative: Atsushi Shiho) from the Ministry of Education, 
Culture, Sports, Science and Technology, Japan 
and 
Grant-in-Aid for Scientific Research (B) 22340001 
(representative: Nobuo Tsuzuki) from 
Japan Society for the Promotion of Science. 

\section*{Convention}
Throughout this paper, $p$ is a fixed prime number. 
Fiber products and tensor products in this paper 
are $p$-adically completed ones, unless otherwise stated. 

\section{Modules with integrable $p^m$-connections}

In this section, we define the notion of modules with 
integrable $p^m$-connections. 
Also, we define the inverse image 
functor of the categories of modules with integrable 
$p^m$-connections, and the `level raising inverse image functor' 
from the categories of modules with integrable 
$p^m$-connections to that of integrable $p^{m-1}$-connections for a lift of 
Frobenius morphism. As a composite of the level raising 
inverse image functors, we define the functor from the category 
of Higgs modules to the category of modules with integrable 
connections, which is a generalization of the inverse of 
local Cartier transform of Ogus-Vologodsky \cite[2.11]{ov}. \par 
First we give a definition of $p^m$-connection, which is a special 
case (the case $\lambda = p^m$) of the notion of $\lambda$-connection 
of Simpson \cite{simpson}. 

\begin{definition}
Let $X \lra S$ be a smooth morphism of schemes over $\Z/p^n\Z$ or 
$p$-adic formal schemes and let 
$m \in \N$. A $p^m$-connection on an $\cO_X$-module $\cE$ is an 
additive map $\nabla: \cE \lra \cE \otimes_{\cO_X} \Omega^1_{X/S}$ 
satisfying $\nabla(fe) = f\nabla(e) + p^me \otimes df$ for 
$e \in \cE, f \in \cO_X$. We call a $p^m$-connection also as 
a connection of level $-m$. 
\end{definition}

When we are given an $\cO_X$-module with $p^m$-connection 
$(\cE,\nabla)$, we can define the additive map 
$\nabla_k: \cE \otimes_{\cO_X} \Omega^k_{X/S} \lra 
\cE \otimes_{\cO_X} \Omega^{k+1}_{X/S}$ which is characterized by 
$\nabla_k(e \otimes \omega) = \nabla(e) \wedge \omega + 
p^me \otimes \omega$. 

\begin{definition}
With the notation above, we call $(\cE,\nabla)$ integrable 
if we have $\nabla_1 \circ \nabla = 0$. We denote the category of 
$\cO_X$-modules with integrable $p^m$-connection by 
$\MIC^{(m)}(X)$. \par 
When $m=0$, the notion of modules
with integrable $p^m$-connection is nothing but that of modules 
with integrable connection. In this case, we denote the category 
$\MIC^{(m)}(X)$ also by $\MIC(X)$. Also, 
when $X, S$ are schemes over $\Z/p^n\Z$ with $n \leq m$, the notion of 
modules with integrable $p^m$-connection is nothing but that of 
Higgs modules. In this case, we denote the category 
$\MIC^{(m)}(X)$ also by $\HIG(X)$. 
\end{definition}

\begin{remark}\label{torsrem1}
For a smooth morphism $f: X \lra S$ of $p$-adic formal schemes 
and $n \in \N$, we denote the full subcategory of $\MIC^{(m)}(X/S)$ 
consisting of $p^n$-torsion objects by $\MIC^{(m)}(X/S)_n$. 
If we denote the morphism $f \otimes \Z/p^n\Z$ by 
$X_n \lra S_n$, the direct image by the canonical closed immersion 
$X_n \hra X$ induces the equivalence of categories 
$\MIC^{(m)}(X_n) \os{=}{\lra} \MIC^{(m)}(X)_n$. 
\end{remark}

Let us assume given a commutative diagram 
\begin{equation}\label{commm}
\begin{CD}
X @>g>> Y \\ 
@VVV @VVV \\ 
S @>>> T 
\end{CD}
\end{equation}
of smooth morphism of schemes over $\Z/p^n\Z$ or 
$p$-adic formal schemes with smooth vertical arrows and 
and an object $(\cE,\nabla)$ in $\MIC^{(m)}(Y)$ 
(where $\MIC^{(m)}(Y)$ is defined for the morphism $Y \lra T$). 
Then we can endow a structure of an integrable $p^m$-connection 
$g^*\nabla$ on $g^*\cE$ by 
$g^*\nabla(fg^*(e)) = fg^*(\nabla(e)) + p^mg^*(e) \otimes df 
\,(e \in \cE, f \in \cO_X)$. 
So we have the inverse image functor 
$$g^*: \MIC^{(m)}(Y) \lra \MIC^{(m)}(X); \,\,\,\, 
(\cE,\nabla) \mapsto g^*(\cE,\nabla) := (g^*\cE,g^*\nabla). $$ 

\begin{remark}\label{torsrem2}
Let us assume given the commutative diagram \eqref{commm} of 
$p$-adic formal schemes with smooth vertical arrows and 
let us denote the morphism $g \otimes \Z/p^n\Z$ by 
$g_n: X_n \lra Y_n$. Then the inverse image functor $g^*$ above 
induces the functor $g^*: \MIC^{(m)}(Y)_n \lra \MIC^{(m)}(X)_n$, 
and it coincides with the inverse image functor 
$g_n^*: \MIC^{(m)}(Y_n) \lra \MIC^{(m)}(X_n)$ associated to $g_n$ 
via the equivalences $\MIC^{(m)}(X_n) \os{=}{\lra} 
\MIC^{(m)}(X)_n, \MIC^{(m)}(Y_n) \os{=}{\lra} 
\MIC^{(m)}(Y)_n$ of Remark \ref{torsrem1}. 
\end{remark}

Next we introduce the notion of quasi-nilpotence. 
Let $X \lra S$ be a smooth morphism of schemes over $\Z/p^n\Z$ 
which admits a local coordinate $t_1,...,t_d$. Then, for 
$(\cE,\nabla) \in \MIC^{(m)}(X)$, we can write $\nabla$ as 
$\nabla(e) = \sum_{i=1}^d \theta_i(e) dt_i$ for some additive maps 
$\theta_i: \cE \lra \cE \,(1 \leq i \leq d)$. Then we have 
$$ 0 = \nabla_1 \circ \nabla(e) 
= \sum_{i<j} (\theta_i\theta_j-\theta_j\theta_i)(e) dt_i \wedge dt_j. $$
So we have $\theta_i \theta_j = \theta_j \theta_i$. Therefore, 
for $a = (a_1,...,a_d) \in \N^d$, the map 
$\theta^a := \prod_{i=1}^d \theta_i^{a_i}$ is well-defined. 

\begin{definition}\label{qn0}
With the above situation, we 
call $(\cE,\nabla)$ quasi-nilpotent with respect to 
$(t_1,...,t_d)$ if, for any section $e \in \cE$, there exists some 
$N \in \N$ such that $\theta^a(e)=0$ for any $a \in \N^d$ with 
$|a| \geq N$. 
\end{definition}

\begin{lemma}\label{qnlem}
The above definition of quasi-nilpotence does not depend on 
the local coordinate $(t_1,...,t_d)$. 
\end{lemma}

\begin{proof}
When $m=0$, this is classical. Here we prove the lemma in the case 
$m>0$. (The proof is easier in this case.) First, let 
us note that, for $f \in \cO_X$, we have the equality 
\begin{align*}
\sum_i \theta_i(fe) dt_i & = \nabla(fe) = f\nabla(e) + p^me df 
= 
\sum_i(f\theta_i(e) + p^m\dfrac{\pa f}{\pa t_i}) dt_i. 
\end{align*}
So we have the equality 
\begin{equation}\label{tra0}
\theta_if = f\theta_i + p^m \dfrac{\pa f}{\pa t_i}. 
\end{equation}
Now let us take another local coordinate $t'_1,...,t'_d$, and write 
$\nabla$ as 
$\nabla(e) = \sum_{i=1}^d \theta'_i(e) dt'_i$. Then we have 
$$ 
\sum_i \theta_i(e) dt_i = 
\sum_{i,j}\theta_i(e) \dfrac{\pa t_i}{\pa t'_j} dt'_j 
= \sum_j (\sum_i \dfrac{\pa t_i}{\pa t'_j} \theta_i(e)) dt'_j. $$
Hence we have $\theta'_j = \sum_i \dfrac{\pa t_i}{\pa t'_j} \theta_i$. \par 
Let us prove that, for any $e \in \cE$ and $a \in \N^d$, there exists some 
$f_{a,b} \in \cO_X \,(b \in \N^d, |b| \leq |a|)$ with 
\begin{equation}\label{tra}
{\theta'}^a(e) = \sum_{|b| \leq |a|} p^{m(|a|-|b|)} f_{a,b} \theta^b(e),  
\end{equation}
by induction on $a$: Indeed, this is trivially true when $a=0$. 
If this is true for $a$, we have 
\begin{align*}
\theta'_j{\theta'}^a(e) & = 
(\sum_i \dfrac{\pa t_i}{\pa t'_j} \theta_i)
(\sum_{|b| \leq |a|} p^{m(|a|-|b|)} f_{a,b} \theta^b)(e) \\ 
& = 
\sum_{i,b} (p^{m(|a|-|b|)}\dfrac{\pa t_i}{\pa t'_j} 
f_{a,b} \theta^{b+e_i}(e) + p^{m(|a|-|b|+1)} \dfrac{\pa f_{a,b}}{\pa t_i} 
\theta^b(e)) \,\,\,\, (\text{by \eqref{tra0}})
\end{align*}
and from this equation, we see that the claim is true for $a+e_j$. \par 
Now let us assume that $(\cE,\nabla)$ is quasi-nilpotent with respect to 
$(t_1,...,t_d)$, and take a local section $e \in \cE$. 
Then there exists some $N \in \N$ such that $\theta^b(e)=0$ for 
any $b \in \N^d, |b| \geq N$. Then, for any $a \in \N^d, |a| \geq N+n$, 
we have either $|b| \geq N$ or $|a| - |b| \geq n$ for any 
$b \in \N^d$. Hence we have either $p^{m(|a|-|b|)}=0$ or 
$\theta^b(e)=0$ on the right hand side of \eqref{tra} and so we have 
${\theta'}^a(e)=0$. So we have shown that $(\cE, \nabla)$ is 
quasi-nilpotent with respect to $(t'_1,...,t'_d)$ and so we are done. 
\end{proof}

\begin{remark}\label{genqn}
By \eqref{tra0}, we have 
$$ \theta^a(fe) = \sum_{0 \leq b \leq a} p^{m|b|}\dfrac{\pa^b f}{\pa t^b}
\theta^{a-b}(e)$$
for $e \in \cE, f \in \cO_X$, and we have 
$p^{m|b|}\dfrac{\pa^b f}{\pa t^b} \in p^{m|b|} b!\cO_X$. 
Hence, if we have $\theta^a(e) = 0$ for any $a \in \N^d, |a| \geq N$, 
we have $\theta^a(fe)=0$ for any $a \in \N^d, |a| \geq N+p^nd$. Therefore, 
to check the quasi-nilpotence of $(\cE,\nabla)$ (with respect to some local 
coordinate $t_1,...,t_d$), it suffices to check that, for some 
generator $e_1, ..., e_r$ of $\cE$, there exists some $N \in \N$ 
such that $\theta^a(e_i)=0$ for any $a \in \N^d, |a| \geq N$ and $1 \leq 
i \leq r$. 
\end{remark}

When a given morphism does not admit a local coordinate globally, 
we define the notion of quasi-nilpotence as follows: 

\begin{definition}\label{qn1}
\begin{enumerate}
\item 
Let $X \lra S$ be a smooth morphism of schemes over $\Z/p^n\Z$. 
Then an object $(\cE,\nabla)$ in $\MIC^{(m)}(X)$ is called 
quasi-nilpotent if, locally on $X$, there exists a local coordinate 
$t_1, ..., t_d$ of $X$ over $S$ such that $(\cE,\nabla)$ is 
quasi-nilpotent with respect to $(t_1,...,t_d)$. 
$($Note that, by Lemma \ref{qnlem}, this definition is independent of 
the choice of $t_1, ..., t_d.)$
\item 
Let $X \lra S$ be a smooth morphism of $p$-adic formal 
schemes. 
Then an object $(\cE,\nabla)$ in $\MIC^{(m)}(X)$ is called 
quasi-nilpotent if it is contained in $\MIC^{(m)}(X)_n$ for some 
$n$ and the object in $\MIC^{(m)}(X_n)$ $($where $X_n := X \otimes 
\Z/p^n\Z)$ corresponding to $(\cE,\nabla)$ via the equivalence 
in Remark \ref{torsrem1} is quasi-nilpotent. 
\end{enumerate}
We denote the full subcategory of $\MIC^{(m)}(X)$ consisting of 
quasi-nilpotent objects by $\MIC^{(m)}(X)^{\qn}$, and in the case of 
$(2)$, we denote the category $\MIC^{(m)}(X)_n \cap \MIC^{(m)}(X)^{\qn}$ 
by $\MIC^{(m)}(X)^{\qn}_n$. 
\end{definition}

Next we prove the functoriality of quasi-nilpotence. 

\begin{proposition}\label{functqn}
Let us assume given a commutative diagram \eqref{commm}
of smooth morphism of $p$-adic formal schemes or 
schemes over $\Z/p^n\Z$ 
with smooth vertical arrows. Then the inverse image functor 
$g^*: \MIC^{(m)}(Y) \lra \MIC^{(m)}(X)$ induces the 
functor $g^{*,\qn}: \MIC^{(m)}(Y)^{\qn} \lra \MIC^{(m)}(X)^{\qn}$, 
that is, $g^*$ sends quasi-nilpotent objects to quasi-nilpotent 
objects. 
\end{proposition}

\begin{proof}
In view of Remark \ref{torsrem2}, it suffices to prove the case 
of schemes over $\Z/p^n\Z$. 
When $m=0$, the proposition is classical (\cite{bo}, \cite{berthelotII}). 
So we may assume $m>0$. 
Since the quasi-nilpotence is a local property, we may assume that 
there exists a local coordinate $t_1, ..., t_d$ (resp. $t'_1, ..., t'_{d'}$) 
of $X$ over $S$ (resp. $Y$ over $T$). 
Let us take 
an object $(\cE,\nabla)$ in $\MIC^{(m)}(Y)^{\qn}$ and write the map 
$\nabla, g^*\nabla$ as 
$\nabla(e) = \sum_j \theta'_j(e) dt'_j$, 
$g^*\nabla(fg^*(e)) = \sum_i \theta_i(fg^*(e)) dt_i$. 
Let us write $g^*(dt'_j) = \sum_i a_{ij}dt_i$. 
Then we have 
\begin{align*}
g^*\nabla(fg^*(e)) & = fg^*(\nabla(e)) + p^mg^*(e) \otimes df \\ 
& = 
\sum_{i,j} (a_{ij} f g^*(\theta'_j(e)) + p^m g^*(e) \dfrac{\pa f}{\pa t_i}) 
dt_i 
\end{align*}
and so we have $\theta_i(fg^*(e)) = 
\sum_j a_{ij} f g^*(\theta'_j(e)) + p^m \dfrac{\pa f}{\pa t_i} g^*(e)$. 
Let us prove that, for any $e \in \cE, f \in \cO_X$ and 
$a \in \N^d$, there exists some 
$f_{a,b} \in \cO_X \,(b \in \N^d, |b| \leq |a|)$ (which depends on $e, f$) 
with 
\begin{equation}\label{tra2}
{\theta}^a(fg^*(e)) = \sum_{|b| \leq |a|} p^{m(|a|-|b|)} f_{a,b} 
g^*({\theta'}^b(e)),  
\end{equation}
by induction on $a$: Indeed, this is trivially true when $a=0$. 
If this is true for $a$, we have 
\begin{align*}
\theta_i{\theta}^a(fg^*(e)) & = 
\theta_i (\sum_{|b| \leq |a|} p^{m(|a|-|b|)} f_{a,b} 
g^*({\theta'}^b(e))) \\ 
& = 
\sum_{j,b} (a_{ij}  p^{m(|a|-|b|)} f_{a,b} g^*({\theta'}^{b+e_j}(e)) + 
p^{m(|a|-|b|+1)} \dfrac{\pa f_{a,b}}{\pa t_i} g^*({\theta'}^b(e)))
\end{align*}
and from this equation, we see that the claim is true for $a+e_i$. 
From \eqref{tra2}, we can prove the quasi-nilpotence of 
$(g^*\cE,g^*\nabla)$ as in the proof of Lemma \ref{qnlem}. 
So we are done. 
\end{proof}

Before we define the level raising inverse image functor, 
we give the following definition to fix the situation. 

\begin{definition}\label{hypabc}
In this definition, `a scheme flat over $\Z/p^{\infty}\Z$' means 
a $p$-adic formal scheme flat over $\Z_p$. \par 
For $a, b, c \in \N \cup\{\infty\}$ with $a \geq b \geq c$, 
we mean by $\Hyp(a,b,c)$ the following hypothesis$:$ 
We are given a scheme $S_a$ flat over $\Z/p^a\Z$, and for 
$j \in \N, j \leq a$, $S_j$ denotes the scheme $S_a \otimes \Z/p^j\Z$. 
We are also given a smooth morphism of finite type 
$f_1: X_1 \lra S_1$, and let 
$F_{X_1}: X_1 \lra X_1, F_{S_1}: S_1 \lra S_1$ be the 
Frobenius endomorphism ($p$-th power endomorphism). Let us put 
$X_1^{(1)} := S_1 \times_{F_{S_1},S_1} X_1$ and denote the 
projection $X_1^{(1)} \lra S_1$ by $f^{(1)}_1$. 
Then the morphism 
$F_{X_1}$ induces the relative Frobenius morphism 
$F_{X_1/S_1}: X_1 \lra X^{(1)}_1$. 
Assume that we are given 
a smooth lift $f_b: X_b \lra S_b$ of $f_1$, 
a smooth lift $f^{(1)}_b: X^{(1)}_b \lra S_b$ of $f_1^{(1)}$, 
and for $j \in \N, j \leq b$, 
denote the morphism $f_b \otimes \Z/p^j\Z$, $f_b^{(1)} \otimes 
\Z/p^j\Z$ by $f_j: X_j \lra S_j$, 
$f^{(1)}_j: X^{(1)}_j \lra S_j$, respectively. Also, assume we are 
given a lift $F_c: X_c \lra X^{(1)}_c$ of the morphism $F_{X_1/S_1}$ 
which is a morphism over $S_c$. For $j \in \N, j \leq c$, 
let $F_{j}: X_{j} \lra X^{(1)}_{j}$ be 
$F_{n+1} \otimes \Z/p^j\Z$. Finally, 
when $a = \infty$ $($resp. $b=\infty$, $c=\infty)$, 
we denote $S_a$ $($resp. $f_b: X_b \lra S_b$ and 
$f^{(1)}_b: X^{(1)}_b \lra S_b$, 
$F_c:  X_c \lra X^{(1)}_c)$ simply by 
$S$ $($resp. $f: X \lra S$ and 
$f^{(1)}: X^{(1)} \lra S$, 
$F:  X \lra X^{(1)})$.
\end{definition}

Roughly speaking, $\Hyp(a,b,c)$ means that $S_1$ is liftable to 
a scheme flat over $\Z/p^a\Z$, $f_1:X_1 \lra S_1$ and $f^{(1)}_1: 
X^{(1)}_1 \lra S_1$ is liftable to morphisms over $S_b$ and 
the relative Frobenius $F_{X_1/S_1}: X_1 \lra X^{(1)}_1$ is liftable 
to a morphism over $S_c$. 

Now we define the level raising inverse image functor for 
a lift of Frobenius. 
Let $n \in \N$ and assume that we are in the situation of $\Hyp(n+1,n+1,n+1)$. 
When we work locally, we can take a local coordinate 
$t_1, ..., t_d$ of $X_{n+1}$ over $S_{n+1}$. If 
us put $t'_i := 1 \otimes t_i \in \cO_{X^{(1)}_{n+1}}\,(1 \leq i \leq d)$, 
they form a local coodinate of $X^{(1)}_{n+1}$ over $S_{n+1}$, and 
we have $F_{n+1}^*(t'_i) = t_i^p + pa_i$ for some $a_i \in \cO_{X_{n}}$. 
Hence we have $F_{n+1}^*(dt'_i) = p(t_i^{p-1}dt_i + da_i)$, that is, 
the image of the homomorphism $F_{n+1}^*: \Omega^1_{X_{n+1}^{(1)}/S_{n+1}} 
\lra \Omega^1_{X_{n+1}/S_{n+1}}$ is contained in 
$p\Omega^1_{X_{n+1}/S_{n+1}}$. So there exists a unique morphism 
$\ol{F}_{n+1}^*: \Omega^1_{X^{(1)}_n/S_n} \lra \Omega^1_{X_n/S_n}$ 
which makes the following diagram commutative: 
\begin{equation*}
\begin{CD}
\Omega^1_{X^{(1)}_{n+1}/S_{n+1}} @>{F_{n+1}^*}>> 
\Omega^1_{X_{n+1}/S_{n+1}} \\ 
@VVV @A{p}AA \\ 
\Omega^1_{X^{(1)}_n/S_n} @>{\ol{F}_{n+1}^*}>> \Omega^1_{X_n/S_n}. 
\end{CD}
\end{equation*}
Using this, we define the level raising inverse image functor 
\begin{equation}\label{lrii}
F_{n+1}^*: \MIC^{(m)}(X_n^{(1)}) \lra \MIC^{(m-1)}(X_n)
\end{equation}
as follows: an object $(\cE,\nabla)$ in $\MIC^{(m)}(X_n^{(1)})$ 
is sent by $F_{n+1}^*$ the object $(F_n^*\cE, F_n^*\nabla)$, where 
$F_n^*\nabla$ is defined by 
$F_n^*\nabla(fF_n^*(e)) = f\ol{F}_{n+1}^*(\nabla(e)) + p^{m-1}
F_n^*(e)df.$ (One can check the integrability by direct calculation, using 
local coordinate.) \par 
Also, in the situation of $\Hyp(\infty,\infty,\infty)$, we have the 
homomorphism 
$\ol{F}^*: \Omega^1_{X^{(1)}/S} \lra \Omega^1_{X/S}$ 
which makes the diagram 
\begin{equation*}
\begin{CD}
\Omega^1_{X^{(1)}/S} @>{F^*}>> 
\Omega^1_{X/S} \\ 
@VVV @A{p}AA \\ 
\Omega^1_{X^{(1)}/S} @>{\ol{F}^*}>> \Omega^1_{X/S} 
\end{CD}
\end{equation*}
commutative, and using this, we can define the level raising inverse 
image functor 
\begin{equation}\label{lriif}
F^*: \MIC^{(m)}(X^{(1)}) \lra \MIC^{(m-1)}(X)
\end{equation}
 in the same way. 

\begin{remark}\label{torsrem3}
Assume we are in the situation of $\Hyp(\infty,\infty,\infty)$ and put 
$X_n^{(1)} := X^{(1)} \otimes \Z/p^n\Z, 
X_n := X \otimes \Z/p^n\Z$. Then the functor \eqref{lriif}
 induces the functor 
$F^*: \MIC^{(m)}(X^{(1)})_n \lra \MIC^{(m-1)}(X)_n$, and this coincides 
with the functor \eqref{lrii} via the equivalences 
$\MIC^{(m)}(X^{(1)}_n) \os{=}{\lra} \MIC^{(m)}(X^{(1)})_n, 
\MIC^{(m-1)}(X_n) \os{=}{\lra} \MIC^{(m-1)}(X)_n$ of Remark \ref{torsrem1}. 
\end{remark}

\begin{remark}\label{qrem}
Assume that we are in the situation of $\Hyp(2,2,2)$. 
Then the level raising inverse image functor for $m=1$ is written as 
$F_2^*: \HIG(X^{(1)}_1) \lra \MIC(X_1)$. 
Let us see how it is calculated locally. 
Let us take local coordinate $t_1,...,t_d$ of $X_2$, let us 
put $t'_i := 1 \otimes t_i \in \cO_{X^{(1)}_2}, 
F_2^*(t'_i) = t_i^p + pa_i$. Then 
$\ol{F}_2^*: \Omega^1_{X_1^{(1)}/S_1} \lra \Omega^1_{X_1/S_1}$ 
is written as $\ol{F}_2^*(dt'_i) = t_i^{p-1} dt_i + da_i$ and 
the functor $F_2^*$ is defined by using it. 
So we obtain the following expression of the functor $F^{*}_2$: 
A Higgs module $(\cE,\theta)$ on $X^{(1)}_1$ of the form 
$\theta(e) = \sum_{i=1}^d \theta_i(e) \otimes dt'_i$ is sent to 
the integrable connection $(F_{X_1/S_1}^*\cE, \nabla)$ 
such that, if we write $\nabla = \sum_{i=1}^d \pa_i dt_i$, 
we have 
$$ \pa_i(1 \otimes e) = t_i^{p-1}\otimes \theta_i(e) + 
\sum_{j=1}^d \dfrac{\pa a_j}{\pa t_i} \otimes \theta_j(e). $$ 
Let $\iota: \HIG(X^{(1)}_1) \lra \HIG(X^{(1)}_1)$ 
be the functor $(\cE,\theta) \mapsto (\cE,-\theta)$. 
Then, by the above expression, 
we see that 
the functor $F_2^* \circ \iota$ coincides with a special case of the 
functor defined in \cite[5.8]{glsq} (the case $m=0$ in the notation 
of \cite{glsq}) for quasi-nilpotent objects. 
(The underlying 
sheaf $F^{*}_1\cE$ is globally the same as 
the image of the functor in \cite[5.8]{glsq}, and the connections 
coincide because they coincide locally.) 
Hence, by \cite[6.5]{glsq}, 
it coincides with the functor in \cite[2.11]{ov} for quasi-nilpotent 
objects. 
\end{remark}

We have the functoriality of quasi-nilpotence with respect to level raising 
inverse image functors, as follows: 

\begin{proposition}\label{functqn2}
Assume that we are in the situation of $\Hyp(n+1,n+1,n+1) \,(n \in \N)$ 
$($resp. $\Hyp(\infty, \infty, \infty)).$ 
Then the level raising inverse image functor 
$F_{n+1}^*: \MIC^{(m)}(X^{(1)}_n) \lra \MIC^{(m-1)}(X_n)$ 
$($resp. $F^*: \MIC^{(m)}(X^{(1)}) \lra \MIC^{(m-1)}(X))$
induces the 
functor $F_{n+1}^{*,\qn}: \MIC^{(m)}(X^{(1)}_n)^{\qn} \lra 
\MIC^{(m-1)}(X_n)^{\qn}$ 
$($resp. $F^{*,\qn}: \MIC^{(m)}(X^{(1)})^{\qn} \lra \MIC^{(m-1)}(X)^{\qn})$, 
that is, $F_{n+1}^*$ $($resp. $F^*)$ 
sends quasi-nilpotent objects to quasi-nilpotent 
objects. 
\end{proposition}

\begin{proof}
In view of Remark \ref{torsrem3}, it suffices to prove the 
proposition for $F_{n+1}^*$. 
In the case $m=n=1$, the functor $F_2^* \circ \iota$ 
($\iota$ is as in Remark \ref{qrem}) coincides with the 
functor in \cite[2.11]{ov}. Hence it sends 
quasi-nilpotent objects to quasi-nilpotent 
objects. Since $\iota$ induces an auto-equivalence of 
$\MIC^{(1)}(X^{(1)}_1)^{\qn}$, we see that $F_2^*$ sends 
quasi-nilpotent objects to quasi-nilpotent 
objects. \par 
Next, let us prove the proposition 
in the case $m=1$ and $n$ general, by induction on $n$. 
Let us take an object $(\cE,\nabla)$ in $\MIC^{(1)}(X_n^{(1)})$. 
Then we have the exact sequence 
$$ 0 \lra (p\cE, \nabla|_{p\cE}) \lra 
(\cE,\nabla) \lra (\cE/p\cE,\ol{\nabla}) \lra 0, $$
where $\ol{\nabla}$ is the $p$-connection on $\cE/p\cE$ induced by 
$\nabla$. Since $F_n: X^{(1)}_n \lra X_n$ is finite flat, 
the above exact sequence induces the following exact sequence: 
$$ 0 \lra F_{n+1}^*(p\cE, \nabla|_{p\cE}) \lra 
F_{n+1}^*(\cE,\nabla) \lra F_{n+1}^*(\cE/p\cE,\ol{\nabla}) \lra 0. $$
Then, since $F_{n+1}^*(p\cE, \nabla|_{p\cE}) = F_{n}^*(p\cE,\nabla|_{p\cE})$ 
and $F_{n+1}^*(\cE/p\cE,\ol{\nabla}) = F_2^*(\cE/p\cE,\ol{\nabla})$, 
they are quasi-nilpotent by induction hypothesis. Hence 
$F_{n+1}^*(\cE,\nabla)$ is also quasi-nilpotent, as desired. \par 
Finally we prove the proposition in the case $m \geq 2$. 
Let us take local coordinates $t_1, ..., t_d$ of $X_n^{(1)}$ and 
put $t'_i := 1 \otimes t_i$. Take an object $(\cE,\nabla)$ in 
$\MIC^{(m)}(X^{(1)}_n)$ and write 
$\nabla(e) := \sum_i \theta'(e) dt'_i$, 
$F_n^*\nabla(fF_n^*(e)) = \sum_i \theta(fF_n^*(e)) dt_i$. 
Then we can prove that, 
for any $e \in \cE, f \in \cO_{X_n}$ and 
$a \in \N^d$, there exists some 
$f_{a,b} \in \cO_{X_n} \,(b \in \N^d, |b| \leq |a|)$ (which depends on $e, f$) 
with 
\begin{equation}\label{tra3}
{\theta}^a(fg^*(e)) = \sum_{|b| \leq |a|} p^{(m-1)(|a|-|b|)} f_{a,b} 
g^*({\theta'}^b(e)), 
\end{equation}
in the same way as the proof of Proposition \ref{functqn}. From this 
we see the quasi-nilpotence of $F_{n+1}^*(\cE,\nabla) = 
(F_n^*\cE, F_n^*\nabla)$ again 
in the same way as the proof of Proposition \ref{functqn}. 
\end{proof}

\begin{remark}\label{qnatode}
In the above proof, we used the results in \cite{ov}. 
Later, we give another proof of Proposition \ref{functqn2} 
which does not use any results in \cite{ov} under 
$\Hyp(\infty, n+1, n+1)$ or $\Hyp(\infty, \infty, \infty)$. 
\end{remark}

Now we define a functor from the category of (quasi-nilpotent) 
Higgs modules to the category of modules with (quasi-nilpotent) 
integrable connections as a composite of level raising inverse 
image functors. Let us consider the following hypothesis. 

\begin{hypothesis}\label{hypo-1}
Let us fix $n \in \N$ and 
let $S_{n+1}$ be a scheme flat over $\Z/p^{n+1}\Z$. 
For $j \in \N, j \leq n+1$, let us put $S_j := S_{n+1} 
\otimes \Z/p^j\Z$. 
Let $f_1: X_1 \lra S_1$ be a smooth morphism of finite type and let 
$F_{S_1}: S_1 \lra S_1$ be the 
Frobenius endomorphism. For $0 \leq m \leq n$, 
let us put 
$X_1^{(m)} := S_1 \times_{F^m_{S_1},S_1} X_1$, denote the 
projection $X_1^{(m)} \lra S_1$ by $f^{(m)}_1$ and for 
$1 \leq m \leq n$, let  
$F^{(m)}_{X_1/S_1}: X^{(m-1)}_1 \lra X^{(m)}_1$ be 
the relative Frobenius morphism for $f^{(m-1)}_1$. \par 
Assume that we are given 
a smooth lift $f_{n+1}: X_{n+1} \lra S_{n+1}$ of $f_1$, 
smooth lifts $f^{(m)}_{n+1}: X^{(m)}_{n+1} \lra S_{n+1}$ of 
$f_1^{(m)} \, (0 \leq m \leq n)$ with $f_{n+1}^{(0)} = f_{n+1}$ and 
lifts $F^{(m)}_{n+1}: X^{(m-1)}_{n+1} \lra X^{(m)}_{n+1}$ of the morphism 
$F^{(m)}_{X_1/S_1} \,(1 \leq m \leq n)$ 
which is a morphism over $S_{n+1}$. Finally, 
let $f_{n}: X_{n} \lra S_{n}$, 
$f^{(m)}_{n}: X^{(m)}_{n} \lra S_{n}$, 
$F^{(m)}_{n}: X^{(m-1)}_{n} \lra X^{(m)}_{n}$ be 
$f_{n+1} \otimes \Z/p^n\Z, f^{(m)}_{n+1} \otimes \Z/p^n\Z, 
F^{(m)}_{n+1} \otimes \Z/p^n\Z$, respectively. 
\end{hypothesis}

Then we define the functor as follows: 

\begin{definition}\label{defdef}
Assume that we are in the situation of Hypothesis \ref{hypo-1}. 
Then we define the canonical functors 
$$ \Psi: \HIG(X^{(n)}_n) \lra \MIC(X_n), \quad 
\Psi^{\qn}: \HIG(X^{(n)}_n)^{\qn} \lra \MIC(X_n)^{\qn} $$
as the composite 
$F^{(n),*}_{n+1} \circ F^{(n-1),*}_{n+1} \circ \cdots \circ 
F^{(1),*}_{n+1}, 
F^{(n),*,\qn}_{n+1} \circ F^{(n-1),*,\qn}_{n+1} 
\circ \cdots \circ F^{(1),*,\qn}_{n+1}$ 
of level raising inverse image functors 
\begin{align*}
& F^{(m),*}_{n+1}: \MIC^{(m)}(X^{(m)}_n) 
\lra \MIC^{(k)}(X^{(m-1)}_n) \quad (1 \leq m \leq n), \\ 
& F^{(m),*,\qn}_{n+1}: \MIC^{(m)}(X^{(m)}_n)^{\qn} 
\lra \MIC^{(m-1)}(X^{(m-1)}_n)^{\qn} \quad (1 \leq m \leq n),  
\end{align*}
respectively. 
\end{definition}

Since the morphisms $F^{(m)}_n$ are finite flat, we see that 
the functors $\Psi, \Psi^{\qn}$ are exact and faithful. 
However, we see in the following example that the functors 
$\Psi, \Psi^{\qn}$ are not so good as one might expect. 

\begin{example}
In this example, let us put $S_{n+1} = \Spec \Z/p^{n+1}\Z$ and let 
$S_j := S_{n+1} \otimes \Z/p^j\Z = \Spec \Z/p^j\Z, 
X_j := \Spec (\Z/p^j\Z)[t^{\pm 1}]$ for 
$j \in \N, j \leq n+1$.  
Also, let 
$X^{(m)}_j := 
\Spec (\Z/p^j\Z)[t^{\pm 1}]$ for all $m \in \N, j \in \N, j \leq n+1$ and 
let $F^{(m)}_j: X^{(m-1)}_j \lra X^{(m)}_j$ be the morphism defined by 
$t \mapsto t^p$. 
Then $\ol{F}_{n+1}^{(m),*}: 
\Omega^1_{X^{(m)}_n/S_n} \lra \Omega^1_{X^{(m-1)}_n/S_n}$ 
sends $f(t)t^{-1}dt$ to $f(t^p)t^{-p}t^{p-1}dt = f(t^p)t^{-1}dt$ 
and the level raising inverse image functor 
$F^{(m),*}_n: \MIC^{(m)}(X^{(m)}_n) \lra \MIC^{(m-1)}(X^{(m-1)}_n)$ 
is defined as `the pull-back by $\ol{F}_{n+1}^{(m),*}$'. \par 
For $m \in \N$ and $f(t) \in \Gamma(X^{(m)}_n,\cO_{X^{(m)}_n})$, we 
define the $p^m$-connection $(\cO_{X^{(m)}_n},\nabla_{f(t)})$ by 
$\nabla_{f(t)} = p^md + f(t)t^{-1}dt$. It is locally free of rank $1$, 
and since any locally free sheaf of rank $1$ on $X^{(m)}_n$ is free, 
any $p^m$-connection on $X^{(m)}_n$ which is locally free of rank $1$ 
has the form $(\cO_{X^{(m)}_n},\nabla_{f(t)})$ for some $f(t)$. 
For a $p^{m}$-connection $(\cO_{X^{(m)}_n},\nabla_{f(t)})$, 
the $p^{m-1}$-connection $F^{(m),*}_{n+1}(\cO_{X^{(m)}_n},\nabla_{f(t)})$ 
is equal to $(\cO_{X^{(m-1)}_n}, \nabla_{f(t^p)})$ thanks to the 
description of $F^{(m),*}_{n+1}$ given in the previous paragraph. \par 
Let us make some more calculation on the $p^m$-connection 
$(\cO_{X^{(m)}_n},\nabla_{f(t)})$. It is easy to see that 
we have an isomorphism 
$(\cO_{X^{(m)}_n},\nabla_{f(t)}) \cong 
(\cO_{X^{(m)}_n},\nabla_0)$ if and only if 
$(\cO_{X^{(m)}_n},\nabla_{f(t)})$ is generated as $\cO_{X^{(m)}_n}$-module 
by a horizontal element. Since 
we have 
$$ \nabla_{f(t)}(g(t)) = p^k \dfrac{dg}{dt}dt + gft^{-1}dt 
= g(p^mtg^{-1}\dfrac{dg}{dt}+f)t^{-1}dt, $$
we see that $(\cO_{X^{(m)}},\nabla_{f(t)})$ is isomorphic to 
$(\cO_{X^{(m)}_n},\nabla_0)$ if and only if there exists an element 
$g \in \Gamma(X^{(m)}_n, \cO_{X^{(m)}_n}^{\times})$ with 
$f = -p^mtg^{-1}\dfrac{dg}{dt}$. If $g$ is an element in 
$\Gamma(X^{(m)}, \cO_{X^{(m)}_n}^{\times})$, it has the form 
$g=t^N + ph_1$ for some $N \in \Z$ and $h_1 \in 
\Gamma(X^{(m)}, \cO_{X^{(m)}})$, and in this case, $g^{-1}$ 
has the form $t^{-N} + ph_2$ for some element $h_2$ in 
$\Gamma(X^{(m)}_n, \cO_{X^{(m)}_n})$. Then we have 
\begin{align*}
-p^mtg^{-1}\dfrac{dg}{dt} & = 
-p^mt(t^{-N}+ph_2)(Nt^{N-1}+p\dfrac{dh_1}{dt}) \\ 
& = -p^mN + p^{m+1}h(t) 
\end{align*}
for some  $h \in \Gamma(X^{(m)}_n, \cO_{X^{(m)}_n})$. Therefore, 
we have shown that if $(\cO_{X^{(m)}_n},\nabla_{f(t)})$ is isomorphic to 
$(\cO_{X^{(m)}_n},\nabla_0)$, $f$ has the form $-p^mN + p^{m+1}h(t)$. \par 
Now, to investigate the functor 
$F^{(1),*}_{n+1}: \MIC^{(1)}(X^{(1)}_n) \lra \MIC(X_n)$, first let us consider
 the $p$-connection $(\cO_{X^{(1)}_n},\nabla_{1})$. 
Then, since there does not exist $N \in \N, h \in 
\Gamma(X^{(1)}_n,\cO_{X^{(1)}_n})$ with $1 = -pN+p^2h$, 
it is not isomorphic to $(\cO_{X^{(1)}_n}, \nabla_0)$. On the other hand, 
we see that the connection $F^{(1),*}_{n+1}(\cO_{X^{(1)}_n}, \nabla_1) = 
(\cO_{X_n},\nabla_1)$ is isomorphic to 
$F^{(1),*}_{n+1}(\cO_{X^{(1)}_n},\nabla_0) = 
(\cO_{X_n},\nabla_0)$ because we have $1 = -tg^{-1}\dfrac{dg}{dt}$ 
when $g=t^{-1}$. So we can conclude that the functor $F^{(1),*}$ 
is not full. Secondly, let us consider the connection 
$(\cO_{X_n},\nabla_t)$. If it is contained in the essential image 
of $F^{(1),*}_{n+1}$, we should have $(\cO_{X_n},\nabla_t) \cong 
F^{(1),*}_{n+1}(\cO_{X^{(1)}_n},\nabla_{f(t)}) = (\cO_{X_n}, \nabla_{f(t^p)})$ 
for some $f(t)$. Then we have 
$f(t^p) - t = -N + ph$ for some $N \in \Z$ and $h \in 
\Gamma(X_n,\cO_{X_n})$, 
but it is impossible. Hence we see that the functor $F^{(1),*}_{n+1}$ 
is not essentially surjective. \par 
Next, let us investigate the functors 
$F^{(m),*}_{n+1}: \MIC^{(m)}(X^{(m)}_n) \lra \MIC^{(m-1)}(X^{(m-1)}_n)$, 
$F^{(m),*,\qn}_{n+1}: \MIC^{(m)}(X^{(m)}_n)^{\qn} \lra 
\MIC^{(m-1)}(X^{(m-1)}_n)^{\qn}$ for $m \geq 2$. First let us consider the 
$p^m$-connection $(\cO_{X^{(m)}_n}, \nabla_{p^{m-1}})$. 
We see as in the previous paragraph that it is not isomorphic to 
$(\cO_{X^{(m)}_n}, \nabla_{0})$, and $F^{(m),*}_{n+1}
(\cO_{X^{(m)}_n}, \nabla_{p^{m-1}}) = (\cO_{X^{(m-1)}_n}, \nabla_{p^{m-1}})$ 
is isomorphic to $F^{(m),*}_{n+1}(\cO_{X^{(m)}_n}, \nabla_{0}) = 
(\cO_{X^{(m-1)}_n}, \nabla_{0})$. If we put $\nabla_{p^{m-1}}(e) = 
\pa(e) dt$, we can see easily by induction that 
$\pa^{l}(1) = (\prod_{i=0}^{l-1}(p^{m-1} - ip^{m}))/t^l$. 
By this and Remark \ref{genqn}, we see that 
$(\cO_{X^{(m)}_n}, \nabla_{p^{m-1}})$ is quasi-nilpotent. 
So we see that the functors $F^{(m),*}_{n+1}$, $F^{(m),*,\qn}_{n+1}$ are not 
full. Secondly, let us consider the connection 
$(\cO_{X^{(m-1)}_n},\nabla_{pt})$. If it is contained in the essential image 
of $F^{(m),*}_{n+1}$, we should have $(\cO_{X^{(m-1)}_n},\nabla_t) \cong 
F^{(m),*}_{n+1}(\cO_{X^{(m)}_n},\nabla_{f(t)}) 
= (\cO_{X^{(m-1)}_n}, \nabla_{f(t^p)})$ 
for some $f(t)$. Then we have 
$f(t^p) - pt = -p^{m-1}N + p^{m}h$ for some $N \in \N$ and $h \in 
\Gamma(X^{(m-1)}_n,\cO_{X^{(m-1)}_n})$, 
but it is impossible. Also, if we put $\nabla_{pt}(e) = 
\pa(e) dt$, we can see easily by induction that 
$\pa^{l}(1) = p^l$. 
By this and Remark \ref{genqn}, we see that 
$(\cO_{X^{(m-1)}_n}, \nabla_{pt})$ is quasi-nilpotent. Hence we see that 
the functors  $F^{(m),*}_{n+1}$, $F^{(m),*,\qn}_{n+1}$ are not 
essentially surjective. \par 
In conclusion, $\Psi$ is not full, not essentially 
surjective for any $n \in \N$, 
and $\Psi^{\qn}$ is not full, not essentially 
surjective for any $n \geq 2$. 
\end{example}

In view of the above example, we would like to ask the following question. 

\begin{question}\label{q}
Is it possible to construct some nice functor $($a fully faithful functor 
or an equivalence$)$ from the functors $F^{(m),*}_{n+1}, 
F^{(m),*,\qn}_{n+1}$ 
possibly under some more assumption? 
\end{question}

Several answers to this question will be given in 
Sections 3, 4 and 5. 

\section{$p$-adic differential operators of negative level}
 
In this section, first 
we introduce the sheaf of $p$-adic differential 
operators of level $-m \,(m \in \N)$, which is a `negative level version' 
of the sheaf of $p$-adic differential 
operators of level $m$ defined by Berthelot, for a smooth morphism 
of $p$-adic formal schemes flat over $\Z_p$. We prove that the 
equivalence of the notion of left $\cD$-modules in this sense 
and the notion of modules with integrable $p^m$-connections. 
We also define the inverse image functors and the level raising 
inverse image functors for left $\cD$-modules, which are compatible 
with the corresponding notion for modules with integrable $p^m$-connections 
over $p$-adic formal schemes. \par 
The definition of the sheaf of $p$-adic differential 
operators of level $-m \,(m \in \N)$ is possible only for 
smooth morphisms of $p$-adic formal schemes, because we use the 
formal blow-up with respect to an ideal containing $p^m$ in the 
definition. In the case of smooth morphisms $X_n \lra S_n$ of schemes flat 
over $\Z/p^n\Z$, we give a similar description by 
considering all the local lifts of $X_n$ to smooth 
$p$-adic formal scheme and consider the `crystalized' category of 
$\cD$-modules. \par 
We also consider a variant of the `crystalized' category of 
$\cD$-modules, which is also related to  
the category of integrable $p^m$-connections. As a consequence, 
we prove certain crystalline property for the category of integrable 
$p^m$-connections: 
When $f_n: X_n \lra S_n$ is a smooth morphism 
of flat $\Z/p^n\Z$-schemes and if we denote $f_n \otimes \Z/p\Z$ by 
$X_1 \lra S_1$, we know that the category $\MIC(X_n)^{\qn}$, which is 
equivalent to the category of crystals on the crystalline site 
$(X_1/S_n)_{\crys}$, depends only on the diagram $X_1 \lra S_1 \hra S_n$. 
We prove here similar results for the categories of integrable 
$p^m$-connections, although the result in the case $m>0$ is 
weaker than that in the case $m=0$. \par 

\subsection{The case of $p$-adic formal schemes}

Let $S$ be a $p$-adic formal scheme flat over $\Spf \Z_p$ and 
let $X$ be a $p$-adic formal scheme smooth of finite type 
over $S$. For a positive integer $r$, we denote the $r$-fold fiber 
product of $X$ over $S$ by $X^r$. For positive integers $m,r$, 
let $\ti{T}_{X,(-m)}(r)$ be the formal blow-up 
of $X^{r+1}$ along the ideal 
$p^m\cO_{X^{r+1}} + \Ker {\Delta(r)^*}$, 
where $\Delta(r)^*: \cO_{X^{r+1}} \lra \cO_X$ denotes the homomorphism 
induced by the diagonal map $\Delta(r): X \hra X^{r+1}$. 
Let $T_{X,(-m)}(r)$ be the open formal subscheme of 
$\ti{T}_{X,(-m)}(r)$ defined by 
\begin{align*}
T_{X,(-m)}(r) & := \{x \in \ti{T}_{X,(-m)}(r)\,|\, \\ 
& \hspace{1.5cm} p^m\cO_{\ti{T}_{X,(-m)}(r),x} = 
((p^m\cO_{X^{r+1}} + \Ker{\Delta(r)^*})
\cO_{\ti{T}_{X,(-m)}(r)})_x\}. 
\end{align*}
Then, since we have $(p^m\cO_{X^{r+1}}+\Ker\Delta(r)^*)|_X = p^m\cO_X$, 
the diagonal map $\Delta(r)$ factors through a morphism 
$\ti{\Delta}(r): 
X \hra T_{X,(-m)}(r)$ by the universality of formal blow-up. 
Let us put $I_{X,(-m)}(r) := \Ker \ti{\Delta}(r)^*$. Let 
$(\cP_{X,(-m)}(r), \ol{I}_{X,(-m)}(r))$ 
be the PD-envelope of $\cO_{T_{X,(-m)}(r)}$ with respect to the 
ideal $I_{X,(-m)}(r)$, and let us put $P_{X,(-m)}(r) := \Spf \cP_{X,(-m)}(r)$. 
Also, for $k \in \N$, let $\cP^k_{X,(-m)}(r), P^k_{X,(-m)}(r)$ be 
$\cP_{X,(-m)}(r), P_{X,(-m)}(r)$ modulo $\ol{I}_{X,(-m)}(r)^{[k+1]}$. 
In the case $r=1$, we drop the symbol $(r)$ from the notation. \par 
Note that $\cP_{X,(-m)}$ admits two $\cO_X$-algebra 
($\cO_X$-module) structure 
induced by the $0$-th and $1$-st projection $X^2 \lra X$, which we call 
the left $\cO_X$-algebra ($\cO_X$-module) structure and the right 
$\cO_X$-algebra ($\cO_X$-module) structure, respectively. Note also that, 
for $m' \leq m$, we have the canonical morphism 
$\cP_{X,(-m')} \lra \cP_{X,(-m)}$. \par 
Locally, $\cP_{X,(-m)}(r)$ is described in the following way. 
Assume that $X$ admits a local parameter $t_1,...,t_d$ over $S$. 
Then, if we denote the $q$-th projection $X^{r+1} \lra X$ by $\pi_i$ 
$\,(0 \leq q \leq r)$ and if we put $\tau_{i,q} := \pi_{q+1}^*t_i - 
\pi_q^*t_i$, $\Ker \Delta(r)^*$ 
is generated by $\tau_{i,q}$'s $\,(1 \leq i \leq d, 
0 \leq q \leq r-1)$ and we have 
$T_{X,(-m)}(r) = \Spf \cO_X\left\{\tau_{i,q}/p^m\right\}_{i,q}$, 
where $\{ - \}$ means the $p$-adically completed 
polynomial algebra. 
So we have 
$\cP_{X,(-m)}(r) =  \cO_X\left\langle\tau_{i,q}/p^m\right\rangle_{i,q}$, 
where $\langle - \rangle$ means the $p$-adically completed 
PD-polynomial algebra. \par 
We see easily that 
the identity map $X^{r+r'+1} \lra X^{r+1} \times_X X^{r'+1}$ 
naturally induces the isomorphism 
$\cP_{X,(-m)}(r) \otimes_{\cO_X} \cP_{X,(-m)}(r') \lra \cP_{X,(-m)}(r+r')$ and 
in the local situation, the element 
$\tau_{i,q}/p^m \otimes 1$ (resp. $1 \otimes \tau_{i,q}/p^m$) on the right 
hand side corresponds to the element $\tau_{i,q}/p^m$ (resp. 
$\tau_{i,q+r}/p^m$) on the left hand side. Then, the projection 
$X^3 \lra X^2$ to the $(0,2)$-th factor induces the homomorphism 
$$\delta: 
\cP_{X,(-m)} \lra \cP_{X,(-m)}(2) \cong \cP_{X,(-m)} \otimes_{\cO_X} 
\cP_{X,(-m)} $$ 
with $\delta(\tau_i/p^m) = \tau_i/p^m \otimes 1 + 1 \otimes \tau_i/p^m$ 
(here we denoted $\tau_{i,0}$ simply by $\tau_i$) and so 
it induces the homomorphism 
$$\delta^{k,k'}: 
\cP^{k+k'}_{X,(-m)} \lra \cP^k_{X,(-m)} \otimes_{\cO_X} 
\cP^{k'}_{X,(-m)}. $$ 
Using these, we define the sheaf of $p$-adic differential operators 
of negative level as follows: 

\begin{definition}
Let $X,S$ be as above. Then we define the sheaf 
$\cD_{X/S,k}^{(-m)}$ of $p$-adic differential operators of level $-m$ 
and order $\leq k$ by 
$\cD_{X/S,k}^{(-m)} := \Inthom_{\cO_X}(\cP^k_{X,(-m)},\cO_X)$ 
and the sheaf 
$\cD_{X/S}^{(-m)}$ of $p$-adic differential operators of level $-m$ by 
$\cD_{X/S}^{(-m)} := \allowbreak \bigcup_{k=0}^{\infty} \allowbreak 
\cD_{X/S,k}^{(-m)}$. 
We define the product 
$$  \cD_{X/S,k}^{(-m)} \times  \cD_{X/S,k'}^{(-m)} \lra 
\cD_{X/S,k+k'}^{(-m)} $$
by sending $(P,P')$ to the homomorphism 
$$ 
\cP_{X,(-m)}^{k+k'} \os{\delta^{k,k'}}{\lra}  
\cP^k_{X,(-m)} \otimes_{\cO_X} 
\cP^{k'}_{X,(-m)} \os{\id \otimes P'}{\lra}
\cP^k_{X,(-m)} \os{P}{\lra} O_X. $$ 
\end{definition}

By definition, $\cD_{X/S}^{(-m)}$ also admits two 
$\cO_X$-module structures, 
which are defined as the multiplication by the elements in 
$\cD_{X/S,0}^{(-m)} = \cO_X$ from left and from right. 
We call these the left and the right 
$\cO_X$-module strucrure of $\cD_{X/S}^{(-m)}$. 
Note that 
$P \in \cD_{X/S,k}^{(-m)}$ acts on $\cO_X$ as the composite 
$$ \cO_X \lra \cP^k_{X,(-m)} \os{P}{\lra} \cO_X $$ 
(where the first map is defined by $f \mapsto 1 \otimes f$), 
and this defines the action of $\cD_{X/S}^{(-m)}$ on $\cO_X$. 
For $m' \leq m$, the canonical map 
$\cP_{X,(-m')} \lra \cP_{X,(-m)}$ induces the homomorphism of rings 
$\rho_{-m',-m}: \cD_{X/S}^{(-m)} \lra \cD_{X/S}^{(-m')}$. \par 
Assume that $X$ admits a local parameter $t_1,...,t_d$ over $S$ and put 
$\tau_i := 1 \otimes t_i - t_i \otimes 1 \in \cO_{X^2}$. Then, as 
we saw before, we have 
$\cP_{X,(-m)} =  \cO_X\left\langle\tau_{i}/p^m\right\rangle_{i}$ and so 
$\cP_{X,(-m)}^k$ admits a basis 
$\{(\tau/p^m)^{[l]}\}_{|l|\leq k}$ as $\cO_X$-module. (Here and after, 
we use multi-index notation.) We denote the dual basis of 
it in $\cD_{X/S,k}^{(-m)}$ by $\{ \pa^{\lara{l}} \}_{|l|\leq k}$. 
When $l=(0,...,1,...,0)$ \,($1$ is placed in the $i$-th entry), 
$\pa^{\lara{l}}$ is denoted also by $\pa_i$. 
When we would like to clarify the level, we denote the element 
$\pa^{\lara{l}}$ by $\pa^{\lara{l}_{-m}}$. 
Since the canonical map $\cP_{X,(-m')} \lra \cP_{X,(-m)}$ sends 
$(\tau/p^{m'})^{[l]}$ to $p^{(m-m')|l|}(\tau/p^{m})^{[l]}$, 
we have $\rho_{-m',-m}(\pa^{\lara{l}_{-m}}) = 
p^{(m-m')|l|} \pa^{\lara{l}_{-m'}}$. \par 
We prove some formulas which are the analogues of the ones in 
\cite[2.2.4]{berthelotI}: 

\begin{proposition}\label{4}
With the above notation, we have the following$:$ 
\begin{enumerate}
\item 
For $f \in \cO_X$, $1 \otimes f = \sum_{|l|\leq k} \pa^{\lara{l}}(f) 
(\tau/p^m)^{[l]}$ in $\cP^k_{X,(-m)}$. 
\item 
$\pa^{\lara{l}}(t^i) = l! \begin{pmatrix} i \\ l \end{pmatrix} 
p^{m|l|} t^{i-l}$. 
\item 
$\pa^{\lara{l}}\pa^{\lara{l'}} = \pa^{\lara{l+l'}}$. 
\item 
$\pa^{\lara{k}}f = 
\sum_{k'+k''=k} \begin{pmatrix} k \\ k' \end{pmatrix} 
\pa^{\lara{k'}}(f) \pa^{\lara{k''}}$. 
\end{enumerate}
\end{proposition}

\begin{proof}
(1) is immediate from definition. By looking the coefficient of 
$(\tau/p^m)^{[l]}$ of $1 \otimes t^i = (t + p^m(\tau/p^m))^i$, 
we obtain (2). From the equality 
\begin{align*}
\pa^{\lara{l}}\pa^{\lara{l'}}((\tau/p^m)^{[i]}) & = 
\pa^{\lara{l}}(\id \otimes \pa^{\lara{l'}}) 
\delta^{|l|,|l'|}((\tau/p^m)^{[i]}) \\ 
& = 
\pa^{\lara{l}}(\id \otimes \pa^{\lara{l'}}) 
(\sum_{a+b=i}(\tau/p^m)^{[a]}\otimes (\tau/p^m)^{[b]}) 
= 
\left\{
\begin{aligned}
& 1, \quad \text{if $i=l+l'$,} \\ 
& 0, \quad \text{otherwise,} 
\end{aligned}
\right. 
\end{align*}
we see the assertion (3). From the equality 
\begin{align*}
(\pa^{\lara{k}}f)((\tau/p^m)^{[i]}) & = 
\pa^{\lara{k}}((1 \otimes f)(\tau/p^m)^{[i]}) \\ 
& = \pa^{\lara{k}}
(\sum_l \pa^{\lara{l}}(f)(\tau/p^m)^{[l]}(\tau/p^m)^{[i]}) \\ 
& = \pa^{\lara{k}}
(\sum_l \begin{pmatrix} l+i \\ l \end{pmatrix} 
\pa^{\lara{l}}(f)(\tau/p^m)^{[l+i]}),  
\end{align*}
we see the assertion (4).
\end{proof}

\begin{remark}
Let $\cD_{X/S}$ be the formal scheme version of the sheaf of 
usual differential operators and let us take a local basis 
$\{\pa^{[l]}\}_{l \in \N^d}$ of $\cD_{X/S}$, which can be 
defined in the same way as $\{\pa^{\lara{l}}\}$ above. 
Then $\cO_X$ admits the natural action of $\cD_{X/S}$ and we see, 
for $l \in \N^d$, $m \in \N$ and $f \in \cO_X$, the equalities 
$$ \pa^{\lara{l}_{-m}}(f) = p^{m|l|}\pa^{\lara{l}_0}(f) = 
l!p^{m|l|}\pa^{[l]}(f). $$
In particular, we have $\pa^{\lara{l}_{-m}}(f) \to 0$ as $|l|\to\infty$. 
\end{remark}

Next we define the notion of $(-m)$-PD-stratification and 
compare it with the notion of left $\cD_{X/S}^{(-m)}$-module. 

\begin{definition}
A $(-m)$-PD-stratification on an $\cO_X$-module $\cE$ is 
a compatible family of $\cP^k_{X,(-m)}$-linear isomorphisms 
$\{\epsilon_k: \cP^k_{X,(-m)} \otimes_{\cO_X} \cE \os{=}{\lra} 
\cE \otimes_{\cO_X} \cP^k_{X,(-m)}\}_k$ with $\epsilon_0 = \id$ 
such that the following diagram is commutative for any $k,k' \in \N:$
\[\xymatrix{
\cP^k_{X,(-m)} \otimes_{\cO_X} \cP^{k'}_{X,(-m)} \otimes_{\cO_X} \cE 
\ar[r]_{\id \otimes \epsilon_{k'}}
\ar[rd]_{\delta^{k,k',*}(\epsilon_{k+k'})} & 
\cP^k_{X,(-m)} \otimes_{\cO_X} \cE \otimes_{\cO_X} 
\cP^{k'}_{X,(-m)} 
\ar[d]_{\epsilon_k \otimes \id} \\ 
& \cE \otimes_{\cO_X} \cP^k_{X,(-m)} \otimes_{\cO_X} \cP^{k'}_{X,(-m)}.
}\]
\end{definition}

The conditions put on $\{\epsilon_k\}_k$ in the above definition is 
called the cocycle condition. 
It is easy to see that the cocycle condition is equivalent to the 
condition $q^{k*}_{02}(\epsilon_k) = q^{k*}_{01}(\epsilon_k) \circ 
q^{k*}_{12}(\epsilon_k)$ for $k \in \N$, where $q_{ij}^k$ denotes 
the homomorphism $\cP^k_{X,(-m)} \lra \cP^k_{X,(-m)}(2)$ induced by 
the $(i,j)$-th projection $X^3 \lra X^2$. \par 
We have the following equivalence, which is an analogue of 
\cite[2.3.2]{berthelotI}: 

\begin{proposition}\label{5}
For an $\cO_X$-module $\cE$, the following three data are equivalent. 
\begin{enumerate}
\item[\rm (a)] 
A left $\cD_{X/S}^{(-m)}$-module structure on $\cE$ which extends the 
given $\cO_X$-module structure. 
\item[\rm (b)] A compatible family of $\cO_X$-linear 
homomorphisms $\{\theta_k: \cE \lra \cE \otimes_{\cO_X} \cP^k_{X,(-m)}\}_k$
 $($where we regard $\cE \otimes_{\cO_X} \cP^k_{X,(-m)}$ as 
$\cO_X$-module by using the right $\cO_X$-module structure of 
$\cP^k_{X,(-m)})$ with $\theta_0=\id$ such that the following diagram 
is commutative for any $k,k' \in \N:$ 
\begin{equation}\label{a7x}
\begin{CD}
\cE \otimes_{\cO_X} \cP^{k+k'}_{X,(-m)} 
@>{\id \otimes \delta^{k,k'}}>> 
\cE \otimes_{\cO_X} \cP^{k}_{X,(-m)} \otimes_{\cO_X} \cP^{k'}_{X,(-m)} \\ 
@A{\theta_{k+k'}}AA @A{\theta_{k} \otimes \id}AA \\
\cE @>{\theta_{k'}}>> \cE \otimes_{\cO_X} \cP^{k'}_{X,(-m)}.
\end{CD}
\end{equation}
\item[\rm (c)] A $(-m)$-PD-stratification $\{\epsilon_k\}_k$ on 
$\cE$. 
\end{enumerate}
\end{proposition}

\begin{proof}
Since the proof is identical with the classical case, 
we only give a brief sketch. The data in (a) is equivalent to 
a compatible family of homomorphisms 
$\mu_k: \cD_{X/S,k}^{(-m)} \otimes_{\cO_X} \cE \lra \cE$ 
satisfying the condition coming from the product structure 
of $\cD_{X/S}^{(-m)}$, and $\mu_k$ induces the homomorphism 
$$ \theta_k: \cE \lra \Inthom_{\cO_X}(\cD_{X/S,k}^{(-m)},\cE) = 
\cE \otimes_{\cO_X} \cP^k_{X,(-m)} $$ 
which satisfies the conditions in (b). So the data in (a) gives the data 
(b), and we see easily that they are in fact equivalent. 
When we are given the data in (b), we obtain the 
$\cP^k_{X,(-m)}$-linear homomorphism 
$\epsilon_k: \cP^k_{X,(-m)} \otimes_{\cO_X} \cE \lra 
\cE \otimes_{\cO_X} \cP^k_{X,(-m)}$ by taking 
$\cP^k_{X,(-m)}$-linearization of $\theta_k$. Since 
$\epsilon_k(1 \otimes x) = \theta_k(x)$ is written locally as 
$\sum_{|l|\leq k} \pa^{\lara{l}}(x) \otimes (\tau/p^m)^{[l]}$, 
we see that $\epsilon_k$ is actually an isomorphism because the inverse 
of it is given locally by 
$x \otimes 1 \mapsto \sum_{|l|\leq k} (-1)^{|l|} (\tau/p^m)^{[l]} \otimes 
\pa^{\lara{l}}(x)$. The cocycle condition for $\{\epsilon_k\}_k$ follows from 
the commutative diagram \eqref{a7x} 
for $\{\theta_k\}_k$ and so the data in (b) gives the data in (c). 
Again we see easily that they are in fact equivalent. 
\end{proof}

Next we relate the notion of left $\cD^{(-m)}_{X/S}$-modules and 
that of $p^m$-connections. Let $X \lra S$ be as above. 
Recall that 
a $p^m$-connection on an $\cO_X$-module $\cE$ 
is an additive map $\nabla: \cE \lra \cE \otimes_{\cO_X} \Omega^1_{X/S}$ 
satisfying $\nabla(fe) = f\nabla(e) + p^me \otimes df$ 
$\,(e \in \cE, f \in \cO_X)$. 
To give another description of $p^m$-connection, let us put 
$J^1_{X/S} := \Ker(\cP^1_{X,(-m)} \lra \cO_X)$. Then we have a 
natural map $\alpha: \Omega^1_{X/S} \lra J^1_{X/S}$ induced by the 
map $\cP^1_{X,(0)} \lra \cP^1_{X,(-m)}$, and locally $\alpha$ is 
given by $dt_i = \tau_i \mapsto p^m(\tau_i/p^m)$. So 
$\alpha$ is injective and the image is equal to $p^mJ^1_{X/S}$. 
Hence we have the unique isomorphism $\beta: \Omega^1_{X/S} \os{=}{\lra} 
J^1_{X/S}$ satisfying $p^m\beta=\alpha$. Via the identification by $\beta$, 
a $p^m$-connection on $\cE$ is equivalent to an additive map 
$\nabla: \cE \lra \cE \otimes_{\cO_X} J^1_{X/S}$ 
satisfying $\nabla(fe) = f\nabla(e) + e \otimes df$ for any 
$e \in \cE, f \in \cO_X$. (Attention: the element $df \in J^1_{X/S}$ here 
is the element $1 \otimes f - f \otimes 1 \in \cP^1_{X,(-m)}$, 
not the element $\beta(1 \otimes f - f \otimes 1)$.) \par 
The following proposition is the analogue of \cite[2.9]{bo}. 

\begin{proposition}\label{7}
Let $S,X$ be as above. 
For an $\cO_X$-module $\cE$, the following data are equivalent$:$ 
\begin{enumerate}
\item[\rm (a)] 
A $p^m$-connection $\nabla: \cE \lra \cE \otimes_{\cO_X} J^1_{X/S}$ on $\cE$. 
\item[\rm (b)] 
A $\cP^1_{X,(-m)}$-linear isomorphism $\epsilon_1: 
\cP^1_{X,(-m)} \otimes_{\cO_X} \cE \os{=}{\lra} \cE \otimes_{\cO_X} 
\cP^1_{X,(-m)}$ which is equal to identity modulo $J^1_{X/S}$. 
\end{enumerate}
\end{proposition}

\begin{proof}
Since the proof is again identical with \cite[2.9]{bo}, we only give 
a brief sketch. First assume that we are given the isomorphism $\epsilon_1$ 
as in (b). Then, if we 
define $\nabla: \cE \lra \cE \otimes_{\cO_X} J^1_{X/S}$ by 
$\nabla(e) := \epsilon(1 \otimes e) - e \otimes 1$, it gives a 
$p^m$-connection. Conversely, if we are given a $p^m$-connection 
$\nabla: \cE \lra \cE \otimes_{\cO_X} J^1_{X/S}$, 
let us define the $\cP^1_{X,(-m)}$-linear homomorphism 
$\epsilon_1: \cP^1_{X,(-m)} \otimes_{\cO_X} \cE \os{=}{\lra} \cE \otimes_{\cO_X} \cP^1_{X,(-m)}$ by $\epsilon(1\otimes e) = \nabla(e) + 
e \otimes 1$. Then it is easy to see that $\epsilon_1$ is 
 equal to identity modulo $J^1_{X/S}$. To show that $\epsilon_1$ is 
an isomorphism, let us consider the isomomorphism 
$t: \cP^1_{X,(-m)} \lra \cP^1_{X,(-m)}$ induced by the morphism 
$X^2 \lra X^2; (x,y) \mapsto (y,x)$ and let $s: 
\cE \otimes_{\cO_X} \cP^1_{X,(-m)} \lra \cP^1_{X,(-m)} \otimes_{\cO_X} \cE$ 
be the isomorphism $x \otimes \xi \mapsto t(\xi) \otimes x$. 
Then we see that $(s\circ \epsilon_1)^2:  \cP^1_{X,(-m)} \otimes_{\cO_X} \cE 
\lra  \cP^1_{X,(-m)} \otimes_{\cO_X} \cE$ is a $\cP^1_{X,(-m)}$-linear 
endomorphism which is equal to the identity modulo $J^1_{X/S}$. Hence 
it is an isomorphism and we see from this that $\epsilon_1$ is also 
an isomorphism. 
\end{proof}

As for the integrability, we have the following proposition. 

\begin{proposition}\label{8}
Let $\cE$ be an $\cO_X$-module and let $\nabla: \cE \lra 
\cE \otimes_{\cO_X} 
\Omega^1_{X/S} \os{\id \otimes \beta, =}{\lra}
\cE \otimes_{\cO_X} 
J^1_{X/S}$ be a $p^m$-connection. Let 
$\epsilon_1: \cP^1_{X,(-m)} \otimes_{\cO_X} \cE 
\lra \cE \otimes_{\cO_X} \cP^1_{X,(-m)}$ be the 
$\cP^1_{X,(-m)}$-linear isomorphism corresponding to $\nabla$ by 
the equivalence in Proposition \ref{7} and let 
$\mu_1: \cD_{X/S,1}^{(-m)} \otimes_{\cO_X} \cE \lra \cE$ be the 
homomorphism induced by the composite 
$$ 
\cE \hra \cP^1_{X,(-1)} \otimes_{\cO_X} \cE 
\os{\epsilon_1}{\lra} 
\cE \otimes_{\cO_X} \cP^1_{X,(-1)} \cong 
\Inthom_{\cO_X}(\cD_{X/S,1}^{(-m)},\cE).
$$ 
Then the following conditions are equivalent. 
\begin{enumerate}
\item 
$(\cE,\nabla)$ is integrable. 
\item 
$\mu_1$ is $($uniquely$)$ extendable to a $\cD_{X/S}^{(-m)}$-module 
structure on $\cE$ which extends the 
given $\cO_X$-module structure. 
\end{enumerate}
\end{proposition}

\begin{proof}
We may work locally. So we can write 
$\nabla(e) = \sum_i \theta_i(e) dt_i \os{\id \otimes \beta}{\mapsto}
\sum_i \theta_i(e) (\tau_i/p^m)$, using local coordinate. 
Then we have $\mu_1(\pa_i \otimes e) = \theta_i(e)$. \par 
First assume the condition (2). Then, since $[\pa_i,\pa_j](e)=0$ for any 
$e \in \cE$, we have $[\theta_i,\theta_j](e)=0$ for any $e$ and so 
$(\cE,\nabla)$ is integrable. So the condition (1) is satisfied. \par 
On the other hand, let us assume the condition (1). Then, 
we define the action of $\pa^{\lara{k}} \in \cD_{X/S}^{(-m)}$ 
on $e \in \cE$ by $\pa^{\lara{k}}(e) := \prod_{i=1}^{d}\theta_i^{k_i}(e)$, 
where $k=(k_1,...,k_d)$. To see that this action acturally defines 
a $\cD_{X/S}^{(-m)}$-module structure on $\cE$, we have to check 
the following equalities for $e \in \cE$ and $f \in \cO_X$ 
(see (3), (4) in Proposition \ref{4}): 
\begin{equation}\label{eq4-3}
\pa^{\lara{k}}\pa^{\lara{k'}}(e) = \pa^{\lara{k+k'}}(e), 
\end{equation}
\begin{equation}\label{eq4-4}
\pa^{\lara{k}}(fe) = 
\sum_{k'+k''=k} \begin{pmatrix} k \\ k' \end{pmatrix} 
\pa^{\lara{k'}}(f) \pa^{\lara{k''}}(e). 
\end{equation}
By the definition of the action of $\pa^{\lara{k}}$'s on $e$ given above, 
the equality \eqref{eq4-3} is reduced to the equality 
$\theta_i\theta_j(e) = \theta_j\theta_i(e)$, that is, the integrability 
of $(\cE,\nabla)$. In view of the equality \eqref{eq4-3}, the proof of 
the equality \eqref{eq4-4} is reduced to the case $|k|=1$, and in this 
case, it is rewritten as 
$$ \theta_i(fe) = f \theta_i(e) + \pa_i(f) e \quad (1 \leq i \leq d)$$ 
and it is equlvalent to the equality 
$\nabla(fe) = f \nabla(e) + e \otimes df$ in 
$\cE \otimes_{\cO_X} J^1_{X/S}$, which is true by the definition 
of $p^m$-connection. So we have the well-defined 
$\cD_{X/S}^{(-m)}$-module structure on $\cE$ and hence the condition (2) 
is satisfied. So we are done. 
\end{proof}

\begin{corollary}\label{8cor}
For an $\cO_X$-module $\cE$, the following three data are equivalent. 
\begin{enumerate}
\item[\rm (a)] An integrable $p^m$-connection on $\cE$. 
\item[\rm (b)] A $\cD_{X/S}^{(m)}$-module structure on $\cE$ 
which extends the 
given $\cO_X$-module structure. 
\item[\rm (c)] A $(-m)$-PD-stratification $\{\epsilon_k\}_k$ on 
$\cE$. 
\end{enumerate}
In particular, we have the equivalence
$$\MIC^{(m)}(X) \os{=}{\lra} \text{\rm (left $\cD^{(-m)}_{X/S}$-modules)}$$ 
and it induces the equivalence 
$$\MIC^{(m)}(X)_n \os{=}{\lra} \text{\rm (left $\cD^{(-m)}_{X/S} \otimes 
\Z/p^n\Z$-modules)}.$$ 
\end{corollary}

\begin{proof}
It suffices to prove the equivalence of (a) and (b). 
When we are given an integrable $p^m$-connection on $\cE$, we have 
the desired $\cD_{X/S}^{(m)}$-module structure on $\cE$ thanks to 
Proposition \ref{8}. Conversely, when we are given a 
$\cD_{X/S}^{(m)}$-module structure on $\cE$ which extends the 
given $\cO_X$-module structure, we have the induced homomorphism 
$\mu_1: \cD_{X/S,1}^{(-m)} \otimes_{\cO_X} \cE \lra \cE$. It gives 
the $\cO_X$-linear homomorphism 
$$ 
\cE \lra \Inthom_{\cO_X}(\cD_{X/S,1}^{(-m)},\cE) \cong 
\cE \otimes_{\cO_X} \cP^1_{X,(-m)}$$ 
(where the $\cO_X$-module structure on the target is induced by 
the right $\cO_X$-module structure on $\cP^1_{X,(-m)}$), and by 
taking the $\cP^1_{X,(-m)}$-linearization of it, we obtain the 
homomorphism $\epsilon_1: \cP^1_{X,(-m)} \otimes_{\cO_X} \cE 
\lra \cE \otimes_{\cO_X} \cP^1_{X,(-m)}$ which is equal to the 
identity modulo $J^1_{X/S}$. It is automatically an isomorphism 
by the last argument in the proof of Proposition \ref{7}, and it gives 
a $p^m$-connection $\nabla$ by Propisition \ref{7}. Then, 
$\nabla$ gives rise to the homomorphism $\mu_1$ by the recipe given in 
the statement of Proposition \ref{8}. Since $\mu_1$ is extendable to 
the $\cD_{X/S}^{(-m)}$-module structure by assumption, we see 
by Proposition \ref{8} that $\nabla$ is integrable. So we obtain 
the integrable $p^m$-connection $\nabla$ and so we are done. 
\end{proof}

Next we give a $\cD$-module theoretic interpretation of 
the quasi-nilpotence for objects in $\MIC^{(m)}(X)_n$. 
The following proposition is the 
analogue of \cite[2.3.7]{berthelotI}. 

\begin{proposition}\label{9}
Let $f: X \lra S$ be a smooth morphism of $p$-adic formal 
schemes flat over $\Z_p$. 
Let $m \in \N$ and 
let $\cE := (\cE,\nabla)$ be an object in $\MIC^{(m)}(X)_n$, regarded as a 
left $\cD_{X/S}^{(-m)} \otimes \Z/p^n\Z$-module. 
Then the following conditions are equivalent. 
\begin{enumerate}
\item[\rm (a)] 
$(\cE, \nabla)$ is quasi-nilpotent as an object in $\MIC^{(m)}(X)_n$. 
\item[\rm (b)] 
Locally on $X$, $f$ admits a local coordinate 
such that 
the following condition is satisfied$:$ 
For any $e \in \cE$, there exists some 
$N \in \N$ such that $\pa^{\lara{k}}(e) = 0$ for any 
$k$ with $|k| \geq N$, where $\pa^{\lara{k}}$ is the element 
in $\cD_{X/S}^{(-m)}$ 
defined by using the fixed local coordinate. 
\item[\rm (c)] 
The condition given in {\rm (b)} is satisfied for any local coordinate. 
\item[\rm (d)] 
There exists $($uniquely$)$ a $\cP_{X,(-m)}$-linear isomorphism 
$\epsilon: \cP_{X,(-m)} \otimes_{\cO_{X}} \cE \os{=}{\lra} 
\cE \otimes_{\cO_{X}} \cP_{X,(-m)}$ satisfying the cocycle condition 
on $\cP_{X,(-m)}(2)$ which induces the $(-m)$-PD stratification 
$\{\epsilon_k\}_k$ on $\cE$ associated to the 
 $\cD_{X/S}^{(-m)}$-module structure on $\cE$ via 
Proposition \ref{5}. 
\end{enumerate} 
$($
We call the isomorphism $\epsilon$ in {\rm (d)} the $(-m)$-HPD-stratification 
 associated to $\cE.)$ 
\end{proposition}

\begin{proof}
The proof is similar to that of 
 \cite[2.3.7]{berthelotI}. 
First, let us work locally on $X$, take a local coordinate 
$t_1, ..., t_d$ of $f$ and write $\nabla$ as 
$\nabla(e) = \sum_i \theta_i(e)dt_i$. Then, in the notation in (b), 
we have $\theta^k = \pa^{\lara{k}}$ for any $k \in \N^d$. Hence 
we have the equivalence of (a) and (b). 
When the condition (b) is satisfied, we can define 
the morphism $\theta: \cE \lra \ \cE \otimes_{\cO_{X}} \cP_{X,(-m)}$ 
by $\theta(e) = \sum_{k} \pa^{\lara{k}}(e) \otimes (\tau/p^m)^{[k]}$ 
and by $\cP_{X,(-m)}$-linearizing it, we obtain the homomorphism 
$\epsilon: \cP_{X,(-m)} \otimes_{\cO_{X}} \cE \lra 
\cE \otimes_{\cO_{X}} \cP_{X,(-m)}$ which induces the stratification 
$\{\epsilon_k\}_k$. The cocycle condition for 
$\epsilon$ follows from that for $\{\epsilon_k\}_k$ and the uniqueness 
is clear. Also, if we define $\theta':  \cE \lra  \cP_{X,(-m)} 
\otimes_{cO_{X}} \cE$ by $\theta'(e) := \sum_{k} (-1)^{|k|} 
(\tau/p^m)^{[k]} \otimes \pa^{\lara{k}}(e)$, we see that the 
$\cP_{X,(-m)}$-linearization of it gives the inverse of $\epsilon$. 
So $\epsilon$ is an isomorphism and thus defines a 
$(-m)$-HPD-stratification. Conversely, if we are given 
a $(-m)$-HPD-stratification $\epsilon$ associated to $\cE$, 
the coefficient of $(\tau/p^m)^{[k]}$ of the elememt 
$\epsilon(1 \otimes e) \in \cE \otimes_{\cO_{X}} \cP_{X,(-m)} = 
\bigoplus_k \cE (\tau/p^m)^{[k]}$ 
is equal to $\pa^{\lara{k}}(e)$, by Proposition \ref{5}. 
Hence the condition (b) is satisfied. Finally, since the 
condition (d) is independent of the choice of the coordinate, 
we have the equivalence of the conditions (c) and (d). 
\end{proof}

\begin{definition}\label{def2.10}
Let $f:X \lra S$ be as above. Then, a left 
$\cD_{X/S}^{(-m)}$-module $\cE$ is said to be quasi-nilpotent 
if it is $p^n$-torsion for some $n$ and that it satisfies the 
condition $(d)$ of Proposition \ref{9}. By Proposition \ref{9}, 
we have the equivalence 
$$ \MIC^{(m)}(X)^{\qn} \os{=}{\lra} 
\text{\rm (quasi-nilpotent left $\cD^{(-m)}_{X/S}$-modules)}, $$ 
which is induced by 
$$ \MIC^{(m)}(X)_n^{\qn} \os{=}{\lra} 
\text{\rm (quasi-nilpotent left $\cD^{(-m)}_{X/S} \otimes \Z/p^n\Z$-modules)} 
\,\,\,\,(n \in \N).$$ 
\end{definition}

Next we give a definition of the inverse image functor for 
left $\cD^{(-m)}_{-/-}$-modules. Let  
\begin{equation}\label{pbdiag1}
\begin{CD}
X' @>f>> X \\ 
@VVV @VVV \\ 
S' @>>> S 
\end{CD}
\end{equation} 
be a commutative diagram of 
$p$-adic formal schemes flat over $\Spf \Z_p$ such that 
the vertical arrows are smooth. 
Then, for any $m,k \in \N$, it induces the commutative 
diagram 
\begin{equation}\label{a}
\begin{CD}
P^k_{X',(-m)} @>{\subset}>> P_{X',(-m)} @>{p'_i}>> X' \\ 
@V{g^k}VV @V{g}VV @VfVV \\ 
P^k_{X,(-m)} @>{\subset}>> P_{X,(-m)} @>{p_i}>> X, 
\end{CD}
\end{equation}
for $i=0,1$, where $p'_i, p_i$ denotes the morphism 
induced by the $i$-th projection ${X'}^2 \lra X', X^2 \lra X$, 
respectively. So, if $\cE$ is an $\cO_X$-module endowed with 
a $(-m)$-PD-stratification $\{\epsilon_k\}_k$, 
$f^*\cE$ is naturally endowed with the 
 $(-m)$-PD-stratification $\{{g^k}^*\epsilon_k\}_k$. Hence, in 
view of Proposition \ref{5}, 
we have the functor 
\begin{align}
f^*: & \text{(left $\cD_{X/S}^{(-m)}$-modules)} \lra 
\text{(left $\cD_{X'/S'}^{(-m)}$-modules)}; \label{iidf} \\ 
& \hspace{1cm}
(\cE,\{\epsilon_k\}_k) \hspace{0.8cm}\mapsto\hspace{0.8cm} 
(f^*\cE, \{{g^k}^*\epsilon_k\}_k), 
\nonumber 
\end{align}
and this induces also the functor 
\begin{equation}\label{iidf_n}
f^*: \text{(left $\cD_{X/S}^{(-m)} \otimes \Z/p^n\Z$-modules)} \lra 
\text{(left $\cD_{X'/S'}^{(-m)} \otimes \Z/p^n\Z$-modules)}. 
\end{equation}

As for the quasi-nilpotence, we have the following: 

\begin{proposition}\label{qnqn}
With the above notation, assume that $\cE$ is quasi-nilpotent. 
Then $f^*\cE$ is also quasi-nilpotent. 
\end{proposition}

\begin{proof}
When $\cE$ is quasi-nilpotent, the 
$(-m)$-PD-stratification $\{\epsilon_k\}_k$ associated to $\cE$ 
is induced from a $(-m)$-HPD-stratification $\epsilon$. Then 
the $(-m)$-PD-stratification $\{{g^k}^*\epsilon_k\}_k$ associated to 
$f^*\cE$ is induced from the $(-m)$-HPD-stratification $g^*\epsilon$ by 
the commutativity of the diagram \eqref{a}. 
So $f^*\cE$ is also quasi-nilpotent. 
\end{proof}

The inverse image functor here is equivalent to the 
inverse image functor in the previous section in the 
following sense. 

\begin{proposition}\label{coin1}
With the above notation, the inverse image functor 
\eqref{iidf} is equal to the inverse image functor 
 $f^*: \MIC^{(m)}(X) \lra \MIC^{(m)}(X')$ defined in the previous 
section via the equivalence in Corollary \ref{8cor}. 
$($Hence the inverse image functor \eqref{iidf_n} is equal to 
the inverse image functor 
 $f^*: \MIC^{(m)}(X)_n \lra \MIC^{(m)}(X')_n$ 
defined in the previous 
section.$)$ 
\end{proposition}

\begin{proof}
Assume given an object $(\cE,\nabla) \in \MIC^{(m)}(X)$ and 
let $(f^*\cE, f^*\nabla) \in \MIC^{(m)}(X')$ be the inverse image 
of it defined in the previous section. On the other hand, 
let $(\cE,\{\epsilon_k\})$ be the $(-m)$-PD-stratification 
associated to  $(\cE,\nabla)$, let 
$(f^*\cE, \{{g^k}^*\epsilon_k\})$ be the inverse image of it defined 
above and let $(f^*\cE, \ti{\nabla})$ be the object in 
$\MIC^{(m)}(X')$ associated to $(f^*\cE, \{{g^k}^*\epsilon_k\})$ 
via the equivalence in Corollary \ref{8cor}. Since 
the underlying $\cO_{X'}$-module $f^*\cE$ of 
$(f^*\cE, f^*\nabla)$ and $(f^*\cE, \ti{\nabla})$ are the same, 
it suffices to prove the coincidence of $p^m$-connections 
$f^*\nabla$ and $\ti{\nabla}$. To see this, we may work locally. \par 
Take a local section $e \in \cE$ and let us put 
$\nabla(e) = \sum_i e_i da_i$. Then we have 
$f^*\nabla(e) = \sum f^*e_i df^*(a_i)$. \par 
On the other hand, 
if we denote the composite 
$\cE \os{\nabla}{\lra} \cE \otimes \Omega^1_{X/S} 
\os{\id \otimes \beta}{\lra}
\cE \otimes J^1_{X/S}$ by $\nabla'$, we have 
$\nabla'(e) = \sum_i e_i (da_i/p^m)$. Hence $\epsilon^1: 
\cP^1_{X,(-m)} \otimes \cE \lra \cE \otimes \cP^1_{X,(-m)}$ 
is written as 
\begin{align*}
\epsilon(1 \otimes e) 
= e \otimes 1 + \nabla'(e) = 
e \otimes 1 + \sum_i e_i \otimes da_i/p^m. 
\end{align*}
Since 
${g^1}^*: \cP^1_{X,(-m)} \lra \cP^1_{X',(-m)}$ sends 
$da_i/p^m$ to $df^*(a_i)/p^m$, we have 
$$
{g^1}^*\epsilon(1 \otimes f^*e) 
= f^*e \otimes 1 + \sum_i f^*e_i \otimes {g^1}^*(da_i/p^m) = 
 f^*e \otimes 1 + \sum_i f^*e_i \otimes df^*(a_i)/p^m, $$
and so we have $\ti{\nabla}(f^*e) = 
\sum_i f^*e_i \otimes df^*(a_i)/p^m$. Therefore we have 
$f^*\nabla = \ti{\nabla}$, as desired. 
\end{proof}

\begin{remark}
Proposition \ref{qnqn} together with Proposition \ref{coin1} 
gives another proof of Proposition \ref{functqn} in the case of 
$p$-adic formal schemes flat over $\Z_p$ (at least in the case 
$m \geq 1$). 
\end{remark}

Next we define the level raising inverse image functor for 
$\cD^{(-)}_{-/-}$-modules. The following proposition is an 
analogue of \cite[2.2.2]{berthelotII}. 

\begin{proposition}\label{13}
In the situation of $\Hyp(\infty,\infty,\infty)$, 
$F$ induces naturally a PD-morphism 
$\Phi: P_{X,(-m+1)} \lra P_{X^{(1)},(-m)}$ $($with respect to the PD-ideal 
on the defining ideal of $X \hra P_{X,(-m+1)}, X^{(1)} \hra 
 P_{X^{(1)},(-m)}$, respectively$)$. 
\end{proposition}

\begin{proof}
We may work locally. So we may assume that 
there exists a local parameter $t_1, ..., t_d$ of $X$. 
Let us put $t'_i := 1 \otimes t_i \in \Gamma(X^{(1)}, 
\cO_{X^{(1)}})$ for $1 \leq i \leq d$. 
Then we have $F^*(t'_i) = t_i^p + pa_i$ for some 
$a_i \in \cO_X \,(1 \leq i \leq d)$. Let us put 
$\tau_i := 1 \otimes t_i - t_i \otimes 1 \in \cO_{X^2}$, 
$\tau'_i := 1 \otimes t'_i - t'_i \otimes 1 \in \cO_{(X^{(1)})^2}$ and 
let us denote the morphism $F \times F: X^2 \lra (X^{(1)})^2$ simply 
by $F^2$. Then we have 
\begin{align}
{F^2}^*(\tau'_i) & = 
1 \otimes t_i^p - t_i^p \otimes 1 + p(1 \otimes a_i - a_i \otimes 1) 
\label{calc} \\ 
& = 
(\tau_i + t_i \otimes 1)^p - 
t_i^p \otimes 1 + p(1 \otimes a_i - a_i \otimes 1) \nonumber \\ 
& = 
\tau_i^p + \sum_{k=1}^{p-1} \begin{pmatrix} p \\ k \end{pmatrix} 
t_i^{p-k}\tau_i^k + p(1 \otimes a_i - a_i \otimes 1). \nonumber 
\end{align}
Hence there exists an element $\sigma_i \in I := \Ker(\cO_{X^2} \ra \cO_X)$ 
such that ${F^2}^*(\tau'_i) = \tau_i^p + p\sigma_i$. So, 
when $m \geq 2$, the image of ${F^2}^*(\tau'_i)$ in 
$\cO_{T_{X,(-m+1)}}$ belongs to $I^p + pI \subseteq 
(p^{m-1}\cO_{T_{X,(-m+1)}})^p + p(p^{m-1}\cO_{T_{X,(-m+1)}}) 
= p^m\cO_{T_{X,(-m+1)}}$. So, by the universality of formal 
blow-up, $F^2$ induces the morphism $T_{X,(-m+1)} \lra T_{X^{(1)},(-m)}$ 
and by the universality of the PD-envelope, it induces the 
PD-morphism $\Phi: P_{X,(-m+1)} \lra P_{X',(-m)}$, as desired. 
When $m=1$, the image of ${F^2}^*(\tau'_i)$ in 
$\cO_{P_{X,(0)}}$ is equal to $p! \tau_i^{[p]} + p\sigma_i$ and so 
it belongs to $p\cO_{P_{X,(0)}}$. Hence $F^2$ induces 
the morphism $P_{X,(0)} \lra T_{X^{(1)},(-1)}$ 
and then it induces 
the PD-morphism $\Phi: P_{X,(0)} \lra P_{X^{(1)},(-1)}$, as desired. 
\end{proof}

\begin{remark}\label{13rem}
In the same way as the above proof, we can prove also that, 
for $r \in \N$, the morphism $F^{r+1}: X^{r+1} \lra {X'}^{r+1}$ 
naturally induces the PD-morphism 
$\Phi: P_{X,(-m+1)}(r) \lra P_{X',(-m)}(r)$.  
\end{remark}

Let the situation be as in $\Hyp(\infty,\infty,\infty)$ 
and let $m \in \N$. 
Then, by Proposition \ref{13}, we have the commutative diagrams 
\begin{equation}\label{dd}
\begin{CD}
P^k_{X,(-m+1)} @>{\subset}>> P_{X,(-m+1)} @>{p_i}>> X \\ 
@V{\Phi^k}VV @V{\Phi}VV @V{F}VV \\ 
P^k_{X^{(1)},(-m)} @>{\subset}>> P_{X^{(1)},(-m)} @>{p^{(1)}_i}>> X^{(1)} 
\end{CD}
\end{equation}
for $i=0,1$, where $p_i, p^{(1)}_i$ denotes the morphism 
induced by the $i$-th projection ${X}^2 \lra X, 
(X^{(1)})^2 \lra X^{(1)}$ 
respectively and $\Phi^k$ is the morphism naturally induced by 
$\Phi$. So, if $\cE$ is an $\cO_{X^{(1)}}$-module endowed with 
a $(-m)$-PD-stratification $\{\epsilon_k\}_k$, $F^*\cE$ is 
endowed with a $(-m+1)$-PD-stratification 
$\{{\Phi^k}^*\epsilon_k\}$. Hence, in view of Proposition \ref{5}, 
we have the functor 
\begin{align}
F^*: & \text{(left $\cD_{X^{(1)}/S}^{(-m)}$-modules)} \lra 
\text{(left $\cD_{X/S}^{(-m+1)}$-modules)}; \label{lriidf} \\ 
& \hspace{1cm}
(\cE,\{\epsilon_k\}_k) \hspace{0.8cm}\mapsto\hspace{0.8cm} 
(f^*\cE, \{{\Phi^k}^*\epsilon_k\}_k), 
\nonumber 
\end{align}
and this induces also the functor 
\begin{equation}\label{lriidf_n}
F^*: \text{(left $\cD_{X^{(1)}/S}^{(-m)} \otimes \Z/p^n\Z$-modules)} \lra 
\text{(left $\cD_{X/S}^{(-m+1)} \otimes \Z/p^n\Z$-modules)}. 
\end{equation}

By the existence of the diagram \eqref{dd}, we can prove the 
following in the same way as Proposition \ref{qnqn} (so we omit the 
proof): 

\begin{proposition}\label{qnqnqnqn}
With the above notation, assume that $\cE$ is quasi-nilpotent. Then 
$F^*\cE$ is also quasi-nilpotent. 
\end{proposition}

The inverse image functor here is equivalent to the 
inverse image functor in the previous section in the 
following sense. 

\begin{proposition}\label{coin2}
With the above notation, the inverse image functor 
\eqref{lriidf} is equal to the level raising inverse image functor 
 $F^*: \MIC^{(m)}(X^{(1)}) \lra \MIC^{(m-1)}(X)$ defined in the previous 
section via the equivalence in Corollary \ref{8cor}. 
$($Hence the inverse image functor \eqref{lriidf_n} is equal to 
the inverse image functor 
 $F^*: \MIC^{(m)}(X^{(1)})_n \allowbreak \lra \MIC^{(m-1)}(X)_n$ 
defined in the previous 
section.$)$ 
\end{proposition}

\begin{proof}
Assume given an object $(\cE,\nabla) \in \MIC^{(m)}(X^{(1)})$ and 
let $(F^*\cE, F^*\nabla)$ be the level raising inverse image 
of it defined in the previous section. On the other hand, 
let $(\cE,\{\epsilon_k\})$ be the $(-m)$-PD-stratification 
associated to $(\cE,\nabla)$, let 
$(F^*\cE, \{{\Phi^k}^*\epsilon_k\})$ be the level raising 
inverse image of it defined above and let 
$(F^*\cE, \ti{\nabla})$ be the object in 
$\MIC^{(m-1)}(X)$ associated to $(F^*\cE, \{{\Phi^k}^*\epsilon_k\})$ 
via the equivalence in Corollary \ref{8cor}. Since 
the underlying $\cO_{X}$-module $F^*\cE$ of 
$(F^*\cE, F^*\nabla)$ and $(F^*\cE, \ti{\nabla})$ are the same, 
it suffices to prove the coincidence of $p^m$-connections 
$F^*\nabla$ and $\ti{\nabla}$. To see this, we may work locally. 
So we can take a local coordinate $t_1, ..., t_d$ of $X^{(1)}$ over 
$S$. Let us put $t'_i, a_i$ as in the proof of 
Proposition \ref{13}. 
Take a local section $e \in \cE$ and let us put 
$\nabla(e) = \sum_i e_i dt'_i$. Then, by definition, we have 
$F^*\nabla(e) = \sum F^*e_i (t_i^{p-1}dt_i + da_i)$. \par 
On the other hand, 
if we denote the composite 
$\cE \os{\nabla}{\lra} \cE \otimes \Omega^1_{X/S} 
\os{\id \otimes \beta}{\lra}
\cE \otimes J^1_{X/S}$ by $\nabla'$, we have 
$\nabla'(e) = \sum_i e_i (dt'_i/p^m)$. Hence $\epsilon^1: 
\cP^1_{X,(-m)} \otimes \cE \lra \cE \otimes \cP^1_{X,(-m)}$ 
is written as 
\begin{align*}
\epsilon(1 \otimes e) 
= e \otimes 1 + \nabla'(e) = 
e \otimes 1 + \sum_i e_i \otimes dt'_i/p^m. 
\end{align*}
Since 
${\Phi^1}^*: \cP^1_{X,(-m)} \lra \cP^1_{X',(-m)}$ sends 
$dt'_i/p^m$ to $t_i^{p-1}(dt_i/p^{m-1}) + da_i/p^{m-1}$ 
by the calculation \ref{calc}, we have 
\begin{align*}
{\Phi^1}^*\epsilon(1 \otimes F^*e) 
& = F^*e \otimes 1 + \sum_i F^*e_i \otimes {\Phi^1}^*(da_i/p^m) \\ 
& = 
 F^*e \otimes 1 + \sum_i F^*e_i \otimes 
(t_i^{p-1}(dt_i/p^{m-1}) + da_i/p^{m-1}), 
\end{align*} 
and so we have $\ti{\nabla}(F^*e) = 
\sum_i f^*e_i \otimes (t_i^{p-1}dt_i + da_i)$. Therefore we have 
$f^*\nabla = \ti{\nabla}$, as desired. 
\end{proof}

\begin{remark}
Proposition \ref{qnqnqnqn} together with Proposition \ref{coin2} 
gives another proof of Proposition \ref{functqn2} under 
$\Hyp(\infty,\infty,\infty)$, as promised in Remark \ref{qnatode}.  
\end{remark}

\subsection{The case of schemes over $\Z/p^n\Z$} 

In the previous subsection, we defined the sheaf of 
$p$-adic differential operators of level $-m$ for smooth 
morphisms of $p$-adic formal schemes flat over $\Z_p$. The construction 
there does not work well for smooth 
morphisms of schemes flat over $\Z/p^n\Z$ because we needed the formal 
blow-up with respect to certain ideal containing $p^m$ in the construction. 
In this subsection, we explain how to interpret the notion of 
the category of integrable $p^m$-connections for smooth 
morphisms of schemes $X_n \lra S_n$ 
flat over $\Z/p^n\Z$ and the (level raising) 
inverse functors between them in terms of $\cD$-modules, under the 
assumption that $S_n$ is liftable to a $p$-adic fomal scheme $S$ 
flat over $\Z_p$. (Note that $X_n$ is not liftable to a smooth 
$p$-adic formal scheme over $S$ globally). 
The key point is to consider all the local lifts of $X_n$ to a smooth 
$p$-adic formal scheme over $S$ and consider the `crystalized' category. 

\begin{definition}\label{defc}
Let $S$ be a $p$-adic formal scheme flat over $\Z_p$, 
let $S_n := S \otimes \Z/p^n\Z$ and let $f: X_n \lra S_n$ be a 
smooth morphism of finite type. Then we define the category 
$C(X_n/S)$ as follows$:$ An object is a triple $(U_n,U,i_U)$ 
consisting of an 
open subscheme $U_n$ of $X_n$, a smooth formal scheme $U$ of finite type 
over $X$ 
and a closed immersion $i_U: U_n \hra U$ which makes the following 
diagram Cartesian$:$ 
\[\xymatrix{
U_n \ar[r]^{\subset} \ar[d]^{i_U} & X_n \ar[r]^{f} & S_n \ar[d]^{\cap} \\ 
U \ar[rr] & & S.}\] 
A morphism $\varphi: (U_n,U,i_U) \lra (V_n,V,i_V)$ in $C(X_n/S)$ is defined to 
be a pair of morphism $\varphi_n: U_n \lra V_n$ over $X$ and morphism 
$\varphi: U \lra V$ over $S$ such that the square 
\begin{equation*}
\begin{CD}
U_n @>{i_U}>> U \\ 
@V{\varphi_n}VV @V{\varphi}VV \\ 
V_n @>{i_V}>> V 
\end{CD}
\end{equation*}
is Cartesian. 
\end{definition}

\begin{lemma}\label{taulemma1}
Let $S$ be a $p$-adic formal scheme flat over $\Z_p$ and 
let $f,f': U \lra V$ be morphisms of smooth $p$-adic 
formal schemes of finite type 
over $S$ which coincide modulo $p^n$. Then, for a 
$\cD^{(-m)}_{V/S} \otimes \Z/p^n\Z$-module $\cE$, there exists a 
canonical isomorphism $\tau_{f,f'}: {f'}^*\cE \os{=}{\lra} f^*\cE$ 
of $\cD^{(-m)}_{U/S} \otimes \Z/p^n\Z$-modules. 
\end{lemma}

\begin{proof}
Let $\{\epsilon_k\}_k$ be the $(-m)$-PD-stratification associated to 
$\cE$ and let $f_k^*: {f'}_k^*: \cP^k_{V,(-m)} \lra \cP^k_{U,(-m)}$ 
be the morphism induced by $f,f'$, respectively. \par 
First let us prove that 
$f_k^*$ is equal to ${f'}_k^*$ modulo $p^n$. Since 
$\cP^k_{V,(-m)}$ is topologically generated by $\cO_V$ and 
the elements of the form $(1 \otimes a - a \otimes 1)/p^m \,(a \in \cO_V)$, 
it suffices to check that the images of these elements by 
$f_k^*$ coincides with the image by ${f'}_k^*$. 
For the elements in $\cO_V$, this is clear since $f$ and $f'$ are equal 
modulo $p^n$. Let us consider the images of the element 
$(1 \otimes a - a \otimes 1)/p^m$. If we put 
$f^*(a) - {f'}^*(a) =: p^nb$, we have 
\begin{align*}
& f_k^*((1 \otimes a - a \otimes 1)/p^m) - 
{f'}_k^*((1 \otimes a - a \otimes 1)/p^m) \\ 
= \, & 
((1 \otimes f^*(a) - f^*(a) \otimes 1)/p^m) - 
(1 \otimes {f'}^*(a) - {f'}^*(a) \otimes 1)/p^m) \\ 
= \, & 
p^n(1 \otimes b - b \otimes 1)/p^m. 
\end{align*}
Hence $f_k^*$ is equal to ${f'}_k^*$ modulo $p^n$, as desired. \par 
Let us put 
$\ol{f} = f \,{\rm mod}\,p^n = f'  \,{\rm mod}\,p^n, 
\ol{f}_k^* = f_k^* \,{\rm mod}\,p^n = {f'}_k^*  \,{\rm mod}\,p^n$. 
Then we have the canonical isomorphism 
$$ \tau_{f,f'}: {f'}^*\cE \os{=}{\lra} \ol{f}^*\cE \os{=}{\lra} f^*\cE, $$
and since we have 
${f'}_k^*\epsilon_k = \ol{f}_k^*\epsilon_k = f_k^*\epsilon_k$, 
$\tau_{f,f'}$ gives an isomorphism as 
$\cD^{(-m)}_{U/S} \otimes \Z/p^n\Z$-modules. 
\end{proof}

Using these, we give the following definition. 

\begin{definition}\label{defd}
Let us take $n,n',m \in \N$ with $n \leq n'$, let 
$S$ be a $p$-adic formal scheme flat over $\Z_p$, 
let $S_{n'} := S \otimes \Z/p^{n'}\Z$ and let $f: X_{n'} \lra S_{n'}$ be a 
smooth morphism of finite type. Then we define the category 
$D^{(-m)}(X_{n'}/S)_n$ as the category of pairs 
\begin{equation}\label{objd}
((\cE_U)_{U := (U_{n'},U,i_U) \in C(X_{n'}/S)}, 
(\alpha_{\varphi})_{\varphi:(U_{n'},U,i_U)  \ra (V_{n'},V,i_V) 
\in {\rm Mor}C(X_{n'}/S)}), 
\end{equation}
where $\cE_U$ is a $\cD^{(-m)}_{U/S} \otimes \Z/p^n\Z$-module and 
$\alpha_{\varphi}$ is an isomorphism 
$\varphi^*\cE_{V} \os{=}{\lra} \cE_{U}$ as 
$\cD^{(-m)}_{U/S} \otimes \Z/p^n\Z$-modules satisfying the following 
condition$:$ 
\begin{enumerate}
\item 
$\alpha_{\id} = \id.$
\item 
$\alpha_{\varphi \circ \psi} = \alpha_{\varphi} \circ 
\varphi^*(\alpha_{\psi})$ for morphisms 
$\varphi:(U_{n'},U,i_U) \lra (V_{n'},V,i_V), \psi: 
(V_{n'},V,i_V) \lra (W_{n'},W,i_W)$ in $C(X_{n'}/S)$. 
\item 
For two morphisms $\varphi, \psi: (U_{n'},U,i_U) \lra (V_{n'},V,i_V)$ 
in $C(X_{n'}/S)$, the isomorphism 
$\alpha_{\psi}^{-1} \circ \alpha_{\varphi}: 
\varphi^*\cE_{V} \os{=}{\lra} \cE_U \os{=}{\lra} \psi^*\cE$ 
coincides with $\tau_{\psi,\varphi}$ defined in Lemma \ref{taulemma1}. 
\end{enumerate}
We denote the object \eqref{objd} simply by 
$((\cE_U)_U,(\alpha_{\varphi})_{\varphi})$ or $(\cE_U)_U$. \par 
We call an object $(\cE_U)_U \in D^{(-m)}(X_{n'}/S)_n$ quasi-nilpotent 
if each $\cE_U$ is a quasi-nilpotent $\cD^{(-m)}_{U/S} \otimes 
\Z/p^n\Z$-module, and denote the full subcategory of 
$D^{(-m)}(X_{n'}/S)_n$ consisting of quasi-nilpotent objects by 
$D^{(-m)}(X_{n'}/S)_n^{\qn}$. 
\end{definition}

Then we have the following, which is a $\cD$-module theoretic 
interpretation of the category $\MIC^{(m)}(X_n)$: 

\begin{proposition}\label{eq_1}
Let us take $n,n',m \in \N$ with $n \leq n'$ and let 
$S$ be a $p$-adic formal scheme flat over $\Z_p$.
Let $S_{n'} := S \otimes \Z/p^{n'}\Z$,  
let $f: X_{n'} \lra S_{n'}$ be a 
smooth morphism of finite type and let 
$X_n \lra S_n$ be $f \otimes \Z/p^n\Z$. 
Then there exists the canonical equivalence of categories 
$$ 
\MIC^{(m)}(X_n) \os{=}{\lra} D^{(-m)}(X_{n'}/S)_n, 
\quad 
\MIC^{(m)}(X_n)^{\qn} \os{=}{\lra} D^{(-m)}(X_{n'}/S)^{\qn}_n.
$$
\end{proposition}

\begin{proof}
Assume we are given an object $(\cE,\nabla)$ in 
$\MIC^{(m)}(X_n)$. Let us take an object $(U_{n'},U,i_U)$ of 
$C(X_{n'}/S)$, let us put $U_{n} := U_{n'} \otimes \Z/p^n\Z$ and 
denote the composite $U_n \hra U_{n'} \os{i_U}{\hra} U$ by 
$\ti{i}_U$. Then, via $\ti{i}_U$, we can regard $\cE_U := 
(\cE,\nabla)|_{U_n}$ as an object in $\MIC^{(m)}(U)_n$ by 
Remark \ref{torsrem1}, and by 
Corollary \ref{8cor}, we can regard it as 
a $\cD^{(-m)}_{U/S} \otimes \Z/p^n\Z$-module. Also, for a 
morphism $\varphi: (U_{n'},U,i_U) \lra (V_{n'},V,i_V)$ in 
$C(X_{n'}/S)$, the induced morphism 
$U_n \lra V_n := V \otimes \Z/p^n\Z$ gives the isomorphism 
$\alpha_{\varphi}: 
\varphi^*\cE_V \os{=}{\lra} \cE_U$ as objects in $\MIC^{(m)}(X_n)$ and 
it induces the isomorphim as $\cD^{(-m)}_{U/S} \otimes \Z/p^n\Z$-modules
 by Remarks \ref{torsrem1}, \ref{torsrem2}, Corollary \ref{8cor} and 
Proposition \ref{coin1}. We can check easily that the pair 
$((\cE_U)_U,(\alpha_{\varphi})_{\varphi})$ defines an object in 
$D^{(-m)}(X_{n'}/S)_n$. (The condition (3) in Definition \ref{defd} 
comes from the fact that the isomorphism 
$\alpha_{\varphi}$ defined above depends only on $\varphi$ modulo $p^n$.) \par 
Conversely, if we are given an object 
$((\cE_U)_U,(\alpha_{\varphi})_{\varphi})$ in 
$D^{(-m)}(X_{n'}/S)_n$, each $\cE_U$ naturally defines an object in 
$\MIC^{(m)}(U_{n'})_n$, and the transition maps $\alpha_{\varphi}$'s depend 
only on $\varphi$ modulo $p^{n'}$ up to canonical isomorphism. 
Hence $\cE_U$'s glue to give an object in 
$\MIC^{(m)}(X_{n'})_n = \MIC^{(m)}(X_n)$. We can check that these 
functors give the inverse of each other, and so obtain the first 
equivalence. \par 
To obtain the second equivalence, it suffices to see the consistence of 
the definitions of quasi-nilpotence. This follows from 
Proposition \ref{9}. 
\end{proof}

Next we define the inverse image functor for 
the categories $D^{(-)}(-/-)$. Let $n,n',m \in \N$ with $n \leq n'$ and 
assume given the following 
commutative diagram 
\begin{equation*}
\begin{CD}
X_{n'} @>>> S_{n'} @>{\subseteq}>> S \\ 
@VfVV @VVV @VVV \\ 
Y_{n'} @>>> T_{n'} @>{\subseteq}>> T, 
\end{CD}
\end{equation*}
where $S, T$ are $p$-adic formal scheme flat over $\Z_p$, 
$S_{n'} = S \otimes \Z/p^{n'}\Z, T_{n'} = T \otimes \Z/p^{n'}\Z$, 
the left top arrow and the left bottom arrow are canonical closed 
immersion, the right top arrow and the right bottom arrow are 
smooth. Under this situation, we define the inverse image functor 
$$ f^*: D^{(-m)}(Y_{n'}/T)_n \lra D^{(-m)}(X_{n'}/S)_n $$ 
as follows: Let us take an object $\cE := (\cE_V)_V$ in 
$D^{(-m)}(Y_{n'}/T)_n$ and $(U_{n'},U,i_U) \in C(X_{n'}/S)$. 
Then, locally on $U$, there exists an object 
$(V_{n'},V,i_V) \in C(Y_{n'}/T)$, a morphism 
$\varphi_{n'}:U_{n'} \lra V_{n'}$ over $f$, a morphism 
$\varphi:U \lra V$ over $S \lra T$ with $\varphi \circ i_U = i_V \circ 
\varphi_{n'}$. Then, we define the $\cD^{(-m)}_{U/S} \otimes \Z/p^n\Z$-module 
$(f^*\cE)_U$ by $(f^*\cE)_U := \varphi^*\cE_V$. When there exists 
another object 
$(V'_{n'},V',i_{V'}) \in C(Y_{n'}/T)$ with morphisms 
$\varphi'_{n'}:U_{n'} \lra V'_{n'}$, 
$\varphi' :U \lra V'$ as above, there exists an isomorphism 
$\iota: V \os{=}{\lra} V'$ locally on $V'$. Then, since 
$\iota \circ \varphi$ and $\varphi'$ are equal modulo $p^{n'}$, 
we have the isomorphism 
$$ \varphi^*\cE_V \os{=}{\lra} (\iota \circ \varphi)^*\cE_{V'} 
\os{\tau_{\varphi',\iota \circ \varphi}}{\lra} {\varphi'}^*\cE_{V'}, $$
and this is indendent of the choice of $\iota$ because, when we are given 
another isomorphism $\iota': V \os{=}{\lra} V'$, we have the 
commutative diagram 
\begin{equation*}
\begin{CD}
\varphi^*\cE_V @>{=}>> 
 (\iota \circ \varphi)^*\cE_{V'} 
@>{\tau_{\varphi',\iota \circ \varphi}}>> 
{\varphi'}^*\cE_{V'} \\ 
@\vert @V{\tau_{\iota' \circ \varphi,\iota \circ \varphi}}VV @\vert \\ 
\varphi^*\cE_V @>{=}>> 
(\iota' \circ \varphi)^*\cE_{V'} 
@>{\tau_{\varphi',\iota' \circ \varphi}}>> {\varphi'}^*\cE_{V'}
\end{CD}
\end{equation*}
(the commutativity of the left square is the pull-back of 
the property 
(3) in Definition 
\ref{defd} by $\varphi^*$ 
and that of the right square comes from the definition 
of $\tau_{-/-}$). Therefore, we can glue the local definition 
 $(f^*\cE)_U := \varphi^*\cE_V$ and define the 
$\cD^{(-m)}_{U/S} \otimes \Z/p^n\Z$-module 
$(f^*\cE)_U$ globally. We can also check that the 
$(f^*\cE)_U$'s for $(U_{n'},U,i_U) \in C(X_{n'}/S)$ forms an 
object $f^*\cE := ((f^*\cE)_U)_U$ in $D^{(-m)}(X_{n'}/S)_n$ in the 
same way. By the correspondence 
$\cE \mapsto f^*\cE$, the inverse image 
functor 
$$ f^*: D^{(-m)}(Y_{n'}/T)_n \lra D^{(-m)}(X_{n'}/S)_n $$ 
is defined. Because this functor is defined locally as the 
inverse image functor of $\cD^{(-m)}_{-/-} \otimes \Z/p^n\Z$-modules, 
it induces the functor 
$$ f^{*,\qn}: D^{(-m)}(Y_{n'}/T)^{\qn}_n \lra 
D^{(-m)}(X_{n'}/S)_n^{\qn}. $$ 
Also, we see by the construction 
that the inverse image functors $f^*. f^{*,\qn}$ 
here is equal to the inverse image functors of integrable 
$p^m$-connections 
\begin{align*}
& f_n^*: \MIC^{(m)}(Y_n) \os{=}{\lra} \MIC^{(m)}(X_n), \\
& f_n^{*,\qn}: \MIC^{(m)}(Y_n)^{\qn} \os{=}{\lra} \MIC^{(m)}(X_n)^{\qn}
\end{align*}
(where $X_n := X_{n'} \otimes \Z/p^n\Z, 
Y_n := Y_{n'} \otimes \Z/p^n\Z, f_n := f \otimes \Z/p^n\Z$) 
defined in Section 1 via the equivalences in 
Proposition \ref{eq_1}. \par 

Next we define the level raising inverse image functor. 
First prove the following lemma, which is an analogue of Lemma 
\ref{taulemma1}. 

\begin{lemma}\label{taulemma2}
Assume we are in the situation of $\Hyp(\infty,n+1,n+1)$ with 
$n'=n+1$. 
Let $f,f': U \lra V$ be morphisms of smooth $p$-adic 
formal schemes over $S$ which 
coincide modulo $p^{n+1}$ and put 
$\ol{f} := f \otimes \Z/p^{n+1}\Z = f' \otimes \Z/p^{n+1}\Z: 
U_{n+1} \lra V_{n+1}$. Assume moreover that this morphism fits 
into the following commutative diagram 
\begin{equation*}
\begin{CD}
U_{n+1} @>{\ol{f}}>> V_{n+1} \\ 
@V{\cap}VV @V{\cap}VV \\ 
X_{n+1} @>{F_{n+1}}>> X^{(1)}_{n+1}, 
\end{CD}
\end{equation*}
where the vertical arrows are open immersions. 
$($So $f,f'$ are local lifts of $F_{n+1}$ and so we can define 
the level raising inverse image functor for $f,f'.)$ 
Then, for a 
$\cD^{(-m)}_{V/S} \otimes \Z/p^n\Z$-module $\cE$, the 
canonical isomorphism $\tau_{f,f'}: {f'}^*\cE \os{=}{\lra} f^*\cE$ 
in Lemma \ref{taulemma1} $($which is a priori an isomorphism as 
$\cD^{(-m)}_{U/S} \otimes \Z/p^n\Z$-modules$)$ is an isomorphism 
of $\cD^{(-m+1)}_{U/S} \otimes \Z/p^n\Z$-modules. 
\end{lemma}

\begin{proof}
Since $F_{n+1}$ is a homeomorphism, we may shrink $V$ so that 
$f,f':U \lra V$ are homeomorphisms. Then we can replace $U_{n+1},V_{n+1}$ 
by $X_{n+1}$, $X^{(1)}_{n+1}$ and then by putting 
$X:= U, X^{(1)} := V$, we may assume that the situation is as in 
$\Hyp(\infty,\infty,\infty)$ (but we have two lifts $f,f'$ instead of 
$F$ there.) \par 
Let $\{\epsilon_k\}_k$ be the $(-m)$-PD-stratification associated to 
$\cE$ and let $\Phi_k^*, {\Phi'}_k^*: \allowbreak \cP^k_{X^{(1)},(-m)} 
\allowbreak \lra \cP^k_{X,(-m)}$ 
be the morphism induced by $f,f'$ respectively, by Proposition 
\ref{13}. (See also the diagram \eqref{dd}.) \par 
Let us prove that 
$\Phi_k^*$ is equal to ${\Phi'}_k^*$ modulo $p^n$. 
To prove this, we may work locally. So we assume that 
there exist a local parameter $t_1,...,t_d$ of $X$ and 
let $t'_i, \tau_i, \tau'_i$ be as in the proof of Proposition \ref{13}. 
Note that $f$ and $f'$ coincide modulo $p^{n+1}$. So, by 
the calculation similar to \eqref{calc}, we see that 
${f^2}^*(\tau'_i) - {{f'}^2}^*(\tau'_i)$ can be written as 
an element of the form $p^{n+1}(1 \otimes b_i - b_i \otimes 1)$. 
Then, by definition of $\Phi_k^*, {\Phi'}_k^*$, we have 
$\Phi_k^*(\tau_i/p^m) - {\Phi'}_k^*(\tau_i/p^m) = 
p^{n}\{(1 \otimes b_i - b_i \otimes 1)/p^{m-1}\}$. 
Since $\cP^k_{X^{(1)},(-m)}$ is topologically generated by 
$\cO_{X^{(1)}}$ and $\tau'_i/p^m$'s, we see the coincidence 
of $\Phi_k^*$ and ${\Phi'}_k^*$ modulo $p^n$, as desired. \par 
Let 
us put 
$\ol{\Phi}_k^* = \Phi_k^* \,{\rm mod}\,p^n = {\Phi'}_k^*  \,{\rm mod}\,p^n$. 
Then, 
since we have 
${\Phi'}_k^*\epsilon_k = \ol{\Phi}_k^*\epsilon_k = \Phi_k^*\epsilon_k$, 
the isomorphism $\tau_{f,f'}$ in Lemma \ref{taulemma1} 
gives an isomorphism as 
$\cD^{(-m)}_{X/S} \otimes \Z/p^n\Z$-modules. So we are done. 
\end{proof}

Under $\Hyp(\infty,n',n')$ with $n' \in \N, n' \geq n+1$, 
we define the level raising inverse image functor 
$$ F_{n'}^*: D^{(-m)}(X^{(1)}_{n'}/T)_n \lra D^{(-m+1)}(X_{n'}/S)_n $$ 
as follows: Let us take an object $\cE := (\cE_V)_V$ in 
$D^{(-m)}(X^{(1)}_{n'}/T)_n$ and $(U_{n'},U,i_U) \in C(X_{n'}/S)$. 
Then, locally on $U$, there exists an object 
$(V_{n'},V,i_V) \in C(X^{(1)}_{n'}/T)$, a morphism 
$\varphi_{n'}:U_{n'} \lra V_{n'}$ over $F_{n'}$, a morphism 
$\varphi:U \lra V$ over $S$ with $\varphi \circ i_U = i_V \circ 
\varphi_{n'}$. Then, we define the 
$\cD^{(-m+1)}_{U/S} \otimes \Z/p^n\Z$-module 
$(f^*\cE)_U$ by $(f^*\cE)_U := \varphi^*\cE_V$, 
where the right hand side denotes the level raising inverse image 
by $\varphi$. When there exists 
another object 
$(V'_{n'},V',i_{V'}) \in C(X^{(1)}_{n'}/T)$ with morphisms 
$\varphi'_{n'}:U_{n'} \lra V'_{n'}$, 
$\varphi' :U \lra V'$ as above, there exists an isomorphism 
$\iota: V \os{=}{\lra} V'$ locally on $V'$. Then, since 
$\iota \circ \varphi$ and $\varphi'$ are equal modulo $p^{n'}$, 
we have the isomorphism 
$$ \varphi^*\cE_V \os{=}{\lra} (\iota \circ \varphi)^*\cE_{V'} 
\os{\tau_{\varphi',\iota \circ \varphi}}{\lra} {\varphi'}^*\cE_{V'}, $$
and this is indendent of the choice of $\iota$ as before. 
Therefore, we can glue the local definition 
 $(f^*\cE)_U := \varphi^*\cE_V$ and define the 
$\cD^{(-m)}_{U/S} \otimes \Z/p^n\Z$-module 
$(f^*\cE)_U$ globally. We can also check that the 
$(f^*\cE)_U$'s for $(U_{n'},U,i_U) \in C(X_{n'}/S)$ forms an 
object $f^*\cE := ((f^*\cE)_U)_U$ in $D^{(-m)}(X_{n'}/S)_n$ in the 
same way. By the correspondence 
$\cE \mapsto f^*\cE$, the level raising inverse image 
functor 
$$ F_{n'}^*: D^{(-m)}(X^{(1)}_{n'}/S)_n \lra D^{(-m+1)}(X_{n'}/S)_n $$ 
is defined. Because this functor is defined locally as the 
level raising 
inverse image functor of $\cD^{(-)}_{-/-} \otimes \Z/p^n\Z$-modules, 
it induces the functor 
$$ F_{n'}^{*,\qn}: D^{(-m)}(X^{(1)}_{n'}/T)^{\qn}_n \lra 
D^{(-m+1)}(X_{n'}/S)_n^{\qn}. $$ 
Also, we see by the construction 
that the inverse image functors $F_{n'}^*, F_{n'}^{*,\qn}$ 
here is equal to the level raising inverse image functors of integrable 
$p^m$-connections 
\begin{align*}
& F_{n+1}^*: \MIC^{(m)}(X_n^{(1)}) \lra \MIC^{(m-1)}(X_n), \\ 
& F_{n+1}^{*, \qn}: \MIC^{(m)}(X_n^{(1)})^{\qn} \lra 
\MIC^{(m-1)}(X_n)^{\qn}
\end{align*} 
defined in Section 1 via the equivalences in 
Proposition \ref{eq_1}. In particular, we have given another proof 
of Proposition \ref{functqn2} under the assumption $\Hyp(\infty,n+1,n+1)$, 
as promised in Remark \ref{qnatode}. 

\subsection{Crystalline property of integrable $p^m$-connections}

In this subsection, we prove a crystalline property for the categories 
$\MIC^{(m)}(X)_n, \allowbreak \MIC^{(m)}(X)_n^{\qn}$ for 
a smooth 
morphism $X_n \lra S_n$ satisfying certain liftability condition. 
We also prove a similar result also for the category 
$\MIC^{(m)}(X)$ for a smooth 
morphism $X \lra S$ 
of $p$-adic formal schemes flat over $\Z_p$. 
The key construction in the former case is a variant 
$\ol{D}^{(-m)}(X_{m+e}/S)_n$ (where $e=1$ or $2$) 
of the category of the form 
$D^{(-)}(-/-)_{-}$ and the (level raising) inverse image functor for it. \par 
The starting point is the following, which is an analogue of 
\cite[2.1.5]{berthelotII}: 

\begin{proposition}\label{11}
Let  
\begin{equation*}
\begin{CD}
X @. Y \\ 
@VVV @VVV \\ 
S @>>> T 
\end{CD}
\end{equation*} 
be a diagram of 
$p$-adic formal schemes flat over $\Spf \Z_p$ such that 
the vertical arrows are smooth. Let $m,e \in \N$, let 
$\cE$ be left $\cD_{Y/T}^{(-m)}$-module and assume one of the 
following$:$ 
\begin{enumerate}
\item[\rm (a)] $p \geq 3, e=1$ or $p=2, e=2$. 
\item[\rm (b)] $\cE$ is quasi-nilpotent and $e=1$. 
\end{enumerate}
Let us put $X_{m+e} := X \otimes \Z/p^{m+e}\Z$ 
and let $f_{m+e}: X_{m+e} \lra Y$ be a morphism over $T$. Assume that 
$f,f': X \lra Y$ are morphisms over $T$ which lift $f_{m+e}$. 
Then there exists a canonical $\cD_{X/S}^{(-m)}$-linear 
isomorphism $\ol{\tau}_{f,f'}: {f'}^*\cE \os{=}{\lra} f^*\cE$. \par 
Moreover, when $\cE$ is $p^n$-torsion and 
$f$ and $f'$ are equal modulo $p^{n+m}$ for some $n \geq 1$, the isomorphism 
$\ol{\tau}_{f,f'}$ is equal to the isomorphism $\tau_{f,f'}$ defined in 
Lemma \ref{taulemma1}. 
\end{proposition}

\begin{proof}
By assumption, we have the commutative diagram 
\begin{equation}\label{gal0}
\begin{CD}
X_{m+e} @>{f_{m+e}}>> Y \\ 
@V{\cap}VV @VV{\text{diag.}}V \\ 
X @>{(f,f')}>> Y \times_T Y, 
\end{CD}
\end{equation}
and by the universality of formal blow-up and PD-envelope, 
the lower horizontal arrow factors as 
\begin{equation}\label{gal1}
X \lra T_{Y,(-m-e)} \lra T_{Y,(-m)} \lra Y \times_T Y.  
\end{equation}
Let us denote the morphism $X \lra T_{Y,(-m)}$ in the diagram 
\eqref{gal1} (the composite of the first two arrows) by $g'$ and 
let us denote the composite $X_e := X \otimes \Z/p^e\Z \hra 
X_{m+e} \os{f_{m+e}}{\lra} Y$ by $f_e$. 
Then the diagram 
\begin{equation}\label{gal2}
\begin{CD}
X_{e} @>{f_{e}}>> Y \\ 
@V{\cap}VV @VVV \\ 
X @>{g'}>> T_{Y,(-m)}
\end{CD}
\end{equation}
(where the right vertical arrow is 
the morphism induced by the diagonal map) 
is commutative: Indeed, using the fact that 
$\cO_{T_{Y,(-m)}}$ is locally topologically generated 
by sections $a$ with $p^ma \in \cO_{Y \times_T Y}$, the commutativity 
of the diagram \eqref{gal0} induces that of \eqref{gal2}. 
Then, by the universality of PD-envelope, the diagram \eqref{gal2}
 induces the commutative diagram 
\begin{equation}\label{b}
\begin{CD}
X_e @>{f_e}>> Y \\ 
@V{\cap}VV @VV{\cap}V \\ 
X @>g>> P_{Y,(-m)}. 
\end{CD}
\end{equation}
In the case (a), the defining ideal $p^e\cO_{X}$ of the 
closed immersion $X_e \hra X$ is topologically PD-nilpotent. So 
the morphism $X \lra P_{Y,(-m)}$ in \eqref{b} factors 
as $X \os{g_r}{\lra} P^r_{Y,(-m)} \hra  P_{Y,(-m)}$ for some $r \in \N$. 
Then we define $\ol{\tau}_{f,f'} := g_r^*(\eta_r)$, 
where $\eta_r$ is the isomorphism of $(-m)$-PD-stratification associated 
to $\cE$ on $P^r_{Y,(-m)}$. 
In the case (b), 
we define $\ol{\tau}_{f,f'} := g^*(\eta)$, 
where $\eta$ is the isomorphism of $(-m)$-HPD-stratification associated 
to $\cE$ on $P_{Y,(-m)}$. Then we see that $\ol{\tau}_{f,f'}$ is an 
$\cO_X$-linear isomorphism. \par 
So it suffices to prove that $\ol{\tau}_{f,f'}$ is $\cD_{X/S}^{(-m)}$-linear. 
Let us put $\cF := f^*\cE, \cF' := {f'}^*\cE$, let us take $l \in \N$, 
let $p_i: P^l_{X,(-m)} \lra X \,(i=0,1)$ be the 
morphism induced by the $i$-th projection $X^2 \lra X$ and 
let $\epsilon_l: p_1^*\cF \os{=}{\lra} p_0^*\cF, 
\epsilon'_l: p_1^*\cF' \os{=}{\lra} p_0^*\cF'$ be 
the isomorphism of $(-m)$-PD-stratification for $\cF, \cF'$ on 
$P^l_{X,(-m)}$. Then it suffice to prove the commutativity 
of the following diagram of sheaves on $P^l_{X,(-m)}$: 
\begin{equation}\label{11kome}
\begin{CD}
p_1^*\cF' @>{p_1^*(\ol{\tau}_{f',f})}>> p_1^*\cF \\ 
@V{\epsilon'_l}VV @V{\epsilon_l}VV \\ 
p_0^*\cF' @>{p_0^*(\ol{\tau}_{f',f})}>> p_0^*\cF. 
\end{CD}
\end{equation}
Let us consider the following commutative diagram
\begin{equation}\label{gal4}
\begin{CD}
X_{m+e} @>{f_{m+e}}>> Y \\ 
@VVV @VVV \\ 
X^2 @>{(f \times f, f' \times f')}>> Y^4, 
\end{CD}
\end{equation}
where the vertical arrows are the diagonal embeddings. 
Then, by the universality of formal blow-up, 
the composite $P_{X,(-m)} \lra X^2 \os{(f \times f, f' \times f')}{\lra} 
Y^4$ factors as $P_{X,(-m)} \os{h'}{\lra} T_{Y,(-m)}(3) \lra Y^4$, and we see 
(in the same way as the proof of the commutativity of \eqref{gal2}) that 
the commutativity of the diagram \eqref{gal4} induces that of 
the following diagram: 
\begin{equation}\label{gal5}
\begin{CD}
X_{e} @>{f_{e}}>> Y \\ 
@VVV @VVV \\ 
P_{X,(-m)} @>{h'}>> T_{Y,(-m)}(3). 
\end{CD}
\end{equation}
Then, 
noting that the defining ideal of 
the closed immersion $X_e \hra X \hra P_{X,(-m)}$ admits 
a PD-structure canonically, we see that the 
diagram \eqref{gal5} gives rise to the morphism 
$h: P^l_{X,(-m)} \lra P_{Y,(-m)}(3)$. Also, 
the defining ideal of the closed immersion 
$X_e \hra P^l_{X,(-m)}$ still admits 
a PD-structure canonically. \par 
Assume that we are in the case (a). Then the PD-structure 
on the defining ideal of $X_e \hra P^l_{X,(-m)}$ is 
topologically PD-nilpotent. So the morphism $h$ factors as 
$P^l_{X,(-m)} \os{h_s}{\lra} P^s_{Y,(-m)}(3) \lra P_{Y,(-m)}(3)$ 
for some $s\in \N, s \geq l, r$. Let 
$q_{ij}: P^s_{Y,(-m)}(3) \lra P^s_{Y,(-m)}$ \,$(0 \leq i < j \leq 3)$ 
be the morphism induced by the $(i,j)$-th projection 
$X^4 \lra X^2$. Then we have 
$p_0^*(\ol{\tau}_{f',f})=p_0^*g_s^*(\eta_s) = h_s^*q_{02}^*(\eta_s)$ and 
$p_1^*(\ol{\tau}_{f',f})= h_s^*q_{13}^*(\eta_s)$. Also, when we denote 
the morphism $P^l_{X,(-m)} \lra P^s_{Y,(-m)}$ induced by 
$f \times f: X^2 \lra Y^2$ by $\varphi$, 
we have 
$\epsilon_l = \varphi^*(\eta_s) = h^*q_{01}^*(\eta_s)$ and 
similarly we have $\epsilon'_l = h^*q_{23}(\eta_s)$. So 
The commutativity of the diagram \eqref{11kome} follows from 
the cocycle condition for $\eta_s$. So we are done in the 
case (a). In the case (b), we can prove the commutativity of 
the diagram \eqref{11kome} in the same way, by replacing 
$\eta_s$ by $\eta$ and $P^s_{Y,(-m)}$ by $P_{Y,(-m)}$. \par 
Finally, let us assume that 
$\cE$ is $p^n$-torsion and that 
$f$ and $f'$ are equal modulo $p^{n+m}$ for some $n \in \N$. 
Let us consider locally and take a local coordinate 
$t_1, ..., t_d$ of $Y$ over $T$. Then, 
by definition, 
the isomorphism $\tau_{f,f'}: {f'}^*\cE \os{=}{\lra} f^*\cE$ is written as 
$${f'}^*(e) \mapsto \sum_{k\in\N^d} (({f'}^*(t)-f^*(t))/p^m)^{[k]} 
f^*(\pa^{\lara{k}}e),$$ 
and the $k$-th term on the right hand side is contained in 
$(k!)^{-1}p^{|k|n} f^*\cE \subseteq p^nf^*\cE = 0$ when 
$k\not=0$. Hence 
we have $\ol{\tau}_{f,f'}({f'}^*(e)) = f^*(e)$ and this implies 
the equality $\ol{\tau}_{f,f'} = \tau_{f,f'}$. 
\end{proof}

\begin{remark}
When $\cE$ is $p^n$-torsion and $f,f'$ are equal only modulo $p^n$, 
the isomorphism $\ol{\tau}_{f,f'}$ above is not necessarily equal to 
the isomorphism $\tau_{f,f'}$ in Lemma \ref{taulemma1} unless $m=0$. 
\end{remark}

By the argument in \cite[2.1.6]{berthelotII}, we have the following 
immediate corollary of Proposition \ref{11} (we omit the proof): 

\begin{corollary}\label{formalcase1}
\begin{enumerate}
\item 
Let us put $e=1$ if $p \geq 3$ and $e=2$ if $p=2$. 
let $f: X \lra S$ be a 
smooth morphism of formal schemes flat over $\Z_p$ and let 
$X_{m+e} \lra S_{m+e}$ be $f \otimes \Z/p^{m+e}\Z$. 
Then the category $$\MIC^{(m)}(X) = 
\text{\rm (left $\cD_{X/S}^{(-m)}$-modules)}$$ 
depends only on the diagram $X_{m+e} \lra S_{m+e} \hra S$ and 
functorial with respect to this diagram. 
\item 
let $f: X \lra S$ be a 
smooth morphism of formal schemes flat over $\Z_p$ and let 
$X_{m+1} \lra S_{m+1}$ be $f \otimes \Z/p^{m+1}\Z$. 
Then the category $$\MIC^{(m)}(X)^{\qn} = 
\text{\rm (quasi-nilpotent left $\cD_{X/S}^{(-m)}$-modules)}$$ 
depends only on the diagram $X_{m+1} \lra S_{m+1} \hra S$ and 
functorial with respect to this diagram. 
\end{enumerate}
\end{corollary}

Next we consider the case of $p^n$-torsion objects. First, 
using Proposition \ref{11}, we give the following definition: 

\begin{definition}\label{defdb}
\begin{enumerate}
\item 
Let us put $e=1$ if $p \geq 3$ and $e=2$ if $p=2$. 
Let us take $n,n',m \in \N$ with $m+e \leq n'$, let 
$S$ be a $p$-adic formal scheme flat over $\Z_p$, 
let $S_{n'} := S \otimes \Z/p^{n'}\Z$ and let $f: X_{n'} \lra S_{n'}$ be a 
smooth morphism of finite type. Then we define the category 
$\ol{D}^{(-m)}(X_{n'}/S)_n$ as the category of pairs 
\begin{equation}\label{objdb}
((\cE_U)_{U := (U_{n'},U,i_U) \in C(X_{n'}/S)}, 
(\alpha_{\varphi})_{\varphi:(U_{n'},U,i_U)  \ra (V_{n'},V,i_V) 
\in {\rm Mor}C(X_{n'}/S)}), 
\end{equation}
where $\cE_U$ is a $\cD^{(-m)}_{U/S} \otimes \Z/p^n\Z$-module and 
$\alpha_{\varphi}$ is an isomorphism 
$\varphi^*\cE_{V} \os{=}{\lra} \cE_{U}$ as 
$\cD^{(-m)}_{U/S} \otimes \Z/p^n\Z$-modules satisfying the following 
conditions$:$ 
\begin{enumerate}
\item 
$\alpha_{\id} = \id.$
\item 
$\alpha_{\varphi \circ \psi} = \alpha_{\varphi} \circ 
\varphi^*(\alpha_{\psi})$ for morphisms 
$\varphi:(U_{n'},U,i_U) \lra (V_{n'},V,i_V), \psi: 
(V_{n'},V,i_V) \lra (W_{n'},W,i_W)$ in $C(X_{n'}/S)$. 
\item 
For two morphisms $\varphi, \psi: (U_{n'},U,i_U) \lra (V_{n'},V,i_V)$ 
in $C(X_{n'}/S)$, the isomorphism 
$\alpha_{\psi}^{-1} \circ \alpha_{\varphi}: 
\varphi^*\cE_{V} \os{=}{\lra} \cE_U \os{=}{\lra} \psi^*\cE$ 
coincides with $\ol{\tau}_{\psi,\varphi}$ defined in 
Proposition \ref{11}. 
\end{enumerate}
We denote the object \eqref{objdb} simply by 
$((\cE_U)_U,(\alpha_{\varphi})_{\varphi})$ or $(\cE_U)_U$. 
\item 
Let us put $e=1$ and 
let us take $n,n',m \in \N$, 
$S, f: X_{n'} \lra S_{n'}$ as in $(1)$. Then we define the category 
$\ol{D}^{(-m)}(X_{n'}/S)^{\qn}_n$ as the category of pairs 
\eqref{objdb}
where $\cE_U$ is a quasi-nilpotent 
$\cD^{(-m)}_{U/S} \otimes \Z/p^n\Z$-module and 
$\alpha_{\varphi}$ is an isomorphism 
$\varphi^*\cE_{V} \os{=}{\lra} \cE_{U}$ as 
$\cD^{(-m)}_{U/S} \otimes \Z/p^n\Z$-modules satisfying the 
conditions $(a), (b), (c)$ in $(1)$. 
\end{enumerate}
\end{definition}

When $n' \geq \max(m+n, m+e)$ (where $e$ is as in Definition \ref{defdb}),  
we have the equalities 
$\ol{D}^{(m)}(X_{n'}/S)_{n} = D^{(m)}(X_{n'}/S)_{n}$, 
$\ol{D}^{(m)}(X_{n'}/S)_{n}^{\qn} = D^{(m)}(X_{n'}/S)_{n}^{\qn}$, 
because 
the isomorphisms $\ol{\tau}_{\psi,\varphi}$ used in 
Definition \ref{defdb} are equal to 
the isomorphisms $\tau_{\psi,\varphi}$ used in 
Definition \ref{defd}. Now let us note that, for 
$n', n'' \in \N$ with $m+e \leq n' \leq n''$, 
the functor 
\begin{align*}
r': & \,\, C(X_{n''}/S) \lra C(X_{n'}/S); \\ & (U_{n''},U,i_U) \mapsto 
(U_{n'}:= U_{n''} \otimes \Z/p^{n'}\Z, U, U_{n'} \hra U_{n''} 
\os{i_U}{\hra} U)
\end{align*}
induces the functor 
$r: \ol{D}^{(m)}(X_{n'}/S)_n \lra \ol{D}^{(m)}(X_{n''}/S)_n$, 
$r: \ol{D}^{(m)}(X_{n'}/S)^{\qn}_n \lra \ol{D}^{(m)}(X_{n''}/S)^{\qn}_n$. 
Hence we obtain the functor 
\begin{align}
& R: \ol{D}^{(m)}(X_{m+e}/S)_n \os{r}{\lra} 
\ol{D}^{(m)}(X_{n'}/S)_{n} \os{=}{\lra}  
D^{(m)}(X_{n'}/S)_{n}, \label{R} \\
& R: \ol{D}^{(m)}(X_{m+e}/S)^{\qn}_n \os{r}{\lra} 
\ol{D}^{(m)}(X_{n'}/S)^{\qn}_{n} \os{=}{\lra}  
D^{(m)}(X_{n'}/S)^{\qn}_{n}, \label{Rqn}. 
\end{align}
for $n' \geq \max(m+n,m+e)$. 
Then we have the following: 

\begin{proposition}\label{eq_1b}
\begin{enumerate}
\item 
Let us put $e=1$ if $p \geq 3$ and $e=2$ if $p=2$. 
Let us take $n,n',m \in \N$ with $n' \geq \max(m+n,m+e)$ and let 
$S$ be a $p$-adic formal scheme flat over $\Z_p$. 
Let $S_{n'} := S \otimes \Z/p^{n'}\Z$,  
let $f: X_{n'} \lra S_{n'}$ be a 
smooth morphism of finite type and let 
$X_{m+e} \lra S_{m+e}$ be $f \otimes \Z/p^{m+e}\Z$. 
Then the functor \eqref{R} is an equivalence of categories. 
\item 
Let us put $e=1$ and let the other notations be as in $(1)$. 
Then 
the functor \eqref{Rqn} is an equivalence of categories. 
\end{enumerate}
\end{proposition}

\begin{proof}
Since the proof is the same, we only prove (1). To do so, 
it suffices to construct the inverse of the functor 
$r: \ol{D}^{(m)}(X_{m+e}/S)_n \lra \ol{D}^{(m)}(X_{n'}/S)_{n}$. 
So let us take $\cE := (\cE_V)_V \in \ol{D}^{(m)}(X_{n'}/S)_{n}$ 
and take an object $(U_{m+e},U,i_U)$ in $C(X_{m+e}/S)$. Then, 
locally on $U$, there exists an object 
$(V_{n'},V,i_V)$ in $C(X_{n'}/S)$ and a morphism 
$\varphi: U \lra V$ over $S$ 
inducing $\ol{\varphi}: U_{m+e} \lra V_{n'}$ which is a morphism 
over the canonical closed immersion $X_{m+e} \hra X_{n'}$: 
Indeed, if we put $U_{n'} := U \otimes \Z/p^{n'}\Z$, it is 
a smooth lift of $U_{m+e} \hra X_{m+e} \lra S_{m+e}$ over 
$S_{n'}$. Hence we have the isomorphism $U_{n'} \isom V_{n'}$ between 
$U_{n'}$ and some open subscheme $V_{n'}$ of $X_{n'}$ locally on $U$. 
Then, locally on $U$, $V_{n'}$ 
admits a smooth lift $i_V: V_{n'} \hra V$ over $S$ and 
the isomorphism  $U_{n'} \isom V_{n'}$ is liftable to an isomorphism 
$\varphi: U \lra V$ over $S$, as desired. \par 
Taking $(V_{n'},V,i_V)$ and $\varphi: U \lra V$ as in the previuos 
paragraph, we define a $\cD^{(-m)}_{U/S} \otimes \Z/p^n\Z$-module 
$r^{-1}(\cE)_U$ by $r^{-1}(\cE)_U := \varphi^*\cE_V$. 
When there exists another object $(V'_{n'},V',i_{V'})$ and 
another isomorphism $\varphi': U \lra V'$, there exists an isomorphism 
$\iota: V \os{=}{\lra} V'$ locally on $V'$. 
Then we have the isomorphism 
$$ \varphi^*\cE_V \os{=}{\lra} (\iota \circ \varphi)^*\cE_{V'} 
\os{\ol{\tau}_{\varphi',\iota \circ \varphi}}{\lra} {\varphi'}^*\cE_{V'}, $$
and this is indendent of the choice of $\iota$ because, when we are given 
another isomorphism $\iota': V \os{=}{\lra} V'$, we have the 
commutative diagram 
\begin{equation*}
\begin{CD}
\varphi^*\cE_V @>{=}>> 
 (\iota \circ \varphi)^*\cE_{V'} 
@>{\ol{\tau}_{\varphi',\iota \circ \varphi}}>> 
{\varphi'}^*\cE_{V'} \\ 
@\vert @V{\ol{\tau}_{\iota' \circ \varphi,\iota \circ \varphi}}VV @\vert \\ 
\varphi^*\cE_V @>{=}>> 
(\iota' \circ \varphi)^*\cE_{V'} 
@>{\ol{\tau}_{\varphi',\iota' \circ \varphi}}>> {\varphi'}^*\cE_{V'}. 
\end{CD}
\end{equation*}
Therefore, we can glue the local definition 
 $r^{-1}(\cE)_U := \varphi^*\cE_V$ and define the 
$\cD^{(-m)}_{U/S} \otimes \Z/p^n\Z$-module 
$r^{-1}(\cE)_U$ globally. We can also check that the 
$r^{-1}(\cE)_U$'s for $(U_{m+e},U, \allowbreak 
i_U) \in C(X_{m+e}/S)$ forms an 
object $r^{-1}\cE := ((r^{-1}\cE)_U)_U$ in $\ol{D}^{(-m)}(X_{m+e}/S)_n$ 
in the same way, and so we have defined the functor 
$r^{-1}: \ol{D}^{(m)}(X_{n'}/S)_n \lra \ol{D}^{(m)}(X_{m+e}/S)_{n}$, 
which is easily seen to be the inverse of the functor $r$. So we are done. 
\end{proof}

Let us put $e=1$ if $p \geq 3$ and $e=2$ if $p=2$. 
Let us take $n,m \in \N$ and assume given the following 
commutative diagram 
\begin{equation}\label{diadiag}
\begin{CD}
X_{m+e} @>>> S_{m+e} @>{\subseteq}>> S \\ 
@VfVV @VVV @VVV \\ 
Y_{m+e} @>>> T_{m+e} @>{\subseteq}>> T, 
\end{CD}
\end{equation}
where $S, T$ are $p$-adic formal scheme flat over $\Z_p$, 
$S_{m+e} = S \otimes \Z/p^{m+e}\Z, T_{m+e} = T \otimes \Z/p^{m+e}\Z$, 
the left top arrow and the left bottom arrow are canonical closed 
immersions, the right top arrow and the right bottom arrow are 
smooth. Under this situation, we can define the inverse image functor 
\begin{equation}\label{fstardb}
f^*: \ol{D}^{(-m)}(Y_{m+e}/T)_n \lra \ol{D}^{(-m)}(X_{m+e}/S)_n 
\end{equation}
in the same way as the inverse image functor 
\begin{equation}\label{fstard}
f^*: D^{(-m)}(Y_{n'}/T)_n \lra D^{(-m)}(X_{n'}/S)_n \quad (n' \geq n)
\end{equation}
defined before. Also, when $n' \geq \max(m+n,m+e)$ and 
the diagram \eqref{diadiag} is liftable to the diagram 
\begin{equation}\label{diadiag2}
\begin{CD}
X_{n'} @>>> S_{n'} @>{\subseteq}>> S \\ 
@VfVV @VVV @VVV \\ 
Y_{n'} @>>> T_{n'} @>{\subseteq}>> T 
\end{CD}
\end{equation}
(where 
$S_{n'} = S \otimes \Z/p^{n'}\Z, T_{n'} = T \otimes \Z/p^{n'}\Z$ 
and the right top arrow and the right bottom arrow are 
smooth), we have the equality $R \circ \text{\eqref{fstardb}} = 
\text{\eqref{fstard}} \circ R$. Also, when $e=1$, we have the 
inverse image functor 
\begin{equation}\label{fstarqndb}
f^{*,\qn}: \ol{D}^{(-m)}(Y_{m+e}/T)^{\qn}_n 
\lra \ol{D}^{(-m)}(X_{m+e}/S)^{\qn}_n 
\end{equation}
in the same way as the inverse image functor 
\begin{equation}\label{fstarqnd}
f^{*,\qn}: D^{(-m)}(Y_{n'}/T)^{\qn}_n 
\lra D^{(-m)}(X_{n'}/S)^{\qn}_n \quad (n' \geq n)
\end{equation}
defined before, and when $n' \geq \max(m+n,m+e)$ and 
the diagram \eqref{diadiag} is liftable to the diagram 
\eqref{diadiag2}, $R \circ \text{\eqref{fstarqndb}} = 
\text{\eqref{fstarqnd}} \circ R$. 
Hence we have the following 
corollary of Proposition \ref{eq_1b}, which is the 
main result in this subsection. 

\begin{corollary}\label{axcor}
\begin{enumerate}
\item 
Let the notations be as in Proposition \ref{eq_1b}$(1)$. 
Then the category $\MIC^{(m)}(X_n) = D^{(-m)}(X_{n'}/S)_n$ 
$($where $X_n := X_{n'} \otimes \Z/p^n\Z)$ 
depends only on the diagram $X_{m+e} \lra S_{m+e} \hra S$ and 
functorial with respect to this diagram. 
\item 
Let the notations be as in Proposition \ref{eq_1b}$(2)$. 
Then the category $\MIC^{(m)}(X_n)^{\qn} \allowbreak = 
D^{(-m)}(X_{n'}/S)^{\qn}_n$ 
$($where $X_n := X_{n'} \otimes \Z/p^n\Z)$ 
depends only on the diagram $X_{m+1} \lra S_{m+1} \hra S$ and 
functorial with respect to this diagram. 
\end{enumerate}
\end{corollary}

\begin{remark}
Let $S, n'$ be as above and put $S_j := S \otimes \Z/p^j\Z$ for 
$j \in \N$. 
Note that the above corollary does not imply that the category 
$\MIC^{(m)}(X_n)$ depends only on $X_{m+e} \lra S_{m+e} \hra S$
for any smooth morphism $X_n \lra S_n$: The above corollary is 
applicable only for the smooth morphism $X_n \lra S_n$ which is 
liftable to a smooth morphism $X_{n'} \lra S_{n'}$. 
\end{remark}

Next we discuss the crystalline property of the level raising 
inverse image functor. To do so, 
we need the following proposition. 

\begin{proposition}\label{taubar2}
Let the notations be as in $\Hyp(\infty,\infty,\infty)$. 
Let $m,e \in \N, \geq 1$, 
let $\cE$ be a left $\cD^{(-m)}_{X^{(1)}/S}$-module 
and assume one of the following$:$ 
\begin{enumerate}
\item[\rm (a)] $p \geq 3, e=1$ or $p=2, e=2$. 
\item[\rm (b)] $\cE$ is quasi-nilpotent and $e=1$. 
\end{enumerate}
Suppose that 
we have another morphism $F': X \lra X^{(1)}$ over $S$ lifting 
the morphism $F_{X_1/S_1}$ which coincides with $F$ modulo $p^{m+e}$. 
Then the 
isomorphism $\ol{\tau}_{F,F'}: {F'}^*\cE \os{=}{\lra} F^*\cE$ defined in 
Proposition \ref{11} is actually $\cD_{X/S}^{(-m+1)}$-linear. 
\end{proposition}

\begin{proof}
Let $t_i, t'_i, \tau_i, \tau'_i$ be as in the proof of 
Proposition \ref{13}. Then we can write 
$F^*(t'_i) = t_i^p + pa_i, {F'}^*(t'_i) = t_i^p + pa_i + p^{m+e}b_i$ 
for some $a_i,b_i \in \cO_X$, and we see by the same calcuation 
as in the proof of Proposition \ref{13} that 
there exist elements $\sigma_i 
\in I := \Ker(\cO_{X^2} \ra \cO_X), \sigma'_i \in \cO_{X^2}$ 
such that 
$(F \times F')^*(\tau'_i) = \tau_i^p + p\sigma_i + p^{m+e}\sigma'_i$. 
So, for $m \geq 2$, the image of this element in 
$\cO_{T_{X,(-m+1)}}$ belongs to $p^m\cO_{T_{X,(-m+1)}}$ and 
in the case $m=1$, the image of this element in 
$\cO_{P_{X,(0)}}$ belongs to $p\cO_{P_{X,(0)}}$. 
Therefore, in both cases, the image of this element in 
$\cO_{P_{X,(-m+1)}}$ belongs to $p^m\cO_{P_{X,(-m+1)}}$. \par 
Now let us consider the morphism $h':=(F \times F, F' \times F'): 
X^2 \lra (X^{(1)})^4$. Then we have the commutative diagram 
\begin{equation}\label{ang1}
\begin{CD}
X_{m+e} @>{F_{m+e}}>> X^{(1)} \\ 
@VVV @VVV \\ 
X^2 @>{h'}>> (X^{(1)})^4, 
\end{CD}
\end{equation}
where $F_{m+e}$ is the composite 
$X_{m+e} \hra X \os{F}{\lra} X^{(1)}$, which is also written as 
the composite $X_{m+e} \hra X \os{F'}{\lra} X^{(1)}$. 
Let us denote the $q$-th projection $(X^{(1)})^{4} \lra X^{(1)}$ by 
$\pi_q$ $\,(0 \leq q \leq 3)$, $(q,q+1)$-th projection 
$(X^{(1)})^{4} \lra (X^{(1)})^2$ by $\pi_{q,q+1}$ $\,(0 \leq q \leq 2)$ 
and put $\tau'_{i,q} := \pi_{q+1}^*t'_i - \pi_q^*t'_i = 
\pi_{q,q+1}^*(\tau'_i)$. Then $\Ker (\cO_{(X^{(1)})^4} \lra 
\cO_{X^{(1)}})$ is generated by $\tau'_{i,q}$'s $(1 \leq i \leq d, 0 \leq q 
\leq 3)$. If we denote the $i$-th projection $X^2 \lra X$ by $p_i \,(i=0,1)$, 
we have 
${h'}^*(\tau_{i,0}) = {h'}^*(\pi_{1}^*t'_i - \pi_0^*t'_i) = 
p_0^*F^*t'_i - p_0^*F^*t'_i =0$ and by similar reason, we also have 
${h'}^*(\tau_{i,2})=0$. Also, we have 
$$ {h'}^*(\tau_{i,1}) = {h'}^*(\pi_{2}^*t'_i - \pi_1^*t'_i)
= p_1^*{F'}^*t'_i - p_0^*F^*t'_i = (F \times F')^*(\tau'_i) $$ 
and the image of this element in belongs in 
$\cO_{P_{X,(-m+1)}}$ belongs to $p^m\cO_{P_{X,(-m+1)}}$. Hence, 
by the universality of formal blow-up, the morphism 
$P_{X,(-m+1)} \lra X^2 \os{h'}{\lra} (X^{(1)})^4$ 
factors as 
$$ P_{X,(-m+1)} \os{h''}{\lra} T_{X^{(1)},(-m)}(3) \lra (X^{(1)})^4. $$
Furthermore, since 
$\cO_{T_{X^{(1)},(-m)}(3)}$ is locally topologically generated by 
the elements in $\cO_{X^{(1)}}$ and 
$\tau'_{i,q}/p^m$, the commutative diagram \eqref{ang1} induces the 
commutative diagram
\begin{equation*}
\begin{CD}
X_e @>{F_e}>> X^{(1)} \\ 
@VVV @VVV \\ 
P_{X,(-m+1)} @>{h''}>> T_{X^{(1)},(-m)}(3), 
\end{CD}
\end{equation*}
where $X_e = X \otimes \Z/p^e\Z$ and $F_e$ is the composite 
$X_e \hra X_{m+e} \os{F_{m+e}}{\lra} X^{(1)}$. 
Noting that the defining ideal of 
the closed immersion $X_e \hra P_{X,(-m+1)}$ admits 
a PD-structure canonically, we see that the 
above diagram gives rise to the morphism 
$h: P^l_{X,(-m+1)} \lra P_{Y,(-m)}(3)$ for any $l \in \N$. 
Then, in the case (b), we can prove the commutativity of the 
diagram \eqref{11kome} on $P^l_{X,(-m+1)}$ by using the 
morphism $h$, in the same way as the proof of Proposition \ref{11}. 
In the case (a), we see that the morphism $h$ factors as 
$P^l_{X,(-m+1)} \os{h_s}{\lra} P^s_{Y,(-m)}(3) \lra P_{Y,(-m)}(3)$ 
for some $s \in \N, s \geq l$ because the the PD-structure on the ideal 
of the defining ideal of $X_e \hra P_{X,(-m+1)}$ is 
topologically PD-nilpotent, and 
then we can prove the commutativity of the 
diagram \eqref{11kome} on $P^l_{X,(-m+1)}$ by using $h_s$. 
So we are done. 
\end{proof}

Again by the argument in \cite[2.1.6]{berthelotII} 
(see also \cite[2.2.6]{berthelotII}), we have the following 
immediate corollary of Proposition \ref{taubar2} (we omit the proof): 

\begin{corollary}\label{formalcase2}
Let the notations be as in $\Hyp(\infty,\infty,\infty)$. 
\begin{enumerate}
\item 
Let us put $e=1$ if $p \geq 3$ and $e=2$ if $p=2$. 
Then the level raising inverse image functor 
$$ F^*: \text{\rm (left $\cD_{X^{(1)}/S}^{(-m)}$-modules)} \lra 
\text{\rm (left $\cD_{X/S}^{(-m+1)}$-modules)}, $$
which is equal to the level raising inverse image functor 
$F^*: \MIC^{(m)}(X^{(1)}) \lra \MIC^{(-m+1)}(X)$, depends only on 
$F_{m+e} := F \,{\rm mod}\, p^{m+e}$. 
\item 
The level raising inverse image functor 
\begin{align*}
F^{*,\qn}: & \,\, \text{\rm (quasi-nilpotent 
left $\cD_{X^{(1)}/S}^{(-m)}$-modules)} \\ & 
\hspace{1cm} \lra 
\text{\rm (quasi-nilpitent left $\cD_{X/S}^{(-m+1)}$-modules)}, 
\end{align*}
which is equal to the level raising inverse image functor 
$F^{*,\qn}: \MIC^{(m)}(X^{(1)})^{\qn} \allowbreak \lra 
\MIC^{(-m+1)}(X)^{\qn}$, depends only on 
$F_{m+1} := F \,{\rm mod}\, p^{m+1}$. 
\end{enumerate}
\end{corollary}

Next we consider the case of $p^n$-torsion objects. 
Let $m \in \N, \geq 1$, 
assume that we are in the situation of $\Hyp(\infty,m+e,m+e)$ 
with $e=1$ if $p \geq 3$ and $e=2$ if 
$p =2$, and let us take $n \in \N$. Then we can define the level raising 
inverse image functor 
\begin{equation}\label{Fstardb}
F_{m+e}^*: \ol{D}^{(-m)}(X^{(1)}_{m+e}/T)_n \lra \ol{D}^{(-m+1)}(X_{m+e}/S)_n 
\end{equation}
in the same way as the level raising inverse image functor 
\begin{equation}\label{Fstard}
F_{n'}^*: D^{(-m)}(X^{(1)}_{n'}/T)_n 
\lra D^{(-m+1)}(X_{n'}/S)_n \quad (n' \geq n+1)
\end{equation}
defined before. Also, when $n' \geq \max(m+n,m+e)$ and 
when we are in the situation of $\Hyp(\infty,n',n')$, 
we have the equality $R \circ \text{\eqref{Fstardb}} = 
\text{\eqref{Fstard}} \circ R$. Also, when $e=1$, we have the 
inverse image functor 
\begin{equation}\label{Fstarqndb}
F^{*,\qn}_{m+e}: \ol{D}^{(-m)}(X^{(1)}_{m+e}/T)^{\qn}_n 
\lra \ol{D}^{(-m+1)}(X_{m+e}/S)^{\qn}_n 
\end{equation}
in the same way as the inverse image functor 
\begin{equation}\label{Fstarqnd}
F^{*,\qn}_{n'}: D^{(-m)}(X^{(1)}_{n'}/T)^{\qn}_n 
\lra D^{(-m)}(X_{n'}/S)^{\qn}_n \quad (n' \geq n+1)
\end{equation}
defined before, and when $n' \geq \max(m+n,m+e)$ and 
when we are in the situation of $\Hyp(\infty,n',n')$, 
we have the equality 
$R \circ \text{\eqref{Fstarqndb}} = 
\text{\eqref{Fstarqnd}} \circ R$. 
Hence we have the following 
corollary, which is the second main result in this subsection. 

\begin{corollary}
\begin{enumerate}
\item 
Let $m \geq 1$ and let us put $e=1$ if $p \geq 3$, $e=2$ if $p=2$. 
Then, under $\Hyp(\infty,n',n')$ with $n' \geq \max(m+n,m+e)$, 
the level raising inverse image functor 
$$ F_{n'}^*: D^{(-m)}(X^{(1)}_{n'}/S)_n \lra D^{(-m+1)}(X_{n'}/S)_n $$
$($which is equal to the level raising inverse image functor 
$F_{n+1}^*: \MIC^{(m)}(X^{(1)}_{n}/S) \allowbreak 
\lra \MIC^{(m-1)}(X_{n}/S))$ 
depends only on $F_{m+e} = F_{n'} \,{\rm mod}\, p^{m+e}$. 
\item 
Let $m \geq 1$.  
Then, under $\Hyp(\infty,n',n')$ with $n' \geq m+n$, 
the level raising inverse image functor 
$$ F_{n'}^{*,\qn}: D^{(-m)}(X^{(1)}_{n'}/S)^{\qn}_n \lra 
D^{(-m+1)}(X_{n'}/S)^{\qn}_n $$
$($which is equal to the level raising inverse image functor 
$F_{n+1}^{*,\qn}: \MIC^{(m)}(X^{(1)}_{n}/S)^{\qn} 
\allowbreak \lra \MIC^{(m-1)}(X_{n}/S)^{\qn})$ 
depends only on $F_{m+1} = F_{n'} \,{\rm mod}\, p^{m+1}$. 
\end{enumerate}
\end{corollary}

\begin{remark}
Let $m, e, n'$ be as above. 
Note that the above corollary does not imply that, under the 
situation $\Hyp(\infty,n+1,n+1)$, the level raising inverse image functor 
\begin{align*}
& F_{n+1}^*: \MIC^{(m)}(X^{(1)}_{n}/S) \lra \MIC^{(m-1)}(X_{n}/S)) \\ 
& \text{(resp. $F_{n+1}^{*,\qn}: \MIC^{(m)}(X^{(1)}_{n}/S)^{\qn} 
\lra \MIC^{(m-1)}(X_{n}/S)^{\qn})$} 
\end{align*} 
depends only on $F_{m+e}$ (resp. $F_{m+1}$): 
The above corollary is 
applicable only in the situation $\Hyp(\infty,n',n')$. 
\end{remark}

\section{Frobenius descent to the level minus one} 

In this section, we prove that the 
level raising inverse image functor for 
relative Frobenius gives an equivalence between 
the category of quasi-nilpotent integrable $p$-connections and 
the category of quasi-nilpotent integrable connections. 
In terms of $\cD$-modules, this is an equivalence of the category 
of quasi-nilpotent left $\cD$-modules of level $-1$ and the 
the category of quasi-nilpotent left $\cD$-modules of level $0$. 
So we can say this result as `the Frobenius descent to the level $-1$'. 
The method of the proof is similar to the proof of 
Frobenius descent due to Berthelot \cite{berthelotII}. \par 
The main result in this section is the following: 

\begin{theorem}[Frobenius descent to the level minus one]\label{fd}
Assume that we are in the situation of $\Hyp(\infty,\infty,\infty)$. 
Then the level raising inverse image functor 
$$F^*: \MIC^{(1)}(X^{(1)})^{\qn} \lra \MIC(X)^{\qn}$$ 
is an equivalence of categories. 
\end{theorem}

We have the following immediate corollaries: 

\begin{corollary}\label{fdcor}
\begin{enumerate}
\item 
Assume that we are in the situation of $\Hyp(\infty,n+1,n+1)$. 
Then the level raising inverse image functors 
\begin{align*}
& F_{n+1}^*: \MIC^{(1)}(X^{(1)}_n)^{\qn} \lra \MIC(X_n)^{\qn} \\ 
& F_{n+1}^*: D^{(1)}(X^{(1)}_{n+1})^{\qn}_n \lra D^{(0)}(X_{n+1})^{\qn}_n 
\end{align*}
are equivalences of categories. 
\item 
Assume that we are in the situation of $\Hyp(\infty,2,2)$. 
Then, for $n \in \N$, the level raising 
inverse image functor 
$$F_2^*: \ol{D}^{(1)}(X^{(1)}_{2})^{\qn}_n \lra \ol{D}^{(0)}(X_{2})^{\qn}_n $$ 
is an equivalence of categories. 
\end{enumerate}
\end{corollary}

\begin{proof}
Since all the categories appearing in the statement satisfy the 
descent property for the Zariski topology, we may work 
Zariski locally. Then we can assume that we are in 
the situation of 
$\Hyp(\infty,\infty,\infty)$, 
and in this case, the level raising inverse image functors 
are interpreted as the $p^n$-torsion part of 
the level raising inverse image functor 
$$F^*: \MIC^{(1)}(X^{(1)})^{\qn} \lra \MIC(X)^{\qn}$$ 
in Theorem \ref{fd}. So the corollary follows from 
Theorem \ref{fd}. 
\end{proof}

Note that this gives a possible answer to Question 
\ref{q}. We prove several lemmas to prove Theorem \ref{fd}. 

\begin{lemma}\label{17}
Let the notations be as in $\Hyp(\infty,\infty,\infty)$. Then 
the morphism $\Phi: P_{X/S,(0)} \lra P_{X^{(1)}/S,(-1)}$ 
defined in Proposition \ref{13} 
is a finite flat morphism of degree $p^{2d}$. 
\end{lemma}

\begin{proof}
It suffices to prove that the morphism $P_{X/S,(0)} 
\lra X \times_{X^{(1)}} P_{X^{(1)}/S,(-1)}$ induced by $\Phi$ 
is a finite flat morphism of degree $p^d$. To show this, we 
may work locally. So we can take a local 
coordinate $t_1, ..., t_d$ of $X$ over $S$. Let us put 
$t'_i := 1 \otimes t_i \in \cO_{X^{(1)}}$, $\tau_i := 
1 \otimes t_i - t_i \otimes 1 \in \cO_{X^2}, 
\tau'_i := 1 \otimes t'_i - t'_i \otimes 1 \in \cO_{(X^{(1)})^2}$. 
Also, let us put $F^*(t'_i) = t_i^p + pa_i$. 
The homomorphism of sheaves corresponding to the morphism 
$P_{X/S,(0)} 
\lra X \times_{X^{(1)}} P_{X^{(1)}/S,(-1)}$ 
has the form 
\begin{equation}\label{18diag1}
\cO_X \lara{\tau'_i/p}_{1 \leq i \leq d} \lra 
\cO_X \lara{\tau_i}_{1 \leq i \leq d}. 
\end{equation}
Since the morphism $F^2: X^2 \lra (X^{(1)})^2$ sends $\tau'_i$ to 
\begin{align*} 
1 \otimes (t_i^p+pa_i) - (t_i^p+pa_i) \otimes 1 
& = 
(\tau_i + t_i \otimes 1)^p - 
t_i^p \otimes 1 + p(1 \otimes a_i - a_i \otimes 1) \\ 
& = 
\tau_i^p + \sum_{k=1}^{p-1} \begin{pmatrix} p \\ k \end{pmatrix} 
t_i^{p-k}\tau_i^k + p(1 \otimes a_i - a_i \otimes 1) 
\end{align*}
and the morphism $\Phi$ is induced by $F^2$, it follows that 
$\tau'_i/p$ is sent by the morphism \eqref{18diag1} to the 
element 
$\sum_{k=1}^{p-1} p^{-1} \begin{pmatrix} p \\ k \end{pmatrix} 
t_i^{p-k} \tau_i^k + \tau_i^{[p]} + 
\sum_{k \in \N^d, k \not= 0} \pa^{\lara{k}_0}(a_i) \tau^{[k]}$. 
For $l \in \N$, let us put 
$I_l := \{k = (k_i)_i \in \N^d \,|\, k \not= 0, \forall i, k_i < p^{l+1}\}$. 
Then, since we have $\pa^{\lara{k}_0}(a_i) \in k!\cO_X \subseteq p\cO_X$ 
for $k \in \N^d \setminus (I_0 \cup \{0\})$, 
we see that 
$\tau'_i/p$ is sent by the morphism \eqref{18diag1} to the 
element of the form 
$\tau_i^{[p]} + 
\sum_{k \in I_0} u_{i,k} \tau^k + pv_{i}$ for some $u_{i,k} \in \cO_X, 
v_i \in \cO_X \lara{\tau_i}_{1 \leq i \leq d}$. 
Hence, for $l \in \N$, 
$(\tau'_i/p)^{[p^l]}$ is sent by the morphism \eqref{18diag1} to 
an element of the form 
$\tau_i^{[p^{l+1}]} + \sum_{k \in I_l} u_{i,l,k} \tau^{[k]} + pv_{i,l}$ 
for some 
$u_{i,l,k} \in \cO_X, v_{i,l} \in 
\cO_X \lara{\tau_i}_{1 \leq i \leq d}$. \par 
To prove the lemma, we may assume that $X$ is affine and it suffices 
to prove that 
the morphism \eqref{18diag1} modulo $p$ is finite flat of defree 
$p^d$. Let us put $A := \cO_X/p\cO_X$. 
Then the morphism \eqref{18diag1} modulo $p$ has the form 
\begin{equation}\label{18diag2}
A[x_{i,l}]_{1 \leq i \leq d, l \geq 0}/(x_{i,l}^p)_{i,l} 
\lra 
A[y_{i,l}]_{1 \leq i \leq d, l \geq 0}/(y_{i,l}^p)_{i,l} 
\end{equation}
($x_{i,l}$ corresponds to 
the element $(\tau'_i/p)^{[p^l]} \,{\rm mod}\, p$ 
and $y_{i,l}$ corresponds to 
the element $\tau^{[p^l]} \,{\rm mod}\, p$) 
and $x_{i,l}$ is sent to an element of the form 
$y_{i,l+1}+a_{i,l}$, where 
$a_{i,l}$ is an element in 
$A[y_{i,l'}]_{1 \leq i \leq d, 0 \leq l' \leq l}/(y_{i,l'}^p)_{i,l'}$ 
with no degree $0$ part. (Here the degree is taken 
with respect to $y_{i,k}$'s.) Let us denote the degree $1$ part of 
$a_{i,l}$ by $b_{i,l}$. Let 
us consider the $A$-algebra homomorphism 
$$ \alpha: A[z_{i,l}]_{1 \leq i \leq d, l \geq -1}/(z_{i,l}^p)_{i,l} \lra 
A[y_{i,l}]_{1 \leq i \leq d, l \geq 0}/(y_{i,l}^p)_{i,l} 
$$ 
defined by $\alpha(z_{i,-1})=y_{i,0}, \alpha(z_{i,l})=
y_{i,l+1} + b_{i,l} \,(l \geq 0)$. Let us define 
the $A$-algebra homomorphism 
$$ \beta: A[y_{i,l}]_{1 \leq i \leq d, l \geq 0}/(y_{i,l}^p)_{i,l} 
\lra A[z_{i,l}]_{1 \leq i \leq d, l \geq -1}/(z_{i,l}^p)_{i,l}$$ 
of the converse direction inductively, in the following way: 
First, let us define $\beta(y_{i,0}) := z_{i,-1}$. 
When we defined $\beta (y_{i,l'})$ for $0 \leq l' \leq l$, 
we can define $\beta (b_{i,l})$ since $b_{i,l}$ is a linear form 
in $y_{i,l'}$'s for $1 \leq i \leq d, 0 \leq l' \leq l$. 
Then we define $\beta(y_{i,l+1}):=
z_{i,l}- \beta(b_{i,l})$. Then $\beta$ is well-defined, and 
it is easy to see that $\alpha$ and $\beta$ are the inverse of 
each other. Hence $\beta$ is an isomorphism, and so 
it suffices to prove that the composite $\beta \circ 
\text{\eqref{18diag2}}$ is finite flat of degree $p^d$. 
Notice that we can factorize 
$\beta \circ 
\text{\eqref{18diag2}}$ as 
\begin{align}
A[x_{i,l}]_{1 \leq i \leq d, l \geq 0}/(x_{i,l}^p)_{i,l} 
& \os{\subset}{\lra} 
A[x_{i,l}]_{1 \leq i \leq d, l \geq -1}/(x_{i,l}^p)_{i,l} 
\label{18diag3} \\ & 
\lra A[z_{i,l}]_{1 \leq i \leq d, l \geq -1}/(z_{i,l}^p)_{i,l} 
\nonumber 
\end{align}
by introducing new variables $x_{i,-1} \,(1 \leq i \leq d)$ and by 
sending them to 
$z_{i,-1}\,(1 \leq i \leq d)$. Also, 
noting the fact that $\beta$ sends $y_{i,l}$ to 
linear forms in $z_{i,l'}$'s ($-1 \leq l' \leq l$), 
we see that $x_{i,l}\,(l \geq 0)$ is sent by 
$\beta \circ 
\text{\eqref{18diag2}}$ (thus by the second homomorphism of 
\eqref{18diag3}) to an element of the form 
$z_{i,l} + c_{i,l}$, where 
$c_{i,l}$ is an element in 
$A[z_{i,l'}]_{1 \leq i \leq d, -1 \leq l' \leq l}/(z_{i,l'}^p)_{i,l'}$ 
whose degree $0$ part and degree $1$ part are zero. 
From this expression, we see easily that the second homomorphism 
in \eqref{18diag3} is an isomorphism. On the other hand, 
it is clear that the first homomorphism in \eqref{18diag3} is 
finite flat of degree $p^d$. So we see that $\beta \circ 
\text{\eqref{18diag2}}$ is finite flat of degree $p^d$ and so the 
proof of the lemma is finished. 
\end{proof}

\begin{lemma}\label{18}
Let the notations be as in $\Hyp(\infty,\infty,\infty)$. Then 
the defining ideal of the 
diagonal morphism $X \hra X\times_{X^{(1)}} X$ 
admits a unique PD-structure. So the closed immersion 
$X \times_{X^{(1)}} X \hra X \times_S X$ induces the 
PD-morphism $X \times_{X^{(1)}} X \lra P_{X,(0)}$. 
\end{lemma}

\begin{proof}
We have the equality 
$\cO_{X \times_{X^{(1)}}X} = \cO_{X \times_S X} / J$ 
(where $J$ is the ideal topologically generated by the 
elements 
$1 \otimes F^*(y) - F^*(y) \otimes 1 \,(y \in \cO_{X^{(1)}})$), 
and the kernel of $\cO_{X \times_{X^{(1)}} X} \lra \cO_X$ is 
topologically generated by the elements $1 \otimes x - x 
\otimes 1 \,(x \in \cO_{X})$. 
For $x \in \cO_X$, let us put $x' := 1 \otimes x \in \cO_{X^{(1)}}$ 
and put $F^*(x') = x^p + pz$. Then we have 
\begin{align*}
0 & = 1 \otimes F^*(x') - F^*(x') \otimes 1 \\ 
& = 1 \otimes x^p - x^p \otimes 1 + p(1 \otimes z - z \otimes 1) \\ 
& = ((1 \otimes x - x \otimes 1) + x \otimes 1)^p - x^p \otimes 1 
+  p(1 \otimes z - z \otimes 1) \\ 
& = 
(1 \otimes x - x \otimes 1)^p + 
\sum_{i=1}^{p-1}\begin{pmatrix} p \\ i \end{pmatrix} 
(x^{p-i} \otimes 1)(1 \otimes x - x \otimes 1)^i + 
p(1 \otimes z - z \otimes 1) 
\end{align*}
in $\cO_{X \times_{X^{(1)}}X}$. So, in $\cO_{X \times_{X^{(1)}}X}$, 
$(1 \otimes x - x \otimes 1)^p$ has the form 
\begin{equation}\label{nocd}
pa(1 \otimes x - x \otimes 1) + p(1 \otimes z - z \otimes 1)
\end{equation}
for some $a \in \cO_{X \times_{X^{(1)}}X}$ and $z \in \cO_X$. \par 
We should prove that the ideal $I$ admits a unique PD-structure. 
Since $\cO_{X \times_{X^{(1)}}X}$ 
is a flat $\Z_p$-algebra, it suffices to 
prove that, for any $x \in \cO_X$ and $k \in \N$, 
$(1 \otimes x - x \otimes 1)^k \in k!\cO_{X \times_{X^{(1)}}X}$. 
First we prove it for $k = p^l \,(l \in \N)$, by induction on $l$. 
In the case $l=0$, it is trivially true. In general, we have 
\begin{align*}
(1 \otimes x - x \otimes 1)^{p^l} & = 
((1 \otimes x -x \otimes 1)^p)^{p^{l-1}} \\ 
& = (pa(1 \otimes x - x \otimes 1) + p(1 \otimes z - z \otimes 1))^{p^{l-1}} 
\qquad \text{(by \eqref{nocd})} \\ 
& = 
\sum_{i=0}^{p^{l-1}} p^{p^{l-1}}
\begin{pmatrix} p^{l-1} \\ i \end{pmatrix} a^{p^{l-1}-i}
(1 \otimes x - x \otimes 1)^{p^{l-1}-i}(1 \otimes z - z \otimes 1)^{i}, 
\end{align*}
and by induction hypothesis, the $i$-th term is contained in 
$$ 
p^{p^{l-1}}\begin{pmatrix} p^{l-1} \\ i \end{pmatrix} 
(p^{l-1}-i)!i!\cO_{X \times_{X^{(1)}}X} = 
p^{p^{l-1}}(p^{l-1})!\cO_{X \times_{X^{(1)}}X} 
= p^l! \cO_{X \times_{X^{(1)}}X}. $$
So we have $(1 \otimes x - x \otimes 1)^{p^l} \in 
p^l! \cO_{X \times_{X^{(1)}}X}$, as desired. 
For $k \in \N$ with $k \not= p^l \,(\forall l \in \N)$, 
let us take the maximal integer $l'$ such that $k$ is divisible by 
$p^{l'}$. Then we have 
$(1 \otimes x - x \otimes 1)^k = (1 \otimes x - x \otimes 1)^{p^{l'}}
(1 \otimes x - x \otimes 1)^{k-p^{l'}}$, and it is contained in 
$p^{l'}!(k-p^{l'})! \cO_{X \times_{X^{(1)}}X} = k! 
\cO_{X \times_{X^{(1)}}X}$ by the result for $p^{l'}$ and the induction 
hypothesis. So we are done. 
\end{proof}

\begin{remark}\label{18rem}
By the same argument, we see the following: 
For $d \in \N$, the kernel $I_1$ of the homomorphism 
$\cO_{X \times_{X^{(1)}} X}\langle y_1,...,y_d \rangle \lra 
\cO_{X} \langle y_1,...,y_d \rangle$ 
(induced by the diagonal morphism) 
admits a unique PD-structure. Moreover, if we put 
\begin{align*}
& I_2 := \Ker (\cO_{X \times_{X^{(1)}} X}\langle y_1,...,y_d \rangle 
\lra \cO_{X \times_{X^{(1)}}X}), \\ 
& I := I_1 + I_2 = \Ker (
\cO_{X \times_{X^{(1)}} X}\langle y_1,...,y_d \rangle \lra \cO_X), 
\end{align*}
we see that $I$ admits a 
unique PD-structure compatible with that on $I_1$ and that on $I_2$, 
since, for $x \in I_1 \cap I_2$, we have 
$n!x^{[n]} = x^n$ for both PD-structures. 
\end{remark} 

\begin{remark}\label{18rem2}
If we put $X(r)$ to be the $(r+1)$-fold fiber product of 
$X$ over $X^{(1)}$, we see in the same way as the proof given above 
that the defining ideal of the 
diagonal morphism $X \hra X(r)$ 
admits a unique PD-structure. So the closed immersion 
$X(r) \hra X^{r+1}$ induces the 
PD-morphism $X(r) \lra P_{X,(0)}(r)$. 
\end{remark}

\begin{lemma}\label{19}
Let the situation be as in $\Hyp(\infty,\infty,\infty)$. 
Then the defining ideal of the diagonal closed immersion 
$X \hra P_{X,(0)} \times_{P_{X^{(1)},(-1)}} P_{X,(0)}$ 
admits a unique PD-structure. 
\end{lemma}

\begin{proof}
The uniqueness follows from the fact that 
$P_{X,(0)} \times_{P_{X^{(1)},(-1)}} P_{X,(0)}$, being flat over 
$P_{X,(0)}$ by Lemma \ref{17}, is flat over $\Z_p$. Let us prove the 
existence of the desired PD-structure. Note that we have 
the isomorphisms 
\begin{equation}\label{19diag1}
\cO_{P_{X,(0)} \times_{P_{X^{(1)},(-1)}} P_{X,(0)}} 
\cong 
\cO_X \lara{\tau_i}_i \otimes_{\cO_{X^{(1)}}\lara{\tau'_i/p}_i} 
\cO_X \lara{\tau_i}_i \cong 
(\cO_{X} \otimes_{\cO_{X^{(1)}}} \cO_X) \lara{\tau_{i,0},\tau_{i,1}}_i /I 
\end{equation}
(where $\tau_{i,0} := \tau_i \otimes 1, \tau_{i,1} := 1 \otimes \tau_i$), 
where $I$ is the closed ideal topologically generated by 
\begin{equation}\label{19diag2}
\Phi_1^*(\tau'_i/p)^{[k]} 
- 
\Phi_0^*(\tau'_i/p)^{[k]} \quad 
(1 \leq i \leq d, k \in \N), 
\end{equation}
where $\Phi_j^* \,(j=0,1)$ is the homomorphism 
\begin{align*}
\cO_{X^{(1)}}\lara{\tau'_i/p}_i = 
\cP_{X^{(1)},(-1)} & \os{\Phi^*}{\lra} \cP_{X,(0)} \\ 
& \os{\text{$j$-th incl.}}{\lra} 
\cP_{X,(0)} \otimes_{\cO_X} \cP_{X,(0)} = 
(\cO_{X} \otimes_{\cO_{X^{(1)}}} \cO_X) \lara{\tau_{i,0},\tau_{i,1}}_i. 
\end{align*}
First, let us note that, 
by Remark \ref{18rem}, the ideal 
\begin{align*}
& \Ker(
(\cO_{X} \otimes_{\cO_{X^{(1)}}} \cO_X) \lara{\tau_{i,0},\tau_{i,1}}_i 
\lra \cO_X) \\ 
= \, & 
\Ker((\cO_{X} \otimes_{\cO_{X^{(1)}}} \cO_X) 
\lara{\tau_{i,0},\tau_{i,1}}_i \lra \cO_{X} \otimes_{\cO_{X^{(1)}}} \cO_X) \\ 
& \hspace{1cm} + 
\Ker((\cO_{X} \otimes_{\cO_{X^{(1)}}} \cO_X) 
\lara{\tau_{i,0},\tau_{i,1}}_i \lra \cO_X \lara{\tau_{i,0},\tau_{i,1}}_i) 
\end{align*}
admits the PD-structure compatible with the canonical 
PD-structure on the ideal 
$\Ker((\cO_{X} \otimes_{\cO_{X^{(1)}}} \cO_X) 
\lara{\tau_{i,0},\tau_{i,1}}_i \lra \cO_{X} \otimes_{\cO_{X^{(1)}}} \cO_X)$. 
So it suffices to 
prove that the ideal $I$ is a PD-subideal of 
$\Ker(
(\cO_{X} \otimes_{\cO_{X^{(1)}}} \cO_X) \lara{\tau_{i,0},\tau_{i,1}}_i 
\lra \cO_X)$ to prove the lemma. If we denote the element 
\eqref{19diag2} simply by $a - b$, it suffices to prove 
that $(a - b)^{[l]} \in I$ for any $l \geq 1$. If we take the maximal 
integer $m$ such that $l$ is divisible by $p^m$, we have 
$(a - b)^{[l]} = (a-b)^{[l-p^m]}(a-b)^{[p^m]}$. So it suffices to 
prove that $(a-b)^{[p^m]} \in I$ for $m \geq 0$. 
To prove it, it suffices to prove the following claim: 
For $m \geq 0$, there exist elements 
$c_j \in \cO_X  \lara{\tau_{i,0},\tau_{i,1}}_i\,(0 \leq j \leq m)$ 
with $(a-b)^{[p^m]} = \sum_{j=0}^{m}c_j(a^{[p^j]}-b^{[p^j]})$, 
because $a^{[p^j]}-b^{[p^j]} \in I$. We prove this claim by 
induction on $m$. When $m=0$, the claim is trivially true. 
Assume that the claim is true for $m-1$ and put 
$(a-b)^{[p^{m-1}]} = \sum_{j=0}^{m-1}c_j(a^{[p^j]}-b^{[p^j]})$. 
Then we have 
\begin{align*}
(a-b)^{[p^{m}]} & = ((a-b)^{[p^{m-1}]})^{[p]} 
= (\sum_{j=0}^{m-1}c_j(a^{[p^j]}-b^{[p^j]}))^{[p]} \\ 
& = \sum_{j=0}^{m-1}c_j^p(a^{[p^j]}-b^{[p^j]})^{[p]} + A 
\end{align*}
for some $A$ of the form $\sum_{j=0}^{m-1}d_j(a^{[p^j]}-b^{[p^j]})$. 
Moreover, we have 
$$ (a^{[p^j]}-b^{[p^j]})^{[p]} = (a^{[p^{j+1}]} - b^{[p^{j+1}]}) 
+ \sum_{s=1}^{p-1} \dfrac{1}{s!(p-s)!}
(a^{[p^j]})^{p-s}(-b^{[p^j]})^s $$ 
in the case $p \geq 3$, and it is easy to see that the second term 
on the right hand side is a multiple of $a^{[p^j]} - b^{[p^j]}$. 
Hence the proof is finished in the case $p \geq 3$. 
The case $p=2$ follows from the equality 
$$ (a^{[2^j]}-b^{[2^j]})^{[2]} = 
a^{[2^{j+1}]} - a^{[2^j]}b^{[2^j]} + b^{[2^{j+1}]} 
= 
a^{[2^{j+1}]} - b^{[2^{j+1}]} - b^{[2^j]}(a^{[2^j]} - b^{[2^j]}). 
$$ 
So we are done. 
\end{proof}

Now we are ready to prove Theorem \ref{fd}. The proof is 
similar to that of \cite[2.3.6]{berthelotII}. 

\begin{proof}[Proof of Theorem \ref{fd}]
In the proof, we freely regard an object in $\MIC^{(1)}(X^{(1)})^{\qn}$ 
(resp. $\MIC(X)^{\qn}$) as a quasi-nilpotent left 
$\cD^{(-1)}_{X^{(1)}/S}$-module (resp. $\cD^{(0)}_{X/S}$-module) or a 
$p$-power torsion module with $(-1)$-HPD-stratification on $X^{(1)}$ 
(resp. $0$-HPD-stratification on $X$). \par 
Since $F:X \lra X^{(1)}$ is finite flat, the functor $F^*$ is 
faithful. Let us prove that $F^*$ is full. 
Let $\Phi: P_{X,(0)} {\lra} P_{X^{(1)},(-1)}$ be the morphism defined in 
Proposition \ref{13}, let $u: X \times_{X^{(1)}} X {\lra} P_{X,(0)}$ be
 the morphism defined in Lemma \ref{18} and let 
$\ol{p}_j:X \times_{X^{(1)}} X \lra X\,(j=0,1)$ be the $j$-th projection. 
Let us take an object $(\cE',\epsilon') \in \MIC^{(1)}(X^{(1)})^{\qn}$ 
and let us put $(\cE,\epsilon) := F^*(\cE',\epsilon') 
= (F^*\cE', \Phi^*\epsilon') \in \MIC(X)^{\qn}$. Then 
$u^*\epsilon = u^*\Phi^*\epsilon'$ is an isomorphism 
$\ol{p}_1^*\cE \os{=}{\lra} \ol{p}_0^*\cE$, and by using Remarks 
\ref{13rem} and \ref{18rem2}, we see that it satisfies the 
cocycle condition on $X \times_{X^{(1)}} X \times_{X^{(1)}} X$. 
So $(\cE,u^*\epsilon)$ is a descent data on $X$ relative to $X^{(1)}$. 
If we take a local coordinate $t_1, ..., t_d$ of $X$ over $S$ and 
if we put $t'_i := 1 \otimes t_i \in \cO_{X^{(1)}}, 
\tau_i := 1 \otimes t_i - t_i \otimes 1 \in \cO_{X^2}, 
\tau'_i := 1 \otimes t'_i - t'_i \otimes 1 \in 
\cO_{(X^{(1)})^2}, F^*(t'_i) = t_i^p + pa_i$, 
the PD-homomorphism of sheaves 
\begin{align*}
(\Phi \circ u)^*: 
\cO_{P_{X^{(1)},(-1)}} \lra \cO_{P_{X,(0)}} \lra 
\cO_{X \times_{X^{(1)}} X}
\end{align*}
associated to $\Phi \circ u$ sends 
$\tau'_i/p$ as 
\begin{align*}
\tau'_i/p & \mapsto 
\tau_i^{[p]}+p^{-1}\sum_{j=1}^{p-1}\begin{pmatrix} p \\ j 
\end{pmatrix} t_i^{p-j}\tau_i^j + (1 \otimes a_i - a_i \otimes 1) \\ 
& \mapsto 
-p^{-1}\sum_{j=1}^{p-1}\begin{pmatrix} p \\ j 
\end{pmatrix} t_i^{p-j}\tau_i^j - (1 \otimes a_i - a_i \otimes 1) \\ 
& \hspace{1cm} + p^{-1}\sum_{j=1}^{p-1}\begin{pmatrix} p \\ j 
\end{pmatrix} t_i^{p-j}\tau_i^j + (1 \otimes a_i - a_i \otimes 1) 
=0, 
\end{align*}
we see that $\Phi \circ u$ factors through $X^{(1)}$. Hence 
$u^*\epsilon$ is the pull-back of the identity map on $\cE'$ by 
$X \times_{X^{(1)}} X \lra X^{(1)}$, that is, 
the descent data $(\cE,u^*\epsilon)$ is equal to the 
one coming canonically from $\cE'$. \par 
Now let us take 
$(\cE',\epsilon'), (\cF',\eta') \in \MIC^{(1)}(X^{(1)}_n)^{\qn}$, 
put $(\cE,\epsilon) := F^*(\cE',\epsilon'), 
(\cF,\eta) \allowbreak 
:= F^*(\cF',\eta) \in \MIC(X_n)^{\qn}$ and assume that we are given 
a morphism $\varphi: (\cE,\epsilon) \lra (\cF,\eta)$. 
For an open subscheme $U$ of $X$, let us put 
$U^{(1)} := U \times_{X} X^{(1)}$, 
$A := \Gamma(U,\cO_X), A' := \Gamma(U^{(1)}, \allowbreak \cO_{X^{(1)}}), 
E := \Gamma(U,\cE), E' := \Gamma(U^{(1)}, \cE'), 
F := \Gamma(U,\cF), F' := \Gamma(U^{(1)}, \cF')$. 
Then we have $A \otimes_{A'} E' = E, A \otimes_{A'} F' = F$ and 
by the argument in the previous paragraph, $E, F$ are naturally 
endowed with the descent data relative to $A'$ coming canonically 
from $E',F'$ and the morphism 
$\Gamma(U,\varphi): E \lra F$ is a morphism of descent data. 
Hence it descents to a morphism $\psi_U: E' \lra F'$. 
By letting $U$ vary, we see that $\{\psi_U\}_U$ defines 
a morphism $\psi:\cE' \lra \cF'$ with $F^*(\psi) = \varphi$. 
To prove that $\psi$ induces a morphism 
$(\cE',\epsilon') \lra (\cF',\eta')$, we should prove the 
compatibility of $\psi$ with $\epsilon', \eta'$. Since 
$\Phi$ is finite flat, it suffices to prove the compatibility of 
$F^*\psi = \varphi$ with $\Phi^*\epsilon' = \epsilon, 
\Phi^*\eta' = \eta$ and it follows from definition. So 
$\psi$ is a morphism in $\MIC^{(1)}(X^{(1)})^{\qn}$ with 
$F^*\psi = \varphi$ and so the functor $F^*$ is full, as desired. \par 
We prove that the functor $F^*$ is essentially surjective. 
Let us take $(\cE,\epsilon) \in \MIC(X)^{\qn}$. Then, as we saw above, 
$u^*\epsilon$ defines a descent data on $\cE$ relative to $X^{(1)}$. 
Hence, for any open subscheme $U$ of $X$ and $A, A', E$ as above, 
$u^*\epsilon$ defines a descent data on the $A$-module $E$ 
relative to $A'$. Hence it descents to a $A'$-module $E'$ satisfying 
$A \otimes_{A'} E' = E$, since $A$ is finite flat over $A'$. 
Next, let $U$ be an open subscheme of $X$, $U = \bigcup_i U_i$ be 
an open covering and put $U_{ij} := U_i \cap U_j$. 
Let $E_i, E_{ij}$ 
(resp. $E'_i, E'_{ij}$) be the module $E$ (resp. $E'$) 
in the case $U=U_i, U=U_{ij}$ respectively. Then we have the 
exact sequence 
$$ 0 \lra E \lra \prod_i E_i \lra \prod_{i,j}E_{ij}, $$
and it implies the exactness of the sequence 
$$ 0 \lra E' \lra \prod_i E'_i \lra \prod_{i,j}E'_{ij}. $$
Hence, by letting $U$ vary, $E'$'s induce a sheaf of 
$\cO_{X^{(1)}}$-module $\cE'$ with $F^*\cE' = \cE$. \par 
Let $p_j: P_{X,(0)} \lra X, p'_j: P_{X^{(1)},(-1)} \lra X^{(1)} \,(j=0,1)$ 
be the morphisms induced by $j$-th projection. Then 
$\epsilon: p_1^*\cE \lra p_0^*\cE$ is rewritten as 
$\epsilon: \Phi^*{p'}_1^*\cE' \lra \Phi^*{p'}_0^*\cE'$. 
We prove that $\epsilon$ descents to a morphism 
$\epsilon': {p'}_1^*\cE' \lra {p'}_0^*\cE'$. Let 
$\pi_{ij}: P_{X,(0)} \times_{P_{X^{(1)},(-1)}} P_{X,(0)} \lra 
 P_{X,(0)}\,((i,j)=(0,1), (2,3))$ be the morphism induced by the 
$(i,j)$-th projection $X^4 \lra X^2$ and let 
$\rho_j\,(j=0,1)$ be the descent data on $p_j^*\cE$ coming from 
${p'}_j^*\cE'$. Then, to see the existence of $\epsilon'$, it suffices 
to prove the commutativity of the following diagram, since 
$\Phi$ is finite flat: 
\begin{equation}\label{21diag}
\begin{CD}
\pi_{23}^*(p_1^*\cE) @>{\rho_1}>> \pi_{01}^*(p_1^*\cE) \\ 
@V{\pi_{23}^*\epsilon}VV @V{\pi_{01}^*\epsilon}VV \\ 
\pi_{23}^*(p_0^*\cE) @>{\rho_0}>> \pi_{01}^*(p_0^*\cE). 
\end{CD}
\end{equation}
Note that, by Lemma \ref{19}, the morphism 
$\ti{P} := P_{X,(0)} \times_{P_{X^{(1)},(-1)}} P_{X,(0)} \lra 
X^4$ induced by two $P_{X,(0)} \lra X^2$'s induces the morphism 
$v: \ti{P} \lra P_{X,(0)}(3)$. Let 
$q_{ij}: P_{X,(0)}(3) \lra 
P_{X,(0)}\,(0 \leq i < j \leq 3)$ be the morphism induced by 
the $(i,j)$-th projection $X^4 \lra X^2$. 
By definition, $\rho_0$ is equal to the pull-back by $p_0 \times p_0: 
\ti{P} \lra X \times_{X^{(1)}} X$ of the descent 
data on $\cE$ relative to $X^{(1)}$ coming from $\cE'$, and 
it is equal to the pull-back of $\epsilon$ by 
$$ u \circ (p_0 \times p_0) = q_{02} \circ v: 
\ti{P} \lra P_{X,(0)}. $$
So we have 
$\rho_0 = (p_0 \times p_0)^*u^*\epsilon = v^*q_{02}^*\epsilon$. 
We see the equality $\rho_1 = v^*q_{13}^*\epsilon$ in the same way. 
On the other hand, we have $\pi_{23}^*\epsilon = v^*q_{23}\epsilon$, 
$\pi_{01}^*\epsilon = v^*q_{01}\epsilon$ by definition. Hence 
the commutativity of the diagram \eqref{21diag} follows from 
the cocycle condition for $\epsilon$. So we have proved the 
existence of the morphism $\epsilon': {p'}_1^*\cE' \lra {p'}_0^*\cE'$. \par 
$\epsilon'$ is an isomorphism because so is $\epsilon$. 
Also, the cocycle condition for $\epsilon'$ is reduced to that for 
$\epsilon$ because $\Phi$ is finite flat. Therefore, 
$\epsilon'$ is a $(-1)$-HPD-stratification and so 
$(\cE',\epsilon')$ forms an object in $\MIC^{(1)}(X^{(1)})^{\qn}$ 
with $F^*(\cE',\epsilon') = (\cE,\epsilon)$. Hence we have shown the 
essential surjectivity of the functor $F^*$ and so the proof of the 
theorem is now finished. 
\end{proof} 

\begin{remark}\label{qrem2}
Assume that we are in the situation of $\Hyp(\infty,2,2)$. 
Also, let 
 $\iota: \HIG(X^{(1)}_1)^{\qn} \lra \HIG(X^{(1)}_1)^{\qn}$ 
be the functor $(\cE,\theta) \mapsto (\cE,-\theta)$. Then, by 
Corollary \ref{fdcor}(1), the functor 
$$ F^*_2 \circ \iota: \HIG(X^{(1)}_1)^{\qn} \lra \MIC(X_1)^{\qn} $$ 
is an equivalence. In view of Remark \ref{qrem}, this reproves 
\cite[2.11]{ov} and a special case of \cite[5.8]{glsq} (the case 
$m=0$ in the notation there). 
So we can regard Theorem \ref{fd} as a generalization of their 
results in some sense. 
Note however 
that our result is slightly weaker than their result in the sense 
that we need the existence 
of the flat $p$-adic formal scheme $S$ with 
$S \otimes \Z/p\Z = S_1$. 
\end{remark}

\section{A comparison of de Rham cohomologies}

In this section, we prove a comparison theorem between the 
de Rham cohomology of certain $p^m$-connections on $p$-adic 
formal schemes and the 
de Rham cohomology of the pull-back of it by the level raising 
inverse image functor associated to 
certain lift of Frobenius. As an application, we prove the 
equivalence between the $\Q$-linearization of the 
category of nilpotent modules with 
integrable $p^m$-connections and the $\Q$-linearization of the 
category of nilpotent modules with 
integrable $p^{m-1}$-connections, under the existence of a nice 
left of Frobenius. \par 
Let $X \lra S$ be a smooth morphism of finite type between 
$p$-adic formal schemes flat over $\Z_p$ and let $m \in \N$. 
We define the category $\MIC^{(m)}(X)^{\qn}_{\b}$ as 
the category of projective systems 
$(\cE_{\b},\nabla_{\b})_{\b \in \N}$ in $\MIC^{(m)}(X)$ with 
$(\cE_{n},\nabla_{n}) \in \MIC^{(m)}(X)^{\qn}_n$ such that 
$(\cE_{n+1}, \nabla_{n+1}) \lra (\cE_n,\nabla_n)$ induces the 
isomorphism 
$(\cE_{n+1}, \allowbreak \nabla_{n+1}) \otimes \Z/p^n\Z \os{=}{\lra}
 (\cE_n,\nabla_n)$ for any $n \in \N$. 
If we put $(\cE,\nabla) := \varprojlim_n (\cE_{n},\nabla_n) \in 
\MIC^{(m)}(X)$, we have 
$(\cE,\nabla) \otimes \Z/p^n\Z = (\cE_n,\nabla_n)$ when 
each $\cE_n$ is quasi-coherent, by \cite[3.3.1]{berthelotI}. 
We define several nilpotent properties for objects in 
$\MIC^{(m)}(X)^{\qn}_{\b}$ 
which are stronger than quasi-nilpotence as follows: 

\begin{definition}\label{4.3}
Let $X \lra S, m$ be as above. 
\begin{enumerate}
\item 
For a smooth scheme $Y$ over $S_n := S \otimes \Z/p^n\Z$, 
an object $(\cE,\nabla)$ in $\MIC^{(m)}(Y)^{\qn}$ 
is called f-constant $($resp. lf-constant$)$ if $\cE$ is a quasi-coherent 
$\cO_{Y}$-module flat over 
$\Z/p^n\Z$ $($resp. locally free $\cO_{Y}$-module of 
finite rank$)$ and it is generated as $\cO_{Y}$-module 
by elements $e$ with $\nabla(e)=0$. For $l \in \N$, 
It is called f-nilpotent of length $\leq l$ 
$($resp. lf-nilpotent of length $\leq l)$
 if it can be written as an iterated extension of length $\leq l$ 
by f-constant $($resp. lf-constant$)$ objects. 
\item 
An object $(\cE_{\b},\nabla_{\b})$ in $\MIC^{(m)}(X)^{\qn}_{\b}$ is 
called f-nilpotent $($resp. lf-nilpotent$)$ 
if there exists some $l \in \N$ such that, for each $n \in \N$, 
there exists some etale surjective morphism $Y_n \lra X_n := X \otimes 
\Z/p^n\Z$ 
such that $(\cE_n,\nabla_n)|_{Y_n} \in \MIC^{(m)}(Y_n)^{\qn}$ is 
 f-nilpotent of length $\leq l$ 
$($resp. lf-nilpotent of length $\leq l)$. 
\item 
An object $(\cE_{\b},\nabla_{\b})$ in $\MIC^{(m)}(X)^{\qn}_{\b}$ is 
called nilpotent if 
it can be written as an iterated extension 
by the object $(\cO_{X_{\b}}, p^md)$. 
\end{enumerate}
We denote the full subcategory of $\MIC^{(m)}(X)^{\qn}_{\b}$ 
consisting of 
f-nilpotent $($resp. lf-nilpotent, nilpotent$)$ objects by 
$\MIC^{(m)}(X)^{\rm fn}_{\b}$ $($resp. 
$\MIC^{(m)}(X)^{\rm lfn}_{\b}$, $\MIC^{(m-1)}(X)^{\rm n}_{\b})$. 
\end{definition}

\begin{definition}
Let the notations be as in Definition \ref{4.3}. 
We call an object $(\cE,\nabla)$ in $\MIC^{(m)}(X)$ 
f-nilpotent $($resp. lf-nilpotent, nilpotent$)$ if there 
exists an object $(\cE_n, \nabla_n)_n$ in 
$\MIC^{(m)}(X)^{\rm fn}_{\b}$ $($resp. 
$\MIC^{(m)}(X)^{\rm lfn}_{\b}$, $\MIC^{(m)}(X)^{\rm n}_{\b})$ 
with $(\cE,\nabla) := \varprojlim_n (\cE_{n}, \allowbreak \nabla_n) \in 
\MIC^{(m)}(X)$, and denote the category of 
f-nilpotent $($resp. lf-nilpotent, nilpotent$)$ objects in 
$\MIC^{(m)}(X)$ by $\MIC^{(m)}(X)^{\rm fn}$ $($resp. 
$\MIC^{(m)}(X)^{\rm lfn}$, $\MIC^{(m)}(X)^{\rm n}).$ 
Since each $\cE_n$ is quasi-coherent for 
any $(\cE_n, \nabla_n)_n \in \MIC^{(m)}(X)^{\rm fn}_{\b}$, 
the functor 
$(\cE_n, \nabla_n) \mapsto (\cE,\nabla) := \varprojlim_n (\cE_{n},\nabla_n)$ 
induces the equivalence 
$$ \MIC^{(m)}(X)^{\rm fn}_{\b} \os{=}{\lra} 
\MIC^{(m)}(X)^{\rm fn}, $$
$$ \text{$($resp. $\MIC^{(m)}(X)^{\rm lfn}_{\b} \os{=}{\lra} 
\MIC^{(m)}(X)^{\rm lfn}, \,\,\,\, \MIC^{(m)}(X)^{\rm n}_{\b} \os{=}{\lra} 
\MIC^{(m)}(X)^{\rm n})$.} $$
$($The inverse is given by $(\cE,\nabla) \mapsto 
((\cE,\nabla) \otimes \Z/p^n\Z)_n.)$ 
\end{definition}

Note that we have implications 
$$ \text{nilpotent} \,\Longrightarrow\, \text{lf-nilpotent} 
\,\Longrightarrow\, \text{f-nilpotent}. $$ 
Recall that an object $(\cE,\nabla)$ in 
$\MIC^{(m)}(X)$ 
induces morphisms 
$$ \nabla_k: \cE \otimes_{\cO_X} \Omega^{k}_{X/S} \lra \cE \otimes_{\cO_X} 
\Omega^{k+1}_{X/S} $$ 
and they form a complex 
$$ 0 \lra \cE \os{\nabla}{\lra} \cE \otimes_{\cO_X} \Omega^1_{X/S} 
\os{\nabla_1}{\lra} \cE \otimes_{\cO_X} \Omega^2_{X/S} 
\os{\nabla_2}{\lra} \cdots, $$ 
which we call the de Rham complex of $(E,\nabla)$. 
We denote the cohomology sheaf of this complex by 
$\cH^i(\cE,\nabla)$ and the hypercohomology of it on $X$ 
by $H^i(X, (\cE,\nabla))$. \par 
To state the main result in this section, 
we give the following definition to fix the situation. 

\begin{definition}\label{hypabcc'}
In this definition, $\G_{m,T}$ denotes the scheme 
$\Spec T[t^{\pm 1}]$ for a scheme $T$ and the $p$-adic 
formal scheme $\Spf T\{t^{\pm 1}\}$ for a $p$-adic formal scheme $T$. 
For $a, b, c, c' \in \N \cup\{\infty\}$ with $a \geq b \geq c \geq c'$, 
we mean by $\Hyp(a,b,c,c')$ the following hypothesis$:$ 
Let $S_j \,(j \leq a), f_j: X_j \lra S_j \,(j \leq b), 
f^{(1)}_j: X^{(1)}_j \lra S_j \,(j \leq b), 
F_j: X_j \lra X^{(1)}_j \,(j \leq c)$ be as in $\Hyp(a,b,c)$. 
Also, we assume that there exists a Cartesian diagram 
\begin{equation}\label{hypo1''diag}
\begin{CD}
X_{c'} @>{F_{c'}}>> X^{(1)}_{c'} \\ 
@VVV @VVV \\ 
\G_{m,S_{c'}}^d @>p>> \G_{m,S_{c'}}^d 
\end{CD}
\end{equation}
Zariski locally $($where $d$ is the dimension of $X)$, 
where the vertical arrows are etale and the map $p$ is the 
morphism over $S$ which sends the coordinates to the $p$-th power of 
them. \par 
When $a = \infty$ $($resp. $b=\infty$, $c=\infty)$, 
we denote $S_a$ $($resp. $f_b: X_b \lra S_b$ and 
$f^{(1)}_b: X^{(1)}_b \lra S_b$, 
$F_c:  X_c \lra X^{(1)}_c)$ simply by 
$S$ $($resp. $f: X \lra S$ and 
$f^{(1)}: X^{(1)} \lra S$, 
$F:  X \lra X^{(1)})$. In the local situation where there exists 
the diagram \eqref{hypo1''diag}, 
let $t_1,...,t_d$ $($resp. $t'_1,...,t'_d)$ 
be the local coordinate on $X$ $($resp. $X^{(1)})$ induced from the 
canonical coordinate of $\fG_{m,S}$ via the left vertical arrow 
$($resp. the right vertical arrow$)$ in 
\eqref{hypo1''diag} and let us put $\tau_i := 1 \otimes t_i - 
t_i \otimes 1 \in \cO_{X^2}, \tau'_i := 1 \otimes t'_i - 
t'_i \otimes 1 \in \cO_{(X^{(1)})^2}$. 
\end{definition}

Roughly speaking, $\Hyp(a,b,c,c')$ means that $S_1$ is liftable to 
a scheme flat over $\Z/p^a\Z$, $f_1:X_1 \lra S_1$ and $f^{(1)}_1: 
X^{(1)}_1 \lra S_1$ is liftable to morphisms over $S_b$, 
the relative Frobenius $F_{X_1/S_1}: X_1 \lra X^{(1)}_1$ is liftable 
to a morphism over $S_c$ and nicely liftable to a morphism 
over $S_{c'}$. \par 
Then the main result in this section is stated in the following way. 

\begin{theorem}\label{25}
Let the notations be as in $\Hyp(\infty,\infty,\infty,\infty)$ and 
let $(\cE,\nabla)$ be an object in $\MIC^{(m)}(X^{(1)})^{\rm fn}$. 
Then the level raising inverse image functor 
$F^*: \MIC^{(m)}(X^{(1)}) \lra \MIC^{(m-1)}(X)$ induces the 
natural isomorphisms 
\begin{align*}
& \cH^0(\cE,\nabla) \os{=}{\lra} \cH^0(F^*(\cE,\nabla)). \\ 
& \cH^i(\cE,\nabla) \otimes {\Q} \os{=}{\lra} 
\cH^i(F^*(\cE,\nabla)) \otimes \Q \quad (i \in \N).
\end{align*}
\end{theorem}

Before the proof, we give several corollaries. 

\begin{corollary}\label{cor1}
Let the notations be as in $\Hyp(\infty,\infty,\infty,\infty)$ 
and 
let $(\cE,\nabla)$ be an object in $\MIC^{(m)}(X^{(1)})^{\rm fn}$. 
Then the level raising inverse image functor 
$F^*: \MIC^{(m)}(X^{(1)}) \lra \MIC^{(m-1)}(X)$ induces the 
isomorphisms 
\begin{align*}
& H^0(X^{(1)}, (\cE,\nabla)) \os{=}{\lra} 
H^0(X, F^*(\cE,\nabla)). \\ 
& H^i(X^{(1)}, (\cE,\nabla)) \otimes {\Q} \os{=}{\lra} 
H^i(X, F^*(\cE,\nabla)) \otimes \Q \quad (i \in \N).
\end{align*}
\end{corollary}

\begin{proof}
This follows from Theorem \ref{25} and the spectral sequences 
\begin{align*}
& E_2^{s,t} = H^s(X^{(1)},\cH^t(\cE,\nabla)) 
\,\Longrightarrow\, H^{s+t}(X^{(1)}, (\cE,\nabla)), \\ 
& E_2^{s,t} = H^s(X,\cH^t(F^*(\cE,\nabla))) 
\,\Longrightarrow\, H^{s+t}(X, F^*(\cE,\nabla)). 
\end{align*}
\end{proof}

To give two more corollaries, we introduce some more 
categories. Recall that, under $\Hyp(\infty,\infty,m+1)$, 
we have an equivalence 
\begin{align*}
\MIC^{(m-1)}(X)^{\qn}_n & \os{=}{\lra} 
\MIC^{(m-1)}(X_{n+m-1})^{\qn}_n 
\os{=}{\lra} 
D^{(m-1)}(X_{n+m-1}/S)^{\qn}_n \\ 
& \os{=}{\lra} 
\ol{D}^{(m-1)}(X_{n+m-1}/S)^{\qn}_n 
\os{=}{\lra} \ol{D}^{(m-1)}(X_{m+1}/S)^{\qn}_n
\end{align*}
for $n \in \N$. 
Let us define the category $\ol{D}^{(m-1)}(X_{m+1}/S)^{\qn}_{\b}$ as 
the category of projective systems 
$(\cE_{\b})_{\b \in \N}$ in 
$\ol{D}^{(m-1)}(X_{m+1}/S)^{\qn} := \bigcup_n 
\ol{D}^{(m-1)}(X_{m+1}/S)^{\qn}_n$ 
with $\cE_{n} \in \ol{D}^{(m-1)}(X_{m+1}/S)^{\qn}_n$ 
such that 
$\cE_{n+1} \lra \cE_n$ induces the 
isomorphism 
$\cE_{n+1} \otimes \Z/p^n\Z \os{=}{\lra}
 \cE_n$ for any $n \in \N$. Then the above equivalence induces 
the equivalence 
\begin{equation}\label{micd}
\MIC^{(m-1)}(X)^{\qn}_{\b} \os{=}{\lra} 
\ol{D}^{(m-1)}(X_{m+1}/S)^{\qn}_{\b}. 
\end{equation}
Then we define the category 
$\ol{D}^{(m-1)}(X_{m+1}/S)^{\rm lfn}_{\b}, 
\ol{D}^{(m-1)}(X_{m+1}/S)^{\rm n}_{\b}$ by the essential image of 
$\MIC^{(m-1)}(X)^{\rm lfn}_{\b}, \MIC^{(m-1)}(X)^{\rm n}_{\b}$ 
by the equivalence \eqref{micd}, respectively. 
Also, we can define the category 
$\ol{D}^{(m)}(X^{(1)}_{m+1}/S)^{\qn}_{\b}$, the equivalence  
\begin{equation*}
\MIC^{(m)}(X^{(1)})^{\qn}_{\b} \os{=}{\lra} 
\ol{D}^{(m)}(X^{(1)}_{m+1}/S)^{\qn}_{\b}, 
\end{equation*}
and the categories $\ol{D}^{(m)}(X^{(1)}_{m+1}/S)^{\rm lfn}_{\b}, 
\ol{D}^{(m)}(X^{(1)}_{m+1}/S)^{\rm n}_{\b}$ in the same way. 
Note that we have the level raising inverse image functor 
$$ F_{m+1}^*: \ol{D}^{(m)}(X_{m+1}^{(1)}/S)^{\qn}_{\b} \lra 
\ol{D}^{(m-1)}(X_{m+1}/S)^{\qn}_{\b} $$ 
induced by the ones 
$F_{m+1}^*: \ol{D}^{(m)}(X^{(1)}_{m+1}/S)^{\qn}_{n} \lra 
\ol{D}^{(m-1)}(X_{m+1}/S)^{\qn}_{n}$ for $n \in \N$. \par 
Then we have the following corollary: 

\begin{corollary}\label{cor2new}
\begin{enumerate}
\item 
Let the sitation be as in $\Hyp(\infty,\infty,m+1,m+1)$. 
Then the level raising inverse image functor 
$$F_{m+1}^*: \ol{D}^{(m)}(X^{(1)}_{m+1}/S)^{\qn}_{\b} \lra 
\ol{D}^{(m-1)}(X_{m+1}/S)^{\qn}_{\b}$$ 
induces the fully faithful functor 
\begin{equation}\label{dbeq}
F_{m+1}^*: \ol{D}^{(m)}(X^{(1)}_{m+1}/S)^{\rm lfn}_{\b} \lra 
\ol{D}^{(m-1)}(X_{m+1}/S)^{\rm lfn}_{\b}, 
\end{equation}
and via the induced fully faithful functor 
\begin{equation}\label{dbeq2}
\ol{D}^{(m)}(X^{(1)}_{m+1}/S)^{\rm lfn}_{\b,\Q} \lra 
\ol{D}^{(m-1)}(X_{m+1}/S)^{\rm lfn}_{\b,\Q} 
\end{equation}
between $\Q$-linealized categories, 
$\ol{D}^{(m)}(X^{(1)}_{m+1}/S)^{\rm lfn}_{\b,\Q}$ is a thick 
full subcategory of 
$\ol{D}^{(m-1)}(X_{m+1}/S)^{\rm lfn}_{\b,\Q}$. 
\item 
Let the sitation be as in $\Hyp(\infty,\infty,\infty,m+1)$. 
Then the level raising inverse image functor 
${F}^*: \MIC(X^{(1)})^{\qn} \lra 
\MIC^{(m-1)}(X)^{\qn}$ 
induces the fully faithful functor 
\begin{equation}\label{miceq}
{F'}^*: \MIC^{(m)}(X^{(1)})^{\rm lfn} \lra 
\MIC^{(m-1)}(X)^{\rm lfn}, 
\end{equation}
and via the induced fully faithful functor 
\begin{equation}\label{miceq2}
\MIC(X^{(1)})^{\rm lfn}_{\Q} \lra 
\MIC^{(m-1)}(X)^{\rm lfn}_{\Q} 
\end{equation}
between $\Q$-linealized categories, 
$\MIC^{(m)}(X^{(1)})^{\rm lfn}_{\Q}$ is a thick 
full subcategory of 
$\MIC^{(m-1)}(X)^{\rm lfn}_{\Q}$. 
\end{enumerate}
\end{corollary}

\begin{proof}
First we prove (2) under $\Hyp(\infty,\infty,\infty,\infty)$. 
Since the source and the 
target in \eqref{miceq} and \eqref{miceq2} are rigid tensor 
categories admitting internal hom's, it suffices to prove that, 
for $(\cE,\nabla) \in \MIC^{(m)}(X^{(1)})^{\rm lfn}$, the 
level raising inverse image functor induces the isomorphisms 
$$ H^0(X^{(1)},(\cE,\nabla)) \os{=}{\lra} H^0(X,F^*(\cE,\nabla)), $$
$$ H^1(X^{(1)},(\cE,\nabla)) \otimes \Q \os{=}{\lra} 
H^1(X,F^*(\cE,\nabla)) \otimes \Q,$$  
and it follows from Corollary \ref{cor1}. \par 
Next we prove (1). First we work locally and assume the 
existence of $F: X \lra X^{(1)}$ which lifts $F_{m+1}$ and 
assume that we are in the situation of $\Hyp(\infty,\infty,\infty,\infty)$ 
with this $F$. 
Then the functor \eqref{dbeq} is identical with the functor 
\eqref{miceq} and the latter is fully faithful by the argument 
in the previous paragraph. Hence the former is also fully faithful. \par 
Since the categories appearing in \eqref{dbeq} satisfies 
the descent property for Zariski coverings, we can deduce from 
the argument in the previous paragraph that the functor \eqref{dbeq} 
is fully faithful globally, under $\Hyp(\infty,\infty,m+1,m+1)$. 
Let us prove that $\ol{D}^{(m)}(X^{(1)}_{m+1}/S)^{\rm lfn}_{\b,\Q}$ is a thick 
full subcategory of 
$\ol{D}^{(m-1)}(X_{m+1}/S)^{\rm lfn}_{\b,\Q}$ via \eqref{dbeq2} globally 
under $\Hyp(\infty,\infty,m+1,m+1)$. 
Let us take an exact sequence 
\begin{equation}\label{extcl}
0 \lra F_{m+1}^*\cE' \os{f}{\lra} \cE \os{g}{\lra} F_{m+1}^*\cE'' \lra 0
\end{equation}
in 
$\ol{D}^{(m-1)}(X_{m+1}/S)^{\rm lfn}_{\b,\Q} = 
\MIC^{(m-1)}(X)^{\rm lfn}_{\Q} \subseteq \MIC^{(m-1)}(X)_{\Q}$ 
with 
$\cE', \cE''$ contained in $\ol{D}^{(m)}(X^{(1)}_{m+1}/S)^{\rm lfn}_{\b}$. 
It suffices to prove that $\cE$ is in the essential image of the 
functor \eqref{dbeq2}. 
We may assume that $f, g$ are morphisms in $\MIC^{(m-1)}(X)$ and that 
$g \circ f = 0$ in $\MIC^{(m-1)}(X)$. Then we have the diagram 
\begin{equation}\label{tuika1}
\begin{CD}
0 @>>>  F_{m+1}^*\cE' @>f>> \cE @>g>> F_{m+1}^*\cE'' @>>> 0 \\ 
@. @V{\alpha}VV @| @A{\beta}AA \\ 
0 @>>> \Ker g @>>> \cE @>>> \im g @>>> 0 
\end{CD}
\end{equation}
in $\MIC^{(m-1)}(X)$. Since $\alpha, \beta$ are isomorphisms in 
$\MIC^{(m-1}(X)_{\Q}$, we have some $a \in \N$ and morphisms 
$$ \alpha': \Ker g \lra F_{m+1}^*\cE', \quad 
\beta': F_{m+1}^*\cE'' \lra \im g$$ 
such that $\alpha \circ \alpha' = \alpha' \circ \alpha = p^a$, 
$\beta \circ \beta' = \beta' \circ \beta = p^a$. If we 
push the lower horizontal line in \eqref{tuika1} by $\alpha'$ and then 
push it by $\beta'$, we obtain the exact sequence of the form 
$$ 0 \lra F_{m+1}\cE' \lra \cE_1 \lra F_{m+1}^*\cE'' \lra 0 $$ 
in $\MIC^{(m-1)}(X)^{\rm lfn}$ with $\cE_1$ isomorphic to 
$\cE$ in the category $\MIC^{(m-1)}(X)_{\Q}$. So, by replacing 
$\cE$ by $\cE_1$, we may assume that the exact sequence 
\eqref{extcl} is the one in the category $\MIC^{(m-1)}(X)^{\rm lfn}$. 
Now let us take a Zariski covering $X = \bigcup_{\alpha} X_{\alpha}$ by 
finite number of open subschemes such that, if we denote the 
corresponding open covering of $X^{(1)}$ by 
$X^{(1)} = \bigcup_{\alpha} X_{\alpha}^{(1)}$, there exists a lift 
$F_{\alpha}: X_{\alpha} \lra X_{\alpha}^{(1)}$ of $F_{m+1}$ for each 
$\alpha$ with which we are in the situation of 
$\Hyp(\infty,\infty,\infty,\infty)$ (for $X_{\alpha}$). 
Then, since the functor $F_{m+1}^*|_{X_{\alpha}^{(1)}}$ is 
equal to $F_{\alpha}^*: \MIC^{(m)}(X^{(1)}_{\alpha})^{\rm lfn} \lra 
\MIC^{(m-1)}(X_{\alpha})^{\rm lfn}$, we have the isomorphism
$$ 
H^1(X^{(1)}_{\alpha}, \Inthom(\cE'',\cE)) \otimes \Q 
\os{=}{\lra} 
H^1(X_{\alpha}, \Inthom(F_{m+1}^*\cE'',F_{m+1}^*\cE')) \otimes \Q. 
$$ 
So, by multiplying the extension class $[\cE]$ of the exact sequence 
\eqref{extcl} by $p^b$ for some $b \in \N$, we may assume 
that there exists an exact sequence 
\begin{equation}\label{extcl2}
0 \lra \cE'|_{X_{\alpha}^{(1)}} \lra \cF_{\alpha} \lra 
\cE''|_{X_{\alpha}^{(1)}} \lra 0
\end{equation}
in $\MIC^{(m)}(X_{\alpha}^{(1)})^{\rm lfn}$ for each $\alpha$ 
such that $F_{m+1}^*\text{\eqref{extcl2}}$ is 
isomorphic to $\text{\eqref{extcl}}|_{X_{\alpha}}$. In particular, 
we have the isomorphism $i_{\alpha}: F_{m+1}^*\cF_{\alpha} \os{=}{\lra} 
\cE|_{X_{\alpha}}$. So, if we put 
$i_{\alpha\beta} := i_{\beta}^{-1} \circ i_{\alpha}: 
F_{m+1}^*\cF_{\alpha} \lra F_{m+1}^*\cF_{\beta}$ on $X_{\alpha} \cap 
X_{\beta}$, it satisfies the cocycle condition. Then, since 
$F_{m+1}^*$ is fully faithful, we see that there exists an object 
$\cF$ in $\ol{D}^{(m-1)}(X^{(1)}_{m+1}/S)^{\rm lfn}_{\b} =  
\MIC^{(m-1)}(X^{(1)})^{\rm lfn}$ with $F_{m+1}^*\cF = \cE$. 
Hence $\cE$ is in the essential image of 
the functor \eqref{dbeq2} and so we have shown that 
 $\ol{D}^{(m)}(X^{(1)}_{m+1}/S)^{\rm lfn}_{\b,\Q}$ is a thick 
full subcategory of 
$\ol{D}^{(m-1)}(X_{m+1}/S)^{\rm lfn}_{\b,\Q}$. So the proof of (1) is 
finished. \par 
Finally, since the functors \eqref{miceq}, 
\eqref{miceq2} are identified with \eqref{dbeq}, \eqref{dbeq2}, 
the assertion (2) in general case is an immediate consequence of 
the assertion (1). So the proof of corollary is finished. 
\end{proof}

\begin{corollary}\label{cor3new}
If we are in the situation of $\Hyp(\infty,\infty,m+1,m+1)$, 
the level raising inverse image functor 
$F_{m+1}^*: \ol{D}^{(m)}(X^{(1)}_{m+1}/S)^{\qn}_{\b} \lra 
\ol{D}^{(m-1)}(X_{m+1}/S)^{\qn}_{\b}$ 
induces the fully faithful functor 
\begin{equation}\label{dbeqn}
F_{m+1}^*: \ol{D}^{(m)}(X^{(1)}_{m+1}/S)^{\rm n}_{\b} \lra 
\ol{D}^{(m-1)}(X_{m+1}/S)^{\rm n}_{\b}, 
\end{equation}
giving the equivalence 
\begin{equation}\label{dbeqn2}
\ol{D}^{(m)}(X^{(1)}_{m+1}/S)^{\rm n}_{\b,\Q} \os{=}{\lra} 
\ol{D}^{(m-1)}(X_{m+1}/S)^{\rm n}_{\b,\Q} 
\end{equation}
between $\Q$-linearized categories. Also, 
if we are in the situation of $\Hyp(\infty,\infty,\infty,m+1)$, 
the level raising inverse image functor 
${F}^*: \MIC(X^{(1)})^{\qn} \lra 
\MIC^{(m-1)}(X)^{\qn}$ 
induces the fully faithful functor 
\begin{equation}\label{miceqn}
{F}^*: \MIC^{(m)}(X^{(1)})^{\rm n} \lra 
\MIC^{(m-1)}(X)^{\rm n}, 
\end{equation}
giving the equivalence 
\begin{equation}\label{miceqn2}
\MIC(X^{(1)})^{\rm n}_{\Q} \os{=}{\lra} 
\MIC^{(m-1)}(X)^{\rm n}_{\Q} 
\end{equation}
between $\Q$-linealized categories. 
\end{corollary}

\begin{proof}
The full faithfulness of \eqref{dbeqn} and \eqref{miceqn} 
follows from that of \eqref{dbeq} and \eqref{miceq}. Also, by the same 
argument as the proof of Corollary \ref{cor2new}, we see that 
$\ol{D}^{(m)}(X^{(1)}_{m+1}/S)^{\rm n}_{\b,\Q}$ 
(resp. $\MIC(X^{(1)})^{\rm n}_{\Q}$) 
is a thick full subcategory of the category 
$\ol{D}^{(m-1)}(X_{m+1}/S)^{\rm n}_{\b,\Q}$ 
(resp. $\MIC^{(m-1)}(X)^{\rm n}_{\Q}$) via the functor 
\eqref{dbeqn2} (resp. \eqref{miceqn2}). Since any object in the categories 
$\ol{D}^{(m-1)}(X_{m+1}/S)^{\rm n}_{\b,\Q}$, 
$\MIC^{(m-1)}(X)^{\rm n}_{\Q}$ is written as an iterated extension of 
trivial objects, thickness above implies the equivalence of the 
functors \eqref{dbeqn2}, \eqref{miceqn2}. 
\end{proof}

Note that Corollaries \ref{cor2new}, \ref{cor3new} give possible answers to 
Question \ref{q}. \par 
We give the proof of Theorem \ref{25}. 

\begin{proof}[Proof of Theorem \ref{25}] 
We have the homomorphism of complexes from the de Rham complex of 
$(\cE,\nabla)$ to that of $F^*(\cE,\nabla)$ induced by the level 
raising inverse image functor. So, to prove the theorem, we may work 
locally. So we may assume that the Cartesian diagram \eqref{hypo1''diag} 
exists globally. Then we have $F^*\cE = \oplus_{a} t^a\cE$, where 
$a$ runs through the set 
$I := \{a = (a_i)_i \in \N^d \,|\, 0 \leq a_i \leq 
p-1\}$. If we express 
$\nabla: \cE \lra \cE \otimes_{\cO_{X^{(1)}}} \Omega^1_{X^{(1)}/S} 
= \oplus_{i=1}^d \cE {t'}_i^{-1}dt'_i$ by 
$\nabla(e) = \sum_i\nabla_i(e) {t'}_i^{-1}dt'_i$, 
$F^*\nabla: F^*\cE \lra F^*\cE \otimes_{\cO_{X}} \Omega^1_{X/S} 
= \oplus_{i=1}^d F^*\cE t_i^{-1}dt_i$ has the property 
$F^*\nabla(e) = \sum_{i}\nabla_i(e)t_i^{-1}dt_i$. 
(Here we wrote the element $F^*(e) \in F^*\cE$ simply by $e$.) 
Hence we have, for $a \in S$, 
\begin{align}
F^*\nabla(t^ae) & 
= t^aF^*\nabla(e) + p^{m-1}\sum_i a_it^ae t_i^{-1}dt_i \label{25eq} \\ 
& = 
t^a(\sum_i\{(1 \otimes \nabla_i)+p^{m-1}a_i\}(e)t_i^{-1}dt_i. \nonumber 
\end{align}
For $a \in S$, let $\theta_a: \cE \lra 
\cE \otimes_{\cO_{X^{(1)}}} \Omega^1_{X^{(1)}/S} 
= \oplus_{i=1}^d \cE {t'}_i^{-1}dt'_i$ be the linear map 
$e \mapsto \sum_i a_ie{t'}_i^{-1}dt'_i$. Then, by \eqref{25eq}, 
we have 
$F^*(\cE,\nabla) = (F^*\cE,F^*\nabla) = \oplus_{a} 
(\cE, \nabla + p^{m-1}\theta_a)$ via 
the identification of $\oplus_{i=1}^d F^*\cE t_i^{-1}dt_i$ and 
$\oplus_a \oplus_{i=1}^d \cE {t'}_i^{-1}dt'_i$ on the target. 
Hence we have 
\begin{equation}\label{ds}
\cH^i(F^*(\cE,\nabla)) = \oplus_{a} \cH^i(\cE,\nabla+p^{m-1}\theta_a). 
\end{equation}
Let us consider the following claim: \\ 
\quad \\
{\bf claim 1.} \,\,\, $\cH^i(\cE,\nabla+p^{m-1}\theta_a) \otimes \Q =0$ 
for all $i \in \N$ and $a \in I, \not= 0$. \\
\quad \\
First we prove that the claim 1 implies the theorem. 
We see easily that the claim 1 and the equality \eqref{ds} implies 
the isomorphism 
$\cH^i(\cE,\nabla) \otimes {\Q} \os{=}{\lra} 
\cH^i(F^*(\cE,\nabla)) \otimes \Q$. Also, since $\cE_{n} := \cE \otimes 
\Z/p^n\Z$ is flat 
over $\Z/p^n\Z$ for each $n$, $\cE$ is flat over $\Z_p$ and hence 
the both hand sides of \eqref{ds} are flat over $\Z_p$ when $i=0$. 
So the claim implies the equalities 
$\cH^0(\cE,\nabla+p^{m-1}\theta_a) =0$ 
for $a \in I, \not= 0$, and this and \eqref{ds} implies the 
isomorphism 
$\cH^0(\cE,\nabla) \os{=}{\lra} \cH^0(F^*(\cE,\nabla))$, as desired. 
So it suffices to prove the claim 1. \par 
In the following, we denote the map $\theta_a \otimes \Z/p^n\Z: \cE_n \lra 
\cE_n \otimes_{\cO_{X^{(1)}}} \Omega^1_{X^{(1)}/S}$ also by 
$\theta_a$, by abuse of notation. 
Let us take $l \in \N$ and etale surjective morphism 
$Y_n \lra X^{(1)}_n \,(n \in \N)$ such that $(\cE_n,\nabla_n) |_{Y_n}$ is 
f-nilpotent of length $\leq l$. 
For an open subset 
$U \subseteq X^{(1)}$, let 
$(\cE, \nabla+p^{m-1}\theta_a)(U)$, 
$(\cE_n, \nabla_n+p^{m-1}\theta_a)(U)$ be 
the complex 
$\Gamma(U, \cE \otimes_{\cO_X} \Omega^{\b}_{X/S})$, 
$\Gamma(U, \cE_n \otimes_{\cO_X} \Omega^{\b}_{X/S})$ 
induced from the de Rham complex of 
$(\cE, \nabla+p^{m-1}\theta_a)$, 
$(\cE_n, \nabla_n+p^{m-1}\theta_a)$ respectively 
and let 
$H^i((\cE, \nabla+p^{m-1}\theta_a)(U))$, 
$H^i((\cE_n, \nabla_n+p^{m-1}\theta_a)(U))$ be the $i$-th 
cohomology of this complex. 
To prove claim 1, first prove the following claim: \\
\quad \\
{\bf claim 2.} \,\,\, 
For $n \in \N$ and an open affine sub formal scheme $U \subseteq X^{(1)}$, 
we have $p^{2l(m-1)(i+1)}H^i((\cE_n, \nabla+p^{m-1}\theta_a)(U)) = 0$ 
for all $i \in \N$ and $a \in I, \not= 0$. \\ 
\quad \\ 
We prove claim 2. Let us put $n' := n + m-1$ and let 
$Y \lra X^{(1)}$ be an etale morphism lifting $Y_{n'} 
\lra X^{(1)}_{n'}$. Take an 
affine etale \v{C}ech hypercovering $U_{\b} \lra U$ such that 
$U_0 \lra U$ is an refinement of $Y \times_{X^{(1)}} U \lra U$. 
Then, since $E_n$ is quasi-coherent, we have 
$$ H^i((\cE_n, \nabla+p^{m-1}\theta_a)(U)) = 
H^i((\cE_n, \nabla+p^{m-1}\theta_a)(U_{\b})) $$
(where the right hand side denotes the $i$-th cohomology of the 
double complex $(\cE_n, \nabla+p^{m-1}\theta_a)(U_{\b})$) and so 
we have the spectral sequence 
$$ E_2^{s,t} := H^t((\cE_n, \nabla+p^{m-1}\theta_a)(U_{s}))
 \,\Longrightarrow\, H^{s+t}((\cE_n, \nabla+p^{m-1}\theta_a)(U)). $$
Hence it suffices to prove 
$p^{2l(m-1)}H^i((\cE_n, \nabla+p^{m-1}\theta_a)(U_s)) = 0$ 
for all $s$. 
By assumption, $(\cE_{n'},\nabla_{n'})|_{U_s}$ is written as 
an iterated extension 
of f-constant objects $(\cE_{n',j}, \nabla_{n',j})\,(1 \leq j \leq l)$ 
in $\MIC^{(m)}(U_s \otimes \Z/p^{n'}\Z)$. 
Then $(\cE_{n'},\nabla_{n'}+p^{m-1}\theta_a)|_{U_s}$ is written as 
an iterated extension 
of $(\cE_{n',j}, \nabla_{n',j}+p^{m-1}\theta_a)\,(1 \leq j \leq l)$, 
and $(\cE_{n},\nabla_{n}+p^{m-1}\theta_a)|_{U_s}$ is written 
an iterated extension 
of $(\cE_{n,j}, \nabla_{n,j}+p^{m-1}\theta_a) := 
(\cE_{n',j},\nabla_{n',j}+p^{m-1}\theta_a) \otimes \Z/p^n\Z 
\,(1 \leq j \leq l)$. So it suffices to prove 
$$p^{2(m-1)}H^i((\cE_{n,j}, \nabla_{n,j}+p^{m-1}\theta_a)(U_s)) = 0.$$ 
Since $(\cE_{n',j}, \nabla_{n',j})$ is f-constant, we have 
$\nabla_{n',j}(\cE_{n',j}) 
\subseteq p^m\cE_{n',j} \otimes \Omega^{1}_{X^{(1)}/S}$. 
Hence we can factorize $\nabla_{n',j}$ as 
$$ \cE_{n',j} \os{\ol{\nabla}_{n',j}}{\lra} 
\cE_{n,j} \otimes \Omega^{1}_{X^{(1)}/S} 
\os{p^{m-1}}{\lra} \cE_{n',j} \otimes \Omega^{1}_{X^{(1)}/S}. $$
Let $\ol{\ol{\nabla}}_j: \cE_{n,j} \lra 
\cE_{n,j} \otimes \Omega^{1}_{X^{(1)}/S}$ be 
$\ol{\nabla}_{n',j} \otimes \Z/p^n\Z$. 
Then it is a f-nilpotent $p$-connection and we have the 
commutative diagram 
\begin{equation}\label{dr1}
\begin{CD}
\cE_{n,j} @>{p^{m-1}}>> \cE_{n,j} @>{=}>> \cE_{n,j} \\ 
@V{\nabla_{n,j}+p^{m-1}\theta_a}VV 
@V{\ol{\ol{\nabla}}_j+\theta_a}VV 
@V{\nabla_{n,j}+p^{m-1}\theta_a}VV \\ 
\cE_{n,j} \otimes \Omega^1_{X^{(1)}/S} @>{=}>> 
\cE_{n,j} \otimes \Omega^1_{X^{(1)}/S} @>{p^{m-1}}>> 
\cE_{n,j} \otimes \Omega^1_{X^{(1)}/S}. 
\end{CD}
\end{equation}
For $k \in \N$, let us denote the homomorphism 
$\cE_{n,j} \otimes \Omega^k_{X^{(1)}/S} \otimes 
\lra \cE_{n,j} \otimes \Omega^{k+1}_{X^{(1)}/S}$ induced by 
$\nabla_{n,j}+p^{m-1}\theta_a$ (resp. 
$\ol{\ol{\nabla}}_j+\theta_a$) by 
$(\nabla_{n,j}+p^{m-1}\theta_a)_k$ (resp. 
$(\ol{\ol{\nabla}}_j+\theta_a)_k$). Then the commutativity of 
\eqref{dr1} implies that of the following diagram: 
\begin{equation}\label{dr2}
\begin{CD}
\cE_{n,j} \otimes \Omega^{i-1}_{X^{(1)}/S} 
@>{p^{2(m-1)}}>> 
\cE_{n,j}  \otimes \Omega^{i-1}_{X^{(1)}/S} 
@>{=}>> 
\cE_{n,j}  \otimes \Omega^{i-1}_{X^{(1)}/S} \\ 
@V{(\nabla_{n,j}+p^{m-1}\theta_a)_{i-1}}VV 
@V{(\ol{\ol{\nabla}}_j+\theta_a)_{i-1}}VV 
@V{(\nabla_{n,j}+p^{m-1}\theta_a)_{i-1}}VV \\ 
\cE_{n,j} \otimes \Omega^i_{X^{(1)}/S} @>{p^{m-1}}>> 
\cE_{n,j} \otimes \Omega^i_{X^{(1)}/S} @>{p^{m-1}}>> 
\cE_{n,j} \otimes \Omega^i_{X^{(1)}/S} \\ 
@V{(\nabla_{n,j}+p^{m-1}\theta_a)_{i}}VV 
@V{(\ol{\ol{\nabla}}_j+\theta_a)_{i}}VV 
@V{(\nabla_{n,j}+p^{m-1}\theta_a)_{i}}VV \\ 
\cE_{n,j} \otimes \Omega^{i+1}_{X^{(1)}/S} @>{=}>> 
\cE_{n,j} \otimes \Omega^{i+1}_{X^{(1)}/S} @>{p^{2(m-1)}}>> 
\cE_{n,j} \otimes \Omega^{i+1}_{X^{(1)}/S}. 
\end{CD}
\end{equation}
From the commutative diagram \eqref{dr2}, 
we see that it suffices to prove 
the equality $H^i((\cE_{n,j}, \ol{\ol{\nabla}}_{j}+\theta_a)(U_s)) = 0$
to prove the claim 2. Since $(\cE_{n,j}, \ol{\ol{\nabla}}_j)$ is 
a f-constant $p$-connection, 
we see that $\ol{\ol{\nabla}}(\cE_{n,j}) \subseteq 
p\cE_{n,j} \otimes \Omega^{1}_{X^{(1)}/S}$. So we can write 
$(\cE_{n,j}, \ol{\ol{\nabla}}_{j}+\theta_a)$ as an iterated extension 
by the Higgs module $(\cE_{1,j},\theta_a)$ (where 
$\cE_{1,j} = \cE_{n,j} \otimes \Z/p\Z$), and by 
\cite[2.2.2]{oguscartier}, we have 
$H^i((\cE_{1,j},\theta_a)(U_s))=0$ for any $i \in \N$. Hence we have 
$H^i((\cE_{n,j}, \ol{\ol{\nabla}}_{j}+\theta_a)(U_s)) = 0 \,(i \in \N)$ 
and so the proof of claim 2 is finished. \par 
Finally we prove the claim 1. For any open affine $U \subseteq X^{(1)}$, 
we have the exact sequence 
\begin{align*}
0 & \lra 
{\varprojlim}^1 H^{i-1}((\cE_n,\nabla_n+p^{m-1}\theta_a)(U)) 
\lra H^i((\cE,\nabla+p^{m-1}\theta_a)(U)) \\ & \lra 
\varprojlim_n H^i((\cE_n,\nabla_n+p^{m-1}\theta_a)(U)) \lra 0 
\end{align*}
and by claim 1, the first and the third term are killed by 
$p^{2l(m-1)(i+1)}$. Hence $H^i((\cE,\nabla+p^{m-1}\theta_a)(U))$ 
is killed by $p^{4l(m-1)(i+1)}$ and so we have 
$\cH^i(\cE,\nabla+p^{m-1}\theta_a) \otimes \Q = 0$, as desired. 
So the proof of the theorem is now finished. 
\end{proof}

\section{Modules with integrable $p^m$-Witt-connections}

In this section, we define the notion of modules with 
integrable $p^m$-Witt-connections. We also define 
the level raising inverse image functor 
from the categories of modules with integrable 
$p^m$-Witt-connections to that of $p^{m-1}$-Witt-connections, 
and prove that it induces the equivalence of categories between 
$\Q$-linearized categories when restricted to nilpotent objects. 
This is the Witt analogue of the equivalence 
\eqref{miceqn2} in Corollary \ref{cor3new}. The categories defined in this 
section is more complicated than the categories defined in Section 1, 
but the equivalance proven in this section 
has the advantage that we need no assumpion on liftability of 
Frobenius. \par 
Throughout in this sction, we fix a perfect scheme $S_1$ of 
characteristic $p>0$. 

\begin{definition}
Let us assume given a smooth morphism $X_1 \lra S_1$ of finite type. 
Let 
$W\cO_{X_1} = \varprojlim_n W_n\cO_{X_1}$ be the ring of Witt vectors 
of $\cO_{X_1}$ and let $W\Omega_{X_1} = \varprojlim_n W_n\Omega_{X_1}^{\b}$ 
be the de Rham-Witt complex of $X_1$. 
Then a $p^m$-Witt-connection on a $W\cO_{X_1}$-module $\cE$ is 
an additive map $\nabla: \cE \lra \cE \otimes_{W\cO_{X_1}} 
W\Omega^1_{X_1}$ 
satisfying $\nabla(fe) = f\nabla(e) + p^me \otimes df$ for 
$e \in \cE, f \in W\cO_{X_1}$. We call a $p^m$-Witt-connection also as 
a Witt-connection of level $-m$. 
\end{definition}

When we are given a $W\cO_{X_1}$-module with $p^m$-Witt-connection 
$(\cE,\nabla)$, we can define the additive map 
$\nabla_k: \cE \otimes_{W\cO_{X_1}} W\Omega^k_{X_1} \lra 
\cE \otimes_{W\cO_{X_1}} W\Omega^{k+1}_{X_1}$ which is characterized by 
$\nabla_k(e \otimes \omega) = \nabla(e) \wedge \omega + 
p^me \otimes \omega$. 

\begin{definition}
With the notation above, we call $(\cE,\nabla)$ integrable 
if we have $\nabla_1 \circ \nabla = 0$. We denote the category of 
$W\cO_{X_1}$-modules with integrable $p^m$-Witt-connection by 
$\MIWC^{(m)}(X_1)$. 
\end{definition}

It is easy to see that the category $\MIWC^{(m)}(X_1)$ is functorial 
with respect to $X_1$. We define the notion of nilpotence as follows: 

\begin{definition}
An object $(\cE,\nabla)$ in $\MIWC^{(m)}(X_1)$ is 
called nilpotent if 
it can be written as an iterated extension 
by the object $(W\cO_{X_1}, p^md)$. We denote the full subcategory of 
$\MIWC^{(m)}(X_1)$ consisting of nilpotent objects by 
$\MIWC^{(m)}(X_1)^{\rm n}$. 
\end{definition}

Let $F_{S_1}: S_1 \lra S_1$ be the Frobenius endomorphism and let 
$X_1^{(1)} := S_1 \times_{F,S_1} X_1$. 
Then we define the homomorphism of differential graded algebras 
\begin{equation}\label{Fm}
(W\Omega_{X_1^{(1)}}^{\b},p^md) \lra (W\Omega_{X_1}^{\b},p^{m-1}d)
\end{equation}
by the composite 
\begin{align*}
(W\Omega_{X_1^{(1)}}^{\b},p^md) 
& \os{=}{\lra} 
\varprojlim_n (W_n\cO_{S_1} \otimes_{W_{n+1}\cO_{S_1}} 
W_{n+1}\Omega_{X_1^{(1)}}^{\b}, \id \otimes p^md) \\ 
& \os{=}{\lra} \varprojlim_n 
(W_n\cO_{S_1} \otimes_{F, W_{n+1}\cO_{S_1}} 
W_{n+1}\Omega_{X_1}^{\b}, \id \otimes p^md) 
\os{\id \otimes F}{\lra} 
(W_{n}\Omega_{X_1}^{\b}, p^{m-1}d) 
\end{align*}
(where $F$ is the Frobenius in de Rham-Witt complexes). 
We denote the homomorphism \eqref{Fm} also by $F$. Also, 
for a $W\cO_{X_1^{(1)}}$-module $\cE$ and $e \in \cE$, denote 
$W\cO_{X_1} \otimes_{F,W\cO_{X_1^{(1)}}} \cE$ and 
$1 \otimes e \in W\cO_{X_1} \otimes_{F,W\cO_{X_1^{(1)}}} \cE$ 
simply by $F_*\cE, F_*(e)$. Then, for $(\cE,\nabla) \in 
\MIWC^{(m)}(X_1^{(1)})$, we define $F_*(\cE,\nabla) \in 
\MIWC^{(m-1)}(X_1)$ by 
$F_*(\cE,\nabla) := (F_*\cE,F_*\nabla)$, where $F_*\cE$ is as above and 
$F_*\nabla: F_*\cE \lra F_*\cE \otimes_{W\cO_{X_1}} W\Omega_{X_1}^1$ 
is the map defined by $F_*\nabla(fF_*(e)) = 
fF_*(\nabla(e)) + p^{m-1}F_*(e)df$ for $e \in \cE, f \in W\cO_{X_1}$. 
(The integrability of $F_*(\cE,\nabla)$ follows from the well-known 
formula $dF = pFd$ for de Rham-Witt complex.) So we have the functor 
$$F_*: \MIWC^{(m)}(X_1^{(1)}) \lra \MIWC^{(m-1)}(X_1), $$
and it is easy to see that it induces the functor 
\begin{equation}\label{Fn}
F_*: \MIWC^{(m)}(X_1^{(1)})^{\rm n} \lra \MIWC^{(m-1)}(X_1)^{\rm n}. 
\end{equation}

An object $(\cE,\nabla)$ in 
$\MIWC^{(m-1)}(X_1)$ (resp. $\MIWC^{(m)}(X^{(1)}_1)$) 
induces the complex
\begin{align*}
& 0 \lra \cE \os{\nabla}{\lra} \cE \otimes_{W\cO_{X_1}} 
W\Omega^1_{X_1} 
\os{\nabla_1}{\lra} \cE \otimes_{W\cO_{X_1}} W\Omega^2_{X_1} 
\os{\nabla_2}{\lra} \cdots \\ 
& \text{(resp. }
0 \lra \cE \os{\nabla}{\lra} \cE \otimes_{W\cO_{X_1^{(1)}}} 
W\Omega^1_{X_1^{(1)}} 
\os{\nabla_1}{\lra} \cE \otimes_{W\cO_{X_1^{(1)}}} W\Omega^2_{X_1^{(1)}} 
\os{\nabla_2}{\lra} \cdots \text{ ),}
\end{align*}
which is 
called the de Rham complex of $(E,\nabla)$. 
We denote the cohomology sheaf of this complex by 
$\cH^i(\cE,\nabla)$ and the hypercohomology of it on $X_1$ (resp. $X_1^{(1)}$) 
by $H^i(X_1, (\cE,\nabla))$ (resp.  $H^i(X^{(1)}_1, (\cE,\nabla))$). \par 
Then we have the following theorem, which is the Witt version 
of Theorem \ref{25}: 

\begin{theorem}\label{25witt}
Let the notations be as above and 
let $(\cE,\nabla) \in \MIWC^{(m)}(X^{(1)}_1)^{\rm n}$. 
Then the level raising inverse image functor 
\eqref{Fn}
induces the natural isomorphisms 
\begin{align*}
& \cH^0(\cE,\nabla) \os{=}{\lra} \cH^0(F_*(\cE,\nabla)). \\ 
& \cH^i(\cE,\nabla) \otimes {\Q} \os{=}{\lra} 
\cH^i(F_*(\cE,\nabla)) \otimes \Q \quad (i \in \N).
\end{align*}
\end{theorem}
 
We have the following corollaries, which we can prove in the same way 
as Corollaries \ref{cor1}, \ref{cor3}. (So we omit the proof.) 

\begin{corollary}\label{cor1witt}
Let $S_1, X_1, X_1^{(1)}$ be as above and 
let $(\cE,\nabla) \in \MIWC^{(m)}(X^{(1)}_1)^{\rm n}$. 
Then the level raising inverse image functor \eqref{Fn}
induces the isomorphisms 
\begin{align*}
& H^0(X^{(1)}_1, (\cE,\nabla)) \os{=}{\lra} 
H^0(X_1, F_*(\cE,\nabla)). \\ 
& H^i(X^{(1)}_1, (\cE,\nabla)) \otimes {\Q} \os{=}{\lra} 
H^i(X_1, F_*(\cE,\nabla)) \otimes \Q \quad (i \in \N).
\end{align*}
\end{corollary}

\begin{corollary}\label{cor3}
Let $S_1, X_1, X_1^{(1)}$ be as above. 
Then the level raising inverse image functor \eqref{Fn} is a 
fully faithful functor which gives the equivalence 
\begin{equation}\label{Fn2}
\MIWC(X^{(1)}_1)^{\rm n}_{\Q} \os{=}{\lra} 
\MIWC^{(m-1)}(X_1)^{\rm n}_{\Q} 
\end{equation}
between $\Q$-linealized categories. 
\end{corollary}

\begin{proof}[Proof of Theorem \ref{25witt}]
Since we may work locally, we may assume that 
$S_1 = \Spec R, \allowbreak X_1 = \Spec B_1$ are affine. 
(Hence $W_nS = \Spec W_n(R)$.) Also, we may assume that 
the homomorphism $R \lra B_1$ corresponding to 
$X_1 \lra S_1$ factors as 
\begin{equation}\label{factor}
R \lra A_1 := R[T_i]_{1 \leq i \leq d} \lra 
A'_1 := R[T_i^{\pm 1}]_{1 \leq i \leq d} \lra 
B_1, 
\end{equation}
where the last map is etale. Recall \cite[3.2]{lz} that 
a Frobenius lift of a smooth $R$-algebra $C_1$ is a projective system 
of data $(C_n, \phi_n, \delta_n)_{n\geq 1}$ 
(with $C_1$ the same as the given one), where 
\begin{itemize}
\item 
$C_n \,(n \geq 1)$ is a smooth lifting over $W_n(R)$ of $C_1$ satisfying 
$W_n(R) \otimes_{W_{n+1}(R)} C_{n+1} \os{=}{\lra} C_n$. 
\item 
$\phi_n \,(n > 1)$ is a map $C_n \lra C_{n-1}$ over 
$F: W_n(R) \lra W_{n-1}(R)$ which is compatible with 
the Frobenius morphism $C_1 \lra C_1$. 
\item 
$\delta_n \,(n \geq 1)$ is a map 
$C_n \lra W_n(C_1)$ such that the diagram 
\begin{equation*}
\begin{CD}
C_n @>{\delta_n}>> W_n(C_1) \\ 
@V{\phi_n}VV @V{F}VV \\ 
C_{n-1} @>{\delta_{n-1}}>> W_{n-1}(C_1) 
\end{CD}
\end{equation*}
is commutative. 
\end{itemize}
Now let us define a Frobenius lift $(A_n, \phi_n, \delta_n)_n, 
(A'_n, \phi'_n, \delta'_n)$ of $A_1, A'_1$ by 
\begin{align*}
& A_n := W_n(R)[T_i]_{1 \leq i \leq d}, \qquad 
\phi_n(T_i) = T_i^p, \qquad \delta_n(T_i) = [T_i], \\ 
& A'_n := W_n(R)[T_i^{\pm 1}]_{1 \leq i \leq d}, \qquad 
\phi_n(T_i^{\pm 1}) = T_i^{\pm p}, \qquad 
\delta_n(T_i^{\pm 1}) = [T_i^{\pm 1}], 
\end{align*}
where $[-]$ denotes the Teichm\"uller lift. Then we have a map 
$(f_n)_n: (A_n,\phi_n,\delta_n)_n \lra (A'_n, \phi'_n, \delta'_n)_n$ over 
the map $A_1 \lra A'_1$ in \eqref{factor} such that the diagrams 
\begin{equation}\label{Flcomp1}
\begin{CD}
A'_n @>{\phi'_n}>> A'_{n-1} \\ 
@A{f_n}AA @A{f_{n-1}}AA \\ 
A_n @>{\phi_n}>> A_{n-1}, 
\end{CD}
\quad \quad \quad \quad \quad 
\begin{CD}
A'_n @>{\delta'_n}>> W_n(A'_1) \\ 
@A{f_n}AA @A{W_n(f_1)}AA \\ 
A_n @>{\delta_n}>> W_n(A_1) 
\end{CD}
\end{equation}
are cocartesian. Moreover, by the proof of \cite[3.2]{lz}, there exists 
a Frobenius lift $(B_n, \psi_n, \epsilon_n)_n$ of $B_1$ and a map 
$(g_n)_n: (A'_n, \phi'_n, \delta'_n)_n \lra 
(B_n, \psi_n, \epsilon_n)_n$ over the map $A'_1 \lra B_1$ in 
\eqref{factor} such that the diagrams 
\begin{equation}\label{Flcomp2}
\begin{CD}
B_n @>{\psi'_n}>> B_{n-1} \\ 
@A{g_n}AA @A{g_{n-1}}AA \\ 
A'_n @>{\phi'_n}>> A'_{n-1}, 
\end{CD}
\quad \quad \quad \quad \quad 
\begin{CD}
B_n @>{\epsilon_n}>> W_n(B_1) \\ 
@A{g_n}AA @A{W_n(g_1)}AA \\ 
A'_n @>{\delta'_n}>> W_n(A'_1) 
\end{CD}
\end{equation}
are cocartesian. Let us define the Frobenius lift 
$(A^{(1)}_n, \phi^{(1)}_n, \delta^{(1)}_n)_n$ of 
$A^{(1)}_1 := R \otimes_{F,R} A_1$ by 
$$ 
A_n^{(1)} := W_n(R) \otimes_{F,W_n(R)} A_n, \qquad 
\phi^{(1)}_n := F \otimes \phi_n, \qquad 
\delta^{(1)}_n := \id \otimes \delta_n, 
$$
and define the Frobenius lift 
$({A'}^{(1)}_n, {\phi'}^{(1)}_n, {\delta'}^{(1)}_n)_n$ of 
${A'}^{(1)}_1 := R \otimes_{F,R} A'_1$ and 
the Frobenius lift 
$(B^{(1)}_n, \psi^{(1)}_n, \epsilon^{(1)}_n)_n$ of 
$B^{(1)}_1 := R \otimes_{F,R} B_1$ in the same way. \par 
Let us put $B := \varprojlim_n B_n, B^{(1)} := \varprojlim_n B^{(1)}_n$. 
Recall that we have the map 
$F: (W\Omega_{B_1^{(1)}}^{\b},p^md) \lra (W\Omega_{B_1}^{\b},p^{m-1}d)$ 
of differential graded algebras defined in \eqref{Fm}. On the other hand, 
we have the map 
$\ol{F}^*: (\Omega^{\b}_{B^{(1)}}, p^md) \lra (\Omega^{\b}_{B},p^{m-1}d)$ 
of differential graded algebras 
such that the composite 
$$\Omega_{B^{(1)}}^i \os{\ol{F}^*}{\lra} \Omega_{B}^i \os{p^i}{\lra} 
 \Omega_{B}^i$$ 
is the one defined by the homomorphism 
$\psi: B^{(1)} \lra B$ induced by $\psi_n$'s. 
(This $\psi$ is a lift of relative Frobenius morphism 
$B^{(1)}_1 \lra B_1$, and the map $\ol{F}^*$ for $\b=1$ is the same 
as the one defined in Section 1.) Furthermore, 
the morphisms 
$\epsilon_n, \epsilon^{(1)}_n \,(n \in \N)$ induce the 
homomorphisms 
\begin{align*}
& \epsilon_n: 
(\Omega^{\b}_{B_n}, p^{m-1}d) \lra (W_n\Omega_{B_1}^{\b}, p^{m-1}d), \\ 
& \epsilon := \varprojlim_n \epsilon_n: 
(\Omega^{\b}_B, p^{m-1}d) \lra (W\Omega_{B_1}^{\b}, p^{m-1}d), \\ 
& \epsilon^{(1)}_n: 
(\Omega^{\b}_{B^{(1)}_n}, p^{m}d) \lra 
(W_n\Omega_{B^{(1)}_1}^{\b}, p^{m}d), \\ 
& \epsilon^{(1)} := \varprojlim_n \epsilon_n^{(1)}: 
(\Omega^{\b}_{B^{(1)}}, p^{m}d) \lra 
(W\Omega_{B^{(1)}_1}^{\b}, p^{m}d),
\end{align*}
and the well-known formula $Fd[x] = [x]^{p-1}d[x]$ of de Rham-Witt 
complexes implies the commutativity of the following diagram: 
\begin{equation}\label{starmaru}
\begin{CD}
(W\Omega_{B^{(1)}_1}^{\b}, p^{m}d) 
@>{F}>>
(W\Omega_{B_1}^{\b}, p^{m-1}d) \\ 
@A{\epsilon^{(1)}}AA @A{\epsilon}AA \\ 
(\Omega^{\b}_{B^{(1)}}, p^{m}d)
@>{\ol{F}^*}>> (\Omega^{\b}_B, p^{m-1}d). 
\end{CD}
\end{equation}
To prove the theorem, it suffices to prove that the homomorphisms 
\begin{align}
& H^0(F): H^0(W\Omega_{B^{(1)}_1}^{\b}, p^{m}d) \lra 
H^0(W\Omega_{B_1}^{\b}, p^{m-1}d), \label{u1} \\ 
& 
H^i(F) \otimes \Q: H^i(W\Omega_{B^{(1)}_1}^{\b}, p^{m}d) \otimes \Q \lra 
H^i(W\Omega_{B_1}^{\b}, p^{m-1}d) \otimes \Q \,\,\,\, (i \in \N) \label{u2}
\end{align}
induced by $F$ are isomorphisms. We prove this by using the commutative 
diagram \eqref{starmaru}. First, by Theorem \ref{25}, the 
homomorphisms 
\begin{align}
& H^0(\ol{F}^*): H^0(\Omega_{B^{(1)}}^{\b}, p^{m}d) \lra 
H^0(\Omega_{B}^{\b}, p^{m-1}d), \label{v1} \\ 
& 
H^i(\ol{F}^*) \otimes \Q: H^i(\Omega_{B^{(1)}}^{\b}, p^{m}d) \otimes \Q \lra 
H^i(\Omega_{B}^{\b}, p^{m-1}d) \otimes \Q \,\,\,\, (i \in \N) \label{v2}
\end{align}
induced by $\ol{F}^*$ are isomorphisms. ($\Hyp(\infty,\infty,\infty,\infty)$ 
is satisfied because of the left cocartesian diagram in \eqref{Flcomp2}.) 
Next we prove that the homomorphisms 
\begin{align}
& H^0(\epsilon): H^0(\Omega_{B}^{\b}, p^{m-1}d) \lra 
H^0(W\Omega_{B_1}^{\b}, p^{m-1}d), \label{w1} \\ 
& 
H^i(\epsilon) \otimes \Q: 
H^i(\Omega_{B}^{\b}, p^{m-1}d) \otimes \Q \lra 
H^i(W\Omega_{B_1}^{\b}, p^{m-1}d) \otimes \Q  \,\,\,\, (i \in \N) \label{w2}
\end{align}
induced by $\epsilon$ are isomorphisms. To prove this, 
we follow the argument in \cite[Theorem 3.5]{lz}: First, let us note the 
commutative diagram proven there 
\begin{equation}\label{lzdiag}
\begin{CD}
(\Omega^{\b}_{B_n},d) @<=<< 
(B_{2n} \otimes_{A_{2n},\phi^n} \Omega^{\b}_{A_n}, \id \otimes d) 
\\ 
@V{\epsilon_n}VV @V{\id \otimes \delta_n}VV \\ 
(W_n\Omega^{\b}_{B_1},d) @<=<< 
(B_{2n} \otimes_{A_{2n},\phi^n} W_n\Omega_{A_1}, \id \otimes d), 
\end{CD}
\end{equation}
where we denoted the maps induced from $\epsilon_n, \delta_n$ by the 
same latters abusively, and 
$\phi^n: A_{2n} \lra A_n$ is the map $\phi_{n+1} \circ \cdots 
\circ \phi_{2n}$. 
Also, in the proof of \cite[Theorem 3.5]{lz}, 
it is shown that we have the decomposition 
$(W_n\Omega^{\b}_{A_1},d) = (C_{\rm int}^{\b},d) \oplus 
(C_{\rm frac}^{\b},d)$ 
such that $\delta_n: (\Omega^{\b}_{A_n}, d) \lra (W_n\Omega_{A_1}, d)$ 
induces the isomorphism of complexes 
$(\Omega^{\b}_{A_n}, d) \os{=}{\lra} (C_{\rm int}^{\b},d)$ and that 
the complex $(C_{\rm frac}^{\b},d)$ is acyclic. Hence 
$\delta_n: (\Omega^{\b}_{A_n}, d) \lra (W_n\Omega_{A_1}, d)$ is 
a quasi-isomorphism and so is 
$\epsilon_n: (\Omega^{\b}_{B_n},d) \lra (W_n\Omega^{\b}_{B_1},d)$ 
by \eqref{lzdiag}. Hence $\epsilon := \varprojlim_n \epsilon_n: 
(\Omega^{\b}_{B},d) \lra (W\Omega^{\b}_{B_1},d)$ is also a 
quasi-isomorphism. Since $\Omega^{i}_B, W\Omega_{B_1}^{i}$'s ($i=0,1$) 
are $p$-torsion free, we have 
$H^0(\Omega_B^{\b}, p^{m-1}d) \allowbreak = H^0(\Omega_B^{\b},d), 
H^0(W\Omega_{B_1}^{\b},p^{m-1}d) = H^0(W\Omega_{B_1}^{\b}, d)$ and so 
the map \eqref{w1} is an isomorphism. Also, 
using the commutative diagrams 
\begin{equation*}
\begin{CD}
\Omega^{i-1}_B @>{=}>> \Omega^{i-1}_B @>{p^{2(m-1)}}>> \Omega^{i-1}_B \\ 
@VdVV @V{p^{m-1}d}VV @VdVV \\ 
\Omega^{i}_B @>{p^{m-1}}>> \Omega^{i}_B @>{p^{m-1}}>> \Omega^{i}_B \\ 
@VdVV @V{p^{m-1}d}VV @VdVV \\ 
\Omega^{i+1}_B @>{p^{2(m-1)}}>> \Omega^{i+1}_B @>{=}>> \Omega^{i+1}_B, 
\end{CD}
\quad \quad \quad 
\begin{CD}
W\Omega^{i-1}_{B_1} @>{=}>> W\Omega^{i-1}_{B_1} @>{p^{2(m-1)}}>> 
W\Omega^{i-1}_{B_1} \\ 
@VdVV @V{p^{m-1}d}VV @VdVV \\ 
W\Omega^{i}_{B_1} @>{p^{m-1}}>> W\Omega^{i}_{B_1} @>{p^{m-1}}>> 
W\Omega^{i}_{B_1} \\ 
@VdVV @V{p^{m-1}d}VV @VdVV \\ 
W\Omega^{i+1}_{B_1} @>{p^{2(m-1)}}>> W\Omega^{i+1}_{B_1} @>{=}>> 
W\Omega^{i+1}_{B_1}, 
\end{CD}
\end{equation*}
\begin{equation*}
\begin{CD}
\Omega^{i-1}_B @>{p^{2(m-1)}}>> \Omega^{i-1}_B @>{\epsilon}>>
 W\Omega^{i-1}_{B_1} @>{=}>> W\Omega^{i-1}_{B_1} \\ 
@V{p^{m-1}d}VV @VdVV @VdVV @V{p^{m-1}d}VV \\ 
\Omega^{i}_B @>{p^{m-1}}>> \Omega^{i}_B @>{\epsilon}>>
 W\Omega^{i}_{B_1} @>{p^{m-1}}>> W\Omega^{i}_{B_1} \\ 
@V{p^{m-1}d}VV @VdVV @VdVV @V{p^{m-1}d}VV \\ 
\Omega^{i+1}_B @>{=}>> \Omega^{i+1}_B @>{\epsilon}>>
 W\Omega^{i+1}_{B_1} @>{p^{2(m-1)}}>> W\Omega^{i+1}_{B_1}, 
\end{CD}
\end{equation*}
we can show that the kernal and the cokernel of the homomorphism 
$$ p^{2(m-1)}H^i(\epsilon):  H^i(\Omega_{B}^{\b}, p^{m-1}d) \lra 
H^i(W\Omega_{B_1}^{\b}, p^{m-1}d) $$ 
are killed by $p^{4(m-1)}$. Hence \eqref{w2} is also an isomorphism. \par 
In the same way as above, we can also prove that the homomorphisms 
\begin{align}
& H^0(\epsilon^{(1)}): H^0(\Omega_{B^{(1)}}^{\b}, p^{m}d) \lra 
H^0(W\Omega_{B^{(1)}_1}^{\b}, p^{m}d), \label{x1} \\ 
& 
H^i(\epsilon^{(1)}) \otimes \Q: 
H^i(\Omega_{B^{(1)}}^{\b}, p^{m}d) \otimes \Q \lra 
H^i(W\Omega_{B_1^{(1)}}^{\b}, p^{m}d) \otimes \Q  \,\,\,\, (i \in \N) 
\label{x2}
\end{align}
induced by $\epsilon^{(1)}$ are isomorphisms. Then, by 
\eqref{v1}, \eqref{v2}, \eqref{w1}, \eqref{w2}, 
\eqref{x1}, \eqref{x2} and the diagram \eqref{starmaru}, we can 
conclude that the homomorphisms \eqref{u1}, \eqref{u2} are 
isomorphisms, as desired. 
\end{proof}

\end{document}